\documentclass[12pt]{article}

\usepackage{amssymb}
\usepackage{amsmath}
\usepackage{color}
\usepackage{graphics}
\usepackage{epsfig}
\usepackage{hyperref}

\newtheorem{claim}{\bf \t}[part]

\newtheorem{Lemma}{Lemma}[part]

\newtheorem{Remark}{Remark}[part]
\newtheorem{Theorem}{Theorem}[part]

\numberwithin{Assumption}{section} \numberwithin{Corollary}{section}
\numberwithin{Definition}{section} \numberwithin{equation}{section}
\numberwithin{Example}{section} \numberwithin{Lemma}{section}
\numberwithin{Proposition}{section} \numberwithin{Remark}{section}
\numberwithin{Theorem}{section}

\def\v{\varepsilon}

\def\x{\xi}
\def\t{\theta}

\def\a{\alpha}
\def\b{\beta}
\def\g{\gamma}

\def\l{\lambda}
\def\f{\frac}

\def\r{\rho}

\def\s{\sigma}

\def\o{\omega}
\def\di{\displaystyle}
\def\i{\infty}

\def\text#1{{\rm #1}}
\begin{document}
\date{}
\title{\Large \bf Global well-posedness of the Cauchy problem of two-dimensional compressible Navier-Stokes equations in weighted spaces}
\author{\small \textbf{Quansen Jiu},$^{1,3}$\thanks{The research is partially
supported by National Natural Sciences Foundation of China (No.
11171229) and Project of Beijing Education Committee. E-mail:
jiuqs@mail.cnu.edu.cn}\quad
  \textbf{Yi Wang}$^{2,3}$\thanks{The research is partially supported by National Natural Sciences Foundation of China (No.
11171326), and by the National Center for Mathematics and Interdisciplinary Sciences, CAS. E-mail: wangyi@amss.ac.cn.}\quad and \textbf{Zhouping
Xin}$^{3}$\thanks{The research is partially supported by Zheng Ge Ru
Funds, Hong Kong RGC Earmarked Research Grant CUHK4042/08P and
CUHK4041/11P, and a Focus Area Grant at The Chinese University of
Hong Kong. Email: zpxin@ims.cuhk.edu.hk}} \maketitle \small $^1$
School of Mathematical Sciences, Capital Normal University, Beijing
100048, P. R. China

\small $^2$ Institute of Applied Mathematics, AMSS, CAS, Beijing 100190, P. R. China

\small $^3$The Institute of Mathematical Sciences, Chinese
University of HongKong, HongKong\\

 {\bf Abstract:} In this paper, we study the global well-posedness of classical solution to 2D Cauchy problem of
 the compressible Navier-Stokes equations with large initial data and vacuum. It is proved that if the shear viscosity
 $\mu$ is a positive constant and the bulk viscosity $\l$ is the power function of the density,
 that is, $\l(\r)=\r^\b$ with $\b>3$, then the 2D Cauchy problem of the compressible Navier-Stokes equations on the
 whole space $\mathbb{R}^2$ admit a unique global classical solution $(\r,u)$ which may contain
  vacuums in an open set of $\mathbb{R}^2$.  Note that the initial data can be arbitrarily large to contain vacuum
  states. Various weighted estimates of the density and velocity are
  obtained in this paper and these self-contained estimates reflect the
  fact that the weighted density and weighted velocity propagate along with the flow.



\section{Introduction and main results} \setcounter{equation}{0}
\setcounter{Assumption}{0} \setcounter{Theorem}{0}
\setcounter{Proposition}{0} \setcounter{Corollary}{0}
\setcounter{Lemma}{0} In this paper, we consider the following
compressible and isentropic Navier-Stokes equations with
density-dependent viscosities
\begin{eqnarray}\label{CNS}
\left\{ \begin{array}{ll}
\partial_t\rho+{\rm div}(\rho u)=0,  \\
\partial_t(\rho u) + {\rm div} (\rho u\otimes u) + \nabla P(\r)=\mu\Delta u+\nabla((\mu+\l(\r)){\rm div}u), &\quad x\in\mathbb{R}^2, t>0,
\end{array}
\right.
\end{eqnarray}
where $\rho (t, x)\geq 0 $, $u(t, x)=(u_1,u_2)(t,x) $ represent the
density and the velocity of the fluid, respectively. And
$x=(x_1,x_2)\in\mathbb{R}^2$ and $t\in[0,T]$ for any fixed $T>0$. We
denote the right hand side of $\eqref{CNS}_2$  by
$$
\mathcal{L}_\r u=\mu\Delta u+\nabla((\mu+\l(\r)){\rm div}u).
$$
Here, it is assumed that
\begin{equation}\label{1.2}
\mu={\rm const.}>0,\qquad \l(\r)=\r^\beta, \quad \beta>3,
\end{equation}
such that the operator $\mathcal{L}_\r$ is strictly elliptic.

For simplicity, we assume that the pressure function is given by
\begin{eqnarray}\label{ga}
P(\r)=A \r^{\gamma },
\end{eqnarray}
where $\gamma >1$ denotes the adiabatic exponent and $A>0$ is the
constant. Without loss of generality, $A$ is normalized to be $1$.
The initial values are imposed as
\begin{equation}\label{initial-v}
(\r,u)(t=0,x)=(\r_0,u_0)(x).
\end{equation}
The system  \eqref{CNS}-\eqref{1.2} was first proposed and studied
by Vaigant-Kazhikhov in \cite{Kazhikhov}  in which the global
well-posedness of the classical solution to  \eqref{CNS}-\eqref{ga}
with general data satisfying periodic boundary conditions was
obtained provided that the initial density is uniformly away from
vacuum. To authors's knowledge,  this is the first result of the
global well-posedness to the multi-dimensional compressible
Navier-Stokes equations with large initial data in the absence of
vacuum.
 Then Perepelitsa \cite{P} studied the global existence and large time behavior of weak solution
 to \eqref{CNS}-\eqref{1.2} with 2D periodic boundary conditions.
Recently, Jiu-Wang-Xin \cite{JWX2} improved the result in
\cite{Kazhikhov} and obtained the global well-posedness of the
classical solution to the periodic problem with general initial data
permitting vacuum. Later on, based on \cite{Kazhikhov}, \cite{P} and
\cite{JWX2}, Huang-Li relaxed the index $\beta$ to be $\beta>\frac
43$ and studied the large time behavior of the solutions in
\cite{HL2012}. However, all these results are concerned with the 2D
periodic problems. In the present paper, we are interested in the
global existence and uniqueness of classical solution to 2D Cauchy
problem with large  data and vacuum.

There are extensive studies on global well-posedness of the
compressible Navier-Stokes equations when both the shear and bulk
viscosities are positive constants.  In particular, the
one-dimensional theory is rather satisfactory, see
\cite{Hoff-Smoller, lius, Ka, KS} and the references therein. In
multi-dimensional case,
 the local well-posedness theory of classical solutions  was
 established   in the absence of
 vacuum (see \cite{Nash},  \cite{Itaya} and  \cite{Tani}) and the global well-posedness theory of classical solutions  was
 obtained  for initial data
  close to a non-vacuum steady state (see \cite{MN},  \cite{Hoff1},   \cite{Dachin},  \cite{CMZ} and references therein).
    For  the large  initial data which may contain vacuums,  the global existence of weak solutions was
    obtained when $\gamma>\frac N2, (N=2,3)$ in general case and $\gamma>1$ if assuming space symmetry (see  \cite{Lions}, \cite{F1}, \cite{JZ}).
 However, the
 uniqueness of such  weak solutions remain completely open in general. By the weak-strong uniqueness of
   \cite{G2011}, this is equivalent to the problem of global (in time) well-posedness of strong solution
  in the
  presence of vacuum.   It should be noted that if the solutions contain possible vacuums, the regularity
  and uniqueness  become difficult and subtle issues.
  In 1998,  Xin showed
\cite{Xin}  that if the initial density has compact support, any smooth solution in $C^1([0,T];H^s(R^N))$ with $s\ge [N/2]+2$ to the
  Cauchy problem of the CNS without heat conduction blows up in finite time for any $N\ge 1$. Then Rozanova \cite{R} generalized the results in \cite{Xin} to the case the data with highly decreasing at infinity.
  Very recently, Xin-Yan \cite{Xin-Yan} improves the blow-up results in
  \cite{Xin} by removing the assumptions that the initial density has
  compact support and the smooth solution has finite energy.
On the other hand, the short time well-posedness
 of either strong or classical solutions containing vacuum was studied recently by Cho-Kim \cite{CK1}
 and Luo\cite{Luo} in 3D and 2D case,  respectively.  A
 natural compatibility condition was imposed in \cite{CK1} to guarantee the local well-posedness  of the classical
 solution for the isentropic CNS
 with general nonnegative initial density. More recently, Huang-Li-Xin \cite{hlx4} proved the global well-posedness
 of classical solutions with small energy but large oscillations which can contain vacuums
  to 3D isentropic compressible Navier-Stokes equations.

The case that the viscosity coefficients depend on the density has
received a lot attention recently, see \cite{BD1, BDL, BD3, Dachin, DWZ,
GJX, J, JXZ, JZ, JWX, JWX1, JX, LLX, LXY, MV, MV2, SZ, YYZ, YZ2, ZF}
and the references therein. When deriving by Chapman-Enskog
expansions from the Boltzmann equation, the viscosity of the
compressible Navier-Stokes equations  depends on the temperature and
thus on the density for isentropic flows (see \cite{LXY}). On the
other hand, in the geophysical flow, the viscous Saint-Venant system
for the shallow water corresponds exactly to a kind of compressible
Navier-Stokes equations with density-dependent viscosities (see
\cite{GP}). Similar to the case of constant viscosities, the
well-posedness theory to the one-dimensional problem with viscosity
coefficients depending on the density has been well-understood.
However, the progress is very limited for multi-dimensional
problems. Even the short time well-posedness of strong or classical solutions
has not been established  in the presence of vacuum. Also, the
global existence of  weak solutions to the compressible
Navier-Stokes equations with density-dependent viscosities  remains
open, except assuming some space symmetry  \cite{GJX}.  One can
refer to \cite{BD3}, \cite{GLX}, \cite{MV} and references therein
for recent developments along this line.

In this paper, we are concerned with the global well-posedness of
the classical solution to the 2D Cauchy problem
\eqref{CNS}-\eqref{initial-v} with  general data permitting  vacuum.
Compared with \cite{Kazhikhov} and \cite{JWX2} for 2D periodic
problems, some new difficulties occur. First, the Poincare-type
inequality fails for the 2D Cauchy problem. In particular, by the
elementary energy estimates in Lemma \ref{lemma-ee} we have $\sqrt\r
u\in L^\i(0,T;L^2(\mathbf{R}^2))$ and $\nabla u\in
L^2((0,T)\times\mathbf{R}^2)$,
 which give no information on the integrability of the velocity $u$ in the whole space. While the $L^p$-integrability $(1\le p< \infty)$
of the velocity $u$ plays a very important role in the arbitrarily
$L^p$-integrability $(1\le p< \infty)$ estimates of the density $\r$
in Lemma \ref{lemma-rho}. One way to encounter this difficulty is to
get the weighted estimates of the velocity like $|x|^{\f\a2}\nabla
u\in L^2((0,T)\times\mathbf{R}^2)$ with suitable $\alpha>0$, which
will lead to $u\in L^2(0,T; L^{\f4\a}(\mathbf{R}^2))$. The weighted
estimates of the velocity $|x|^{\f\a2}\nabla u\in
L^2((0,T)\times\mathbf{R}^2)$
 are strongly coupled
with the higher integrability estimates of the density function, see
\eqref{We-l1} in Lemma \ref{lemma-wee}. Therefore, the delicate
combination of the new weighted estimates to the velocity and the
techniques in \cite{Kazhikhov} and \cite{JWX2} for 2D periodic
problems yields the expected $L^p$-integrability $(1\le p< \infty)$
estimates of the density $\r$ in Lemma \ref{lemma-rho}. With Lemma
\ref{lemma-rho} in hand, one can get the higher order estimates as
in \cite{Kazhikhov} and \cite{JWX2} to get the upper bound of the
density. For this, the weighted estimates for the density (Lemma
\ref{density-wee}) and weighted estimates of $\nabla \dot u$ (see
Lemma \ref{lemma-u-sec}) will be established. In particular, it is
highly nontrivial to get the weighted estimates $\||x|^{\f\a2}\nabla
\dot u\|_{L^2_{x,t}}$ or $\||x|^{\f\a2}\sqrt\r \dot
u\|_{L^\i_tL^2_x}$ (see \eqref{HL-e1}) where  we will make full use
of the sharp Cafferelli-Kohn-Nirenberg inequality \cite{CKN} and
Gagliardo-Nirenberg inequaltiy \cite{CW} for the best constants and
the weighted estimates for the singular integral operators
\cite{Stein}. Furthermore, the higher order estimates in Lemma
\ref{lemma4.7}  also involve the weighted estimates and are also
crucial to get our results. It should be noted that in \cite{Luo}
Luo studied the Cauchy problem of the 2D compressible Navier-Stokes
equations with constant viscosities by  using the homogeneous
Sobolev spaces and weighted estimates and obtained the local
existence and uniqueness of the classical solution for large data
with vacuum. The main reason that  the global existence and
uniqueness of the classical solution in this paper is that the bulk
viscosity $\lambda(\r)=\r^\b$ with $\b>3$ will provide higher $L^p
(1<p<\infty)$ estimates of the density and based on this observation
we can furthermore obtain the upper bound of the density, and then
get our main results.

It is also interesting to obtain various weighted estimates of the
density and velocity itself in $L^P (1<p<\infty)$
spaces. These self-contained
estimates reflect the
  fact that the weighted density and weighted velocity propagate along with the flow. Moreover,
the weighted estimates will provide an appropriate approach to deal
with the two-dimensional Cauchy problem of other fluid models having
similar structure. As an example, it is possible that the methods
here can be applied to 2D Cauchy problem of MHD systems as in
\cite{BG}. Very recently, we just learned that Huang-Li \cite{hl2}
studied the Cauchy problem \eqref{CNS}-\eqref{initial-v}
independently and obtained the global well-posedness of strong and
classical solution in quite different weighted spaces and by using
different approaches, in which the index $\beta>\frac 43$.

The main results of the present paper can be stated in the
following.

\begin{Theorem}\label{theorem2}
Suppose that the initial values $(\r_0,u_0)(x)$ satisfy
\begin{equation}\label{in-d1}
\begin{array}{ll}
\di 0\leq(\r_0(x), P(\r_0)(x))\in
W^{2,q}(\mathbb{R}^2)\times W^{2,q}(\mathbb{R}^2),\quad u_0(x)\in D^1\cap D^2(\mathbb{R}^2),\\
\di \r_0(1+|x|^{\a_1})\in L^1(\mathbb{R}^2), \quad
\sqrt{\r_0}u_0(1+|x|^{\f\a2})\in L^2(\mathbb{R}^2),\quad \nabla
u_0|x|^{\f\a2} \in L^2(\mathbb{R}^2),
\end{array}
\end{equation}
for some $q>2$ and the weights $0<\a<2\sqrt{\sqrt2-1}$, $\a<\a_1$,
and the compatibility condition
\begin{equation}\label{cc}
\mathcal{L}_{\r_0}u_0-\nabla P(\r_0)=\sqrt\r_0  g(x)
\end{equation}
with some $g$ satisfying $g(1+|x|^{\f\a2})\in L^2(\mathbb{R}^2)$. If
one of the  following restrictions holds:
\begin{eqnarray}
&& 1) \quad 1<\a<2\sqrt{\sqrt2-1}, \ \beta>3,\
\gamma>1, \label{June9-1}\\[3mm]
 &&2)\quad  0<\alpha\le 1, \ \beta>3,\ 1<\gamma\le 2\b, \label{June9-2}
\end{eqnarray}
 then there exists a unique global classical solution
$(\r,u)(t,x)$ to the Cauchy problem \eqref{CNS}-\eqref{initial-v}
with
\begin{equation}\label{Jan-16-1}
\begin{array}{ll}
\di 0\leq \r\leq C,\quad (\r,P(\r))(t,x)\in C([0,T]; W^{2,q}(\mathbb{R}^2)),\quad \r(1+|x|^{\a_1})\in C([0,T]; L^1(\mathbb{R}^2)),\\
\di \sqrt\r u  (1+|x|^{\f\a2}), \sqrt\r \dot u(1+|x|^{\f\a2}), \nabla u|x|^{\f\a2}\in C([0,T]; L^2(\mathbb{R}^2)),\\
u\in C([0,T];L^{\f4\a}\cap D^2(\mathbb{R}^2))\cap
L^2(0,T;L^{\f4\a}\cap D^3(\mathbb{R}^2)),~~ \sqrt t u\in L^\i(0,T;
D^{3}(\mathbb{R}^2)),\\t u\in L^\i(0,T; D^{3,q}(\mathbb{R}^2)),~~
u_t\in L^2(0,T;L^{\f4\a}(\mathbb{R}^2)\cap D^1(\mathbb{R}^2))\\
\di \sqrt tu_t\in L^2(0,T; D^2(\mathbb{R}^2))\cap
L^\i(0,T;L^{\f8\a}\cap D^1(\mathbb{R}^2)),~~ t u_t\in L^\i(0,T;
D^2(\mathbb{R}^2)),\\ \sqrt t\sqrt\r u_{tt}\in
L^2(0,T;L^2(\mathbb{R}^2)),~~
 t\sqrt \r u_{tt}\in
L^\i(0,T;L^2(\mathbb{R}^2)),~~t\nabla u_{tt}\in
L^2(0,T;L^2(\mathbb{R}^2)),
\end{array}
\end{equation}
where $\dot u$ is the material derivative of $u$ defined in
\eqref{dot-u}.
\end{Theorem}

\begin{Remark}
From the regularity of the solution $(\r,u)(t,x)$, it can be shown
that $(\r,u)$ is a classical solution of the system \eqref{CNS} in
$[0,T]\times\mathbb{R}^2$ (see the details in Section 5).
\end{Remark}

\begin{Remark}
If the initial data contains vacuum, then  the compatibility
condition \eqref{cc} is necessary for the existence of the classical solution,
just as the case of constant viscosity coefficients in \cite{CK1}.
\end{Remark}

If the initial values are more regular, based on Theorem
\ref{theorem2}, we can prove

\begin{Theorem}\label{theorem}
Under assumptions of \eqref{in-d1}-\eqref{June9-2}, assume further
that
\begin{equation}\label{in-d}
0\leq(\r_0(x), P(\r_0)(x))\in H^3(\mathbb{R}^2)\times
H^{3}(\mathbb{R}^2),\quad u_0(x)\in D^1\cap D^3(\mathbb{R}^2)
\end{equation}
and the compatibility condition \eqref{cc}, then there exists a
unique global classical solution $(\r,u)(t,x)$ to the Cauchy problem
\eqref{CNS}-\eqref{initial-v} satisfying all the properties listed
in \eqref{Jan-16-1} in Theorem \ref{theorem2} with any $2<q<\i$.
Furthermore, it holds that
\begin{equation}\label{h-regu}
\begin{array}{ll}
\di u\in L^2(0,T;D^4(\mathbb{R}^2)),~~(\r,P(\r))\in C([0,T]; H^3(\mathbb{R}^2)),\\
\di \r u\in C([0,T];D^1\cap D^3(\mathbb{R}^2)),~~\sqrt\r \nabla^3
u\in C([0,T]; L^2(\mathbb{R}^2)).
\end{array}
\end{equation}
\end{Theorem}
\begin{Remark}
In fact, the conditions on the initial velocity $u_0$ in
\eqref{in-d} can be weakened to $u_0\in D^1\cap D^2(\mathbb{R}^2)$
and $\sqrt{\r_0}\nabla^3 u_0\in L^2(\mathbb{R}^2)$ to get
\eqref{h-regu}.
\end{Remark}

\begin{Remark}
In Theorem \ref{theorem}, it is not clear whether or not $u\in
C([0,T]; D^3(\mathbb{R}^2))$ even though one has $\r u\in C([0,T];
D^3(\mathbb{R}^2))$.
\end{Remark}

The plan of the paper is as follows. In Section 2, we present the reformulations of the system, some
elementary facts and inequalities. In
Sections 3-4, we derive a priori estimates which are needed to
extend the local solution to global one. Finally, in Section 5, we
prove our main results.

 \vskip 2mm
\noindent\emph{Notations.}Throughout this paper, positive generic
constants are denoted by $c$ and $C$, which are independent of
$m$ and $t\in[0,T]$, without confusion, and $C(\cdot)$ stands for
some generic constant(s) depending only on the quantity listed in
the parenthesis. For functional spaces, $L^{p}(\mathbb{R}^2), 1\leq
p\leq \infty$, denote the usual Lebesgue spaces on $\mathbb{R}^2$
and $\|\cdot\|_p$ denotes its $L^p$ norm. $W^{k,p}(\mathbb{R}^2)$
denotes the standard $k^{th}$ order Sobolev space and
$H^{k}(\mathbb{R}^2):=W^{k,2}(\mathbb{R}^2)$. For $1<p<\i$, the
homogenous Sobolev space $D^{k,p}(\mathbb{R}^2)$ is defined by
$D^{k,p}(\mathbb{R}^2)=\{u\in
L^1_{loc}(\mathbb{R}^2)|\|\nabla^ku\|_p<+\i\}$ with
$\|u\|_{D^{k,p}}:=\|\nabla^k u\|_p$ and
$D^k(\mathbb{R}^2):=D^{k,2}(\mathbb{R}^2)$.


\section{Preliminaries}

As in \cite{Kazhikhov}, we introduce the following variables. First denote the effective viscous flux by
\begin{equation}\label{flux}
F=(2\mu+\l(\r)){\rm div}u-P(\r),
\end{equation}
and the vorticity by
$$
\o=\partial_{x_1}u_2-\partial_{x_2}u_1.
$$
Also, we define that
$$
H=\frac{1}{\r}(\mu\o_{x_1}+F_{x_2}),\qquad\qquad
L=\frac{1}{\r}(-\mu\o_{x_2}+F_{x_1}).
$$
Then the momentum equation $\eqref{CNS}_2$ can be rewritten as
\begin{equation}\label{ns1}
\left\{
\begin{array}{ll}
\dot{u}_1=u_{1t}+u\cdot \nabla u_1=\frac{1}{\r}(-\mu\o_{x_2}+F_{x_1})=L,\\
\dot{u}_2=u_{2t}+u\cdot \nabla u_2=\frac{1}{\r}(\mu\o_{x_1}+F_{x_2})=H,
\end{array}
\right.
\end{equation}
that is,
\begin{equation}\label{dot-u}
\dot{u}=(\dot{u}_1,\dot{u}_2)^t=(L,H)^t.
\end{equation}
Then the effective viscous flux $F$ and the vorticity $\o$ solve the following system:
\begin{equation}
\left\{
\begin{array}{ll}
\o_{t}+u\cdot \nabla \o+\o{\rm div}u=H_{x_1}-L_{x_2},\\
(\f{F+P(\r)}{2\mu+\l(\r)})_{t}+u\cdot \nabla
(\f{F+P(\r)}{2\mu+\l(\r)})+(u_{1x_1})^2+2u_{1x_2}u_{2x_1}+(u_{2x_2})^2=H_{x_2}+L_{x_1}.
\end{array}
\right.
\end{equation}
Due to the continuity equation $\eqref{CNS}_1$, it holds that
\begin{equation}\label{F-omega}
\left\{
\begin{array}{ll}
\o_{t}+u\cdot \nabla \o+\o{\rm div}u=H_{x_1}-L_{x_2},\\
F_{t}+u\cdot \nabla
F-\r(2\mu+\l(\r))[F(\f{1}{2\mu+\l(\r)})^\prime+(\f{P(\r)}{2\mu+\l(\r)})^\prime]{\rm
div}u\\
\qquad+(2\mu+\l(\r))[(u_{1x_1})^2+2u_{1x_2}u_{2x_1}+(u_{2x_2})^2]=(2\mu+\l(\r))(H_{x_2}+L_{x_1}).
\end{array}
\right.
\end{equation}
Furthermore,  the system for $(H,L)$ can be derived as
\begin{equation}\label{H-L}
\left\{
\begin{array}{ll}
\r H_{t}+\r u\cdot \nabla H-\r H{\rm div}u+u_{x_2}\cdot\nabla F+\mu
u_{x_1}\cdot\nabla\o+\mu(\o{\rm
div}u)_{x_1}\\
\qquad
-\big\{\r(2\mu+\l(\r))[F(\f{1}{2\mu+\l(\r)})^\prime+(\f{P(\r)}{2\mu+\l(\r)})^\prime]{\rm
div}u\big\}_{x_2}\\
\qquad+\big\{(2\mu+\l(\r))[(u_{1x_1})^2+2u_{1x_2}u_{2x_1}+(u_{2x_2})^2]\big\}_{x_2}\\
\qquad=[(2\mu+\l(\r))(H_{x_2}+L_{x_1})]_{x_2}+\mu(H_{x_1}-L_{x_2})_{x_1},\\
\r L_{t}+\r u\cdot \nabla L-\r L{\rm div}u+u_{x_1}\cdot\nabla F-\mu
u_{x_2}\cdot\nabla\o-\mu(\o{\rm
div}u)_{x_2}\\
\qquad
-\big\{\r(2\mu+\l(\r))[F(\f{1}{2\mu+\l(\r)})^\prime+(\f{P(\r)}{2\mu+\l(\r)})^\prime]{\rm
div}u\big\}_{x_1}\\
\qquad+\big\{(2\mu+\l(\r))[(u_{1x_1})^2+2u_{1x_2}u_{2x_1}+(u_{2x_2})^2]\big\}_{x_1}\\
\qquad=[(2\mu+\l(\r))(H_{x_2}+L_{x_1})]_{x_1}-\mu(H_{x_1}-L_{x_2})_{x_2}.
\end{array}
\right.
\end{equation}
In the following, we will utilize the above systems in different
steps. Note that these systems are equivalent to each other for the
smooth solution to the original system \eqref{CNS}.
We first state the local existence and uniqueness of classical
solution when the initial data may contain vacuum.
\begin{Lemma}\label{lemma0}
Under assumptions of Theorem \ref{theorem2}, there exists a $T_*>0$
and a unique classical solution $(\rho,u)$  to the Cauchy problem
\eqref{CNS}-\eqref{initial-v} satisfying \eqref{Jan-16-1} with $T$
replaced by $T_*$.

\end{Lemma}

Lemma \ref{lemma0} can be proved in a similar way as in \cite{CK1}
and \cite{Luo}, by using the linearization method, Schauder fixed point
theorem and borrowing a priori estimates in Sections 3-4 of this
paper. We omit the details here.

 Several elementary Lemmas are needed later. The first
one is the various Gagliardo-Nirenberg inequalities.

\begin{Lemma}\label{lemma1}

\begin{itemize}
\item[(1)]
$\forall h\in W^{1,m}(\mathbb{R}^2)\cap L^r(\mathbb{R}^2)$, it holds that
\begin{equation}
\|h\|_q\leq C\|\nabla h\|_m^\t\|h\|_r^{1-\t},
\end{equation}
where $\t=(\f1r-\f1q)(\f1r-\f1m+\f12)^{-1}$, and if $m<2,$ then $q$
is between $r$ and $\f{2m}{2-m}$, that is, $q\in[r,\f{2m}{2-m}]$ if
$r<\f{2m}{2-m}$, $q\in[\f{2m}{2-m},r]$ if $r\geq\f{2m}{2-m},$ if
$m=2,$ then $q\in[r,+\i)$, if $m>2$, then $q\in[r,+\i].$
\item[(2)](Best constant for Gagliardo-Nirenberg inequality)

$\forall h\in \mathbb{D}^m(\mathbb{R}^2)\doteq\Big\{h\in L^{m+1}(\mathbb{R}^2)\Big|\nabla h\in L^2(\mathbb{R}^2),h\in L^{2m}(\mathbb{R}^2)\Big\}$ with $m>1$, it holds that
\begin{equation}
\|h\|_{2m}\leq A_m\|\nabla h\|_2^\t\|h\|_{m+1}^{1-\t},
\end{equation}
where $\t=\f12-\f{1}{2m}$ and
$$
A_m=\Big(\f{m+1}{2\pi}\Big)^{\f{\t}{2}}\Big(\f{2}{m+1}\Big)^{\f{1}{2m}}\leq C m^{\f14}
$$
with the positive constant $C$ independent of $m$, and $A_m$ is the optimal constant.

\item[(3)] $\forall h\in W^{1,m}(\mathbb{R}^2)$ with
$1\leq m<2,$ then
\begin{equation}
\|h\|_{\f{2m}{2-m}}\leq C(2-m)^{-\f12}\|\nabla h\|_m,
\end{equation}
where the positive constant $C$ is independent of $m.$

\item[(4)] $\forall h\in W^{1,\f{2m}{m+\eta}}(\mathbb{R}^2)$ with $m\geq2$ and
$0<\eta\leq1$, we have
\begin{equation}
\|h\|_{2m}\leq Cm^{\f12}\| h\|_{2(1-\v)}^s\|\nabla
h\|_{\f{2m}{m+\eta}}^{1-s},
\end{equation}
where $\v\in[0,\f12], s=\f{(1-\v)(1-\eta)}{m-\eta(1-\v)}$  and the
positive constant $C$ is independent of $m.$
\end{itemize}
\end{Lemma}

{\bf Proof:} The proof of (1) can be found in \cite{NS} while the proof of (2) can be found in \cite{Del-Pino}.
The proof of (3) can be found in \cite{GT} and the proof of (4) is a direct consequence of (2) and the interpolation inequality. $\hfill\Box$

The following Lemma is about the Caffarelli-Kokn-Nirenberg inequalities, which are crucial to the weighted estimates in 2D Cauchy problem.

\begin{Lemma}\label{lemma2}
\begin{itemize}
\item [(1)] $\forall h\in C^\i_0(\mathbb{R}^2)$, it holds that
\begin{equation}
\||x|^\kappa h\|_r\leq C\||x|^\a |\nabla h|\|^\t_p~\||x|^\b h\|^{1-\t}_q
\end{equation}
where $1\leq p,q<\i, 0<r<\i, 0\leq \t\leq 1, \f1p+\f{\a}{2}>0,\f1q+\f\b2>0,\f1r+\f\kappa2>0$ and satisfying
\begin{equation}
\f1r+\f\kappa2=\t(\f1p+\f{\a-1}{2})+(1-\t)(\f1q+\f\b2),
\end{equation}
and
$$
\kappa=\t\sigma+(1-\t)\b,
$$
with $0\leq \a-\sigma$ if $\t>0$ and $0\leq \a-\sigma\leq 1$ if $\t>0$ and $\f1p+\f{\a-1}{2}=\f1r+\f\kappa2.$

\item [(2)](Best constant for Caffarelli-Kohn-Nirenberg inequality)

$\forall h\in C^\i_0(\mathbb{R}^2)$, it holds that
\begin{equation}\label{b-CKN}
\||x|^b h\|_p\leq C_{a,b}\||x|^a\nabla h\|_2
\end{equation}
where $a>0, a-1\leq b\leq a$ and $p=\f{2}{a-b}$.  If $b=a-1$, then $p=2$ and the best constant in the inequality \eqref{b-CKN} is $$C_{a,b}=C_{a,a-1}=a.$$
\end{itemize}
\end{Lemma}
{\bf Proof:} The proof of (1) can be found in \cite{CKN} while the proof of (2) can be found in \cite{CW}.$\hfill\Box$

\begin{Lemma}\label{lemma3}
\begin{itemize}
\item [(1)] It holds that for $1<p<\i$ and $u\in C_0^\i(\mathbb{R}^2)$,
\begin{equation}
\|\nabla u\|_p\leq C(\|{\rm div} u\|_p+\|\o\|_p);
\end{equation}

\item [(2)] It holds that for $1<p<\i$, $-2<\a<2(p-1)$ and $u\in C_0^\i(\mathbb{R}^2)$,
\begin{equation}
\||x|^{\f\a p}|\nabla u|\|_p\leq C(\||x|^{\f\a p}{\rm div} u\|_p+\||x|^{\f\a p}\o\|_p).
\end{equation}
\end{itemize}
\end{Lemma}

\noindent{\bf Proof:} (1) Since
\begin{equation*}
\Delta u=\nabla ({\rm div} u)-\nabla\times\nabla\times u=\nabla ({\rm div} u)-\nabla\times \o,
\end{equation*}
where $\nabla\times$ denotes the  3-dimensional {\it curl} operator,
and
$$
\nabla\times\o=(\partial_{x_2}\o,-\partial_{x_1}\o,0)
$$
is regarded as the 2-dimensional vector
$(\partial_{x_2}\o,-\partial_{x_1}\o)^t$, then it holds that
\begin{equation}\label{singu}
\nabla u=\nabla \Delta^{-1}\nabla ({\rm div} u)-\nabla \Delta^{-1}\nabla\times \o:=\mathcal{T}_1({\rm div} u)+\mathcal{T}_2~\o,
\end{equation}
where $\mathcal{T}_1=\nabla \Delta^{-1}\nabla$ and $\mathcal{T}_2=-\nabla \Delta^{-1}\nabla\times$ both are the singular integral operators of the convolution type which are bounded in $L^p(\mathbb{R}^2)$. Thus Lemma \ref{lemma3} (1) is proved.

(2) If $-2<\a<2(p-1)$, then $|x|^\a$ is in the class $A_p$ (cf. p. 194 in \cite{Stein}), that is,
$$
\f{1}{|B|}\int_B|x|^\a dx\cdot\Big[\f{1}{|B|}\int_B(|x|^\a)^{-\f{p^\prime}{p}}dx\Big]^{\f{p}{p^\prime}}<\i,
$$
for all balls $B$ in $\mathbb{R}^2$, where $p^\prime$ is the dual to $p$, i.e., $\f{1}{p}+\f{1}{p^\prime}=1.$ Then by the Corollary in p. 205 of \cite{Stein}, there exist positive constants $C_1, C_2$ such that for $u\in C_0^\i(\mathbb{R}^2)$,
$$
\int_{\mathbb{R}^2}|\mathcal{T}_1({\rm div} u)|^p |x|^\a dx\leq C_1 \int_{\mathbb{R}^2}|{\rm div} u|^p |x|^\a dx,
$$
and
$$
\int_{\mathbb{R}^2}|\mathcal{T}_2~\o|^p |x|^\a dx\leq C_2 \int_{\mathbb{R}^2}|\o|^p |x|^\a dx.
$$
Therefore, it follows from \eqref{singu} that
$$
\begin{array}{ll}
\di \||x|^{\f\a p}|\nabla u|\|_p^p=\int_{\mathbb{R}^2}|\nabla u|^p |x|^\a dx\leq C_p\Big[\int_{\mathbb{R}^2}|\mathcal{T}_1({\rm div} u)|^p |x|^\a dx+\int_{\mathbb{R}^2}|\mathcal{T}_2~\o|^p |x|^\a dx\Big]\\
\di\qquad\quad\qquad~ \leq C_p\Big[\int_{\mathbb{R}^2}|{\rm div} u|^p |x|^\a dx+\int_{\mathbb{R}^2}|\o|^p |x|^\a dx\Big]=C_p\Big[\||x|^{\f\a p}{\rm div} u\|^p_p+\||x|^{\f\a p}\o\|_p^p\Big].
\end{array}
$$
Thus the proof of Lemma \ref{lemma3} (2) is completed.$\hfill\Box$

\section{A priori estimates (I)}
\setcounter{equation}{0}

In this section,  we obtain various a priori estimates and weighted
estimates on the classical solution $(\r,u)$ on the time interval $[0,T]$.
Denote
\begin{equation}
\label{M}
M=\|(\r_0,P(\r_0))\|_{W^{2,q}}+\|\r_0(1+|x|^{\a_1})\|_1+\|u_0\|_{D^1\cap D^2}+\||x|^{\f\a2}(\sqrt{\r_0} u_0,\nabla u_0)\|_2+\|g(1+|x|^{\f\a2})\|_2.
\end{equation}

\underline{Step 1. Elementary energy estimates:}
\begin{Lemma}\label{lemma-ee}
There exists a positive constant $C$ only depending on $M$, such that
\begin{equation*}
\sup_{t\in[0,T]}\big(\|\sqrt\r u\|^2_2+\|\r\|^\g_\g+\|\r\|_1\big)+\int_0^T\big(\|\nabla u\|_2^2+\|\o\|_2^2+\|(2\mu+\l(\r))^\f12{\rm div} u\|_2^2\big)dt\leq C(M).
\end{equation*}
\end{Lemma}
{\bf Proof:} Multiplying the equation $\eqref{CNS}_2$ by $u$ and the
continuity equation $\eqref{CNS}_1$ by $\f{\g}{\g-1}\r^{\g-1}$, then
summing the resulting equations, we have
\begin{equation}\label{e1}
\begin{array}{ll}
\di (\r\f{|u|^2}{2}+\f{\r^\g}{\g-1})_t+ {\rm div}(\r u\f{|u|^2}{2}+\f{\g\r^\g u}{\g-1})\\
\di ={\rm div}\big[\mu\nabla\f{|u|^2}{2}+(\mu+\l(\r))({\rm div} u)u\big]-\mu|\nabla u|^2-(\mu+\l(\r))({\rm div} u)^2.
\end{array}
\end{equation}
Integrating the above equality over $[0,t]\times\mathbb{R}^2$ with
respect to $t$ and $x$ yields that
\begin{equation}\label{Eee}
\begin{array}{ll}
  \di\int(\f12\r |u|^2+\f{1}{\g-1}\r^\g)dx+\int_0^t\int
\big[\mu|\nabla u|^2+(\mu+\l(\r))({\rm div} u)^2\big] dxdt\\
\di=\int(\f12\r_0 |u_0|^2+\f{1}{\g-1}(\r_0)^\g)dx\leq C.
\end{array}
\end{equation}
Note that
\begin{equation}\label{July-19}
\begin{array}{ll}
&\di \int \big[\mu|\nabla u|^2+(\mu+\l(\r))({\rm div} u)^2\big]dx=\int \big[\mu\o^2+(2\mu+\l(\r))({\rm div} u)^2\big]dx.
\end{array}
\end{equation}
Integrating the continuity equation $\eqref{CNS}_1$ with respect to
$t,x$ over $[0,t]\times\mathbb{R}^2$ yields that
$$
\int\r(t,x)dx=\int \r_0(x)dx.
$$
Thus the proof of Lemma \ref{lemma-ee} is completed. $\hfill\Box$

 \underline{Step 2. Weighted energy
estimates:}

The following weighted energy estimates is fundamental and crucial in our paper.
\begin{Lemma}\label{lemma-wee}
If one of the  following restrictions holds:
\begin{eqnarray}
&& 1) \quad 1<\a<2\sqrt{\sqrt2-1}, \ \beta>0,\
\gamma>1, \label{June21-1}\\[3mm]
 &&2)\quad  0<\alpha\le 1, \ \beta>\frac12,\ 1<\gamma\le 2\b, \label{June21-2}
\end{eqnarray}
then it holds that for sufficiently large $m>1$ and $\forall
t\in[0,T]$,
\begin{equation}\label{We-l1}
\begin{array}{ll}
\di \int_{\mathbb{R}^2}|x|^\a(\r|u|^2+\r^\g)(t,x)dx+\int_0^t\big[\||x|^{\f\a2}\nabla u\|_2^2(s)+\||x|^{\f\a2}{\rm div} u\|_2^2(s)+\||x|^{\f\a2}\sqrt{\l(\r)}{\rm div} u\|_2^2(s)\big]ds\\
\di \leq C_\a(M)\Big[1+\int_0^t(\|\r\|^\b_{2m\b+1}(s)+1)(\|\nabla
u\|_2^2(s)+1)ds\Big],
\end{array}
\end{equation}
where the positive constant $C_\a(M)$ depend on $\a$ and $M$ but is
independent of $m$.
\end{Lemma}
{\bf Proof:} Multiplying the equality \eqref{e1} by $|x|^\a$ yields
that
\begin{equation}\label{e2}
\begin{array}{ll}
\di \big[|x|^\a(\r\f{|u|^2}{2}+\f{\r^\g}{\g-1})\big]_t+\big[\mu|\nabla u|^2+(\mu+\l(\r))({\rm div} u)^2\big]|x|^\a\\
\di =-{\rm div}\big[|x|^\a(\r u\f{|u|^2}{2}+\f{\g\r^\g u}{\g-1})\big]+{\rm div}\big[\big(\mu\nabla\f{|u|^2}{2}+(\mu+\l(\r))({\rm div} u)u\big)|x|^\a\big]\\
\di \quad+(\r
\f{|u|^2}{2}+\f{\g\r^\g}{\g-1})u\cdot\nabla(|x|^\a)-\big[\mu\nabla\f{|u|^2}{2}+(\mu+\l(\r))({\rm
div} u)u\big]\cdot\nabla(|x|^\a).
\end{array}
\end{equation}
Integrating the above equation \eqref{e2} with respect to $x$ over
$\mathbb{R}^2$ yields that
\begin{equation}\label{e3}
\begin{array}{ll}
\di \f{d}{dt}\int |x|^\a(\r\f{|u|^2}{2}+\f{\r^\g}{\g-1})(t,x)dx+\mu\||x|^{\f\a2}\nabla u\|_2^2(t)+\mu\||x|^{\f\a2}{\rm div} u\|_2^2(t)+\||x|^{\f\a2}\sqrt{\l(\r)}{\rm div} u\|_2^2(t) \\
\di =\int(\r
\f{|u|^2}{2}+\f{\g\r^\g}{\g-1})u\cdot\nabla(|x|^\a)dx-\int
\big[\mu\nabla\f{|u|^2}{2}+(\mu+\l(\r))({\rm div}
u)u\big]\cdot\nabla(|x|^\a)dx.
\end{array}
\end{equation}
Note that the conservation terms in \eqref{e2} is vanished, which
can be proved rigourously by multiplying  a smooth cutting-off
function $\phi_R(x)=\phi(\frac{x}{R})$ on both sides of the equation
\eqref{e2}, where
$$
\phi(x)=\phi(|x|)=\left\{
\begin{array}{ll}
 & 1, \quad |x|\le 1,\\
& 0, \quad |x|\ge 2,
\end{array}
\right.
$$
satisfying $|D\phi(x)|\le 2$ and then taking the limit $R\rightarrow
+\infty$.

Now we estimate the terms on the right hand side of \eqref{e3}.
First, it holds that
\begin{equation}\label{WE1}
\begin{array}{ll}
\di |\int\r \f{|u|^2}{2}u\cdot\nabla(|x|^\a) dx|\di =\a|\int\r \f{|u|^2}{2}|x|^{\a-2}u\cdot x dx|\\
\di\quad \leq \f{\a}{2}\int\r |u|^3|x|^{\a-1} dx\leq \f{\a}{2}\|\sqrt\r u\|_2\|\sqrt\r\|_{p_1}\||x|^{\a-1}|u|^2\|_{q_1}\\
\di \quad\leq C\|\sqrt\r\|_{p_1}\||x|^{\f{\a-1}{2}}u\|^2_{2q_1}\leq C\|\sqrt\r\|_{p_1}\|\nabla u\|_2^{2\t_1}\||x|^{\b_1}u\|^{2(1-\t_1)}_{\bar q_1}\\
\di \quad \leq C\|\r\|^{\f12}_{\f{p_1}{2}}\|\nabla u\|_2^{2\t_1}\||x|^{\f{\a}{2}}\nabla u\|^{2(1-\t_1)}_{2} \leq \sigma\||x|^{\f{\a}{2}}\nabla
u\|^{2}_{2}+C_\sigma\|\r\|^{\f1{2\t_1}}_{\f{p_1}{2}}\|\nabla
u\|_2^2,
\end{array}
\end{equation}
where and in the sequel $\sigma >0$ is a small constant to be
determined, $C_\sigma$ is a positive constant depending on $\sigma$.
By the H${\rm\ddot{o}}$lder inequality and the
Caffarelli-Kohn-Nirenberg inequality in  Lemma \ref{lemma2} (1), the
positive constants $p_1>2, q_1>2, \bar q_1>1,\b_1>0, \t_1\in(0,1)$
in the above inequality  \eqref{WE1} satisfying
$$
\f{1}{p_1}+\f{1}{q_1}=\f12,
$$
$$
\f{1}{2q_1}+\f{\f{\a-1}{2}}{2}=\t_1(\f12+\f{0-1}{2})+(1-\t_1)(\f{1}{\bar
q_1}+\f{\b_1}{2}),
$$
and
$$
\f{1}{\bar q_1}+\f{\b_1}{2}=\f{1}{2}+\f{\f{\a}{2}-1}{2}.
$$
The combination of the above three equalities yields that
\begin{equation}\label{W-p1}
p_1=\f{2}{\a\t_1}.
\end{equation}
Note that one should choose the parameters $\a>0$ and $0<\t_1<1$
such that $p_1>2$ in \eqref{W-p1}.  Now choose $m>1$ sufficiently
large such that $2m\b+1>\f{p_1}{2}.$ Therefore, by the interpolation
inequality, it holds that
\begin{equation}\label{Im1}
\|\r\|_{\f{p_1}{2}}^{\f{1}{2\t_1}}\leq
\|\r\|_{1}^{\f{a_1}{2\t_1}}\|\r\|^{\f{1-a_1}{2\t_1}}_{2m\b+1},
\end{equation}
with $a_1\in(0,1)$ satisfying
$$
\f{a_1}{1}+\f{1-a_1}{2m\b+1}=\f{2}{p_1}=\t_1\a,
$$
which implies that
$$
a_1=\t_1\a\big(1+\f{1}{2m\b}\big)+\f{1}{2m\b}.
$$
To close the estimates in Lemma \ref{lemma-rho}, the
following restriction should be imposed to \eqref{Im1}
$$
\f{1-a_1}{2\t_1}\leq \b, ~~{\rm i.~e.}~~\t_1\geq \f{1-a_1}{2\b},
$$
For definiteness, we can choose $\t_1=\f12$ and then
$a_1=\f\a2-\f{2-\a}{4m\b}\in(0,1)$ . Obviously,  the restriction
$\t_1=\f12\geq \f{1}{2\b}(1-a_1)=\f{1}{2\b}(1-\f\a2+\f{2-\a}{4m\b})$
is satisfied if $m\gg1$. Then it follows from \eqref{Im1} that
\begin{equation}\label{Im1-1}
\|\r\|_{\f{p_1}{2}}^{\f{1}{2\t_1}}\leq C(\|\r\|^{\b}_{2m\b+1}+1),
\end{equation}
with the positive constant $C$ independent of $m$.

Then it holds that
\begin{equation}\label{WE2}
\begin{array}{ll}
\di |\int\f{\g\r^\g}{\g-1}u\cdot\nabla(|x|^\a) dx|\di =\f{\g\a}{\g-1}|\int\r^\g|x|^{\a-2}u\cdot x dx|\\
\di\quad \leq \f{\g\a}{\g-1}\int\r^\g |u||x|^{\a-1} dx\leq \f{\g\a}{\g-1}\|\r^\g\|_{p_2}\||x|^{\a-1}u\|_{q_2}\\
\di \quad\leq C\|\r^\g\|_{p_2}\|\nabla
u\|_2^{\t_2}\||x|^{\b_2}u\|^{1-\t_2}_{\bar q_2}
\leq C\|\r\|^{\g}_{ p_2\g}\|\nabla u\|_2^{\t_2}\||x|^{\f{\a}{2}}\nabla u\|^{1-\t_2}_{2}\\
\di \quad \leq \sigma\||x|^{\f{\a}{2}}\nabla u\|^{2}_{2}+C_\sigma\|\r\|^{\f{2\g}{1+\t_2}}_{p_2\g}\|\nabla u\|_2^{\f{2\t_2}{1+\t_2}} \leq \sigma\||x|^{\f{\a}{2}}\nabla
u\|^{2}_{2}+C_\sigma\|\r\|^{\f{2\g}{1+\t_2}}_{p_2\g}(\|\nabla
u\|_2^2+1).
\end{array}
\end{equation}
By the H${\rm\ddot{o}}$lder inequality and the
Caffarelli-Kohn-Nirenberg inequality in  Lemma \ref{lemma2} (1), the
positive constants $p_2>1, q_2>1, \bar q_2>1,\b_2>0, \t_2\in(0,1)$
in the above inequality \eqref{WE2} satisfying
$$
\f{1}{p_2}+\f{1}{q_2}=1,
$$
$$
\f{1}{q_2}+\f{\a-1}{2}=\t_2(\f12+\f{0-1}{2})+(1-\t_2)(\f{1}{\bar
q_2}+\f{\b_2}{2}),
$$
and
$$
\f{1}{\bar q_2}+\f{\b_2}{2}=\f{1}{2}+\f{\f{\a}{2}-1}{2}.
$$
The combination of the above three equalities yields that
\begin{equation}\label{W-p2}
p_2=\f{4}{2+\a(1+\t_2)}.
\end{equation}
Note that one should choose the parameters $\a>0$ and $0<\t_2<1$
such that $p_2>1$ in \eqref{W-p2}.  Now choose $m>1$ sufficiently
large such that $2m\b+1>p_2\g.$ Therefore, by the interpolation
inequality, it holds that
\begin{equation}\label{Im2}
\|\r\|_{p_2\g}^{\f{2\g}{1+\t_2}}\leq
\|\r\|_{\g}^{\f{2\g}{1+\t_2}a_2}\|\r\|^{\f{2\g}{1+\t_2}(1-a_2)}_{2m\b+1},
\end{equation}
with $a_2\in(0,1)$ satisfying
$$
\f{a_2}{\gamma}+\f{1-a_2}{2m\b+1}=\f{1}{p_2\g}=\f{2+\a(1+\t_2)}{4\g},
$$
which implies that
$$
a_2\rightarrow\f{2+\a(1+\t_2)}{4},~~{\rm as}~m\rightarrow+\i.
$$
The following restriction should be imposed to \eqref{Im2}
\begin{equation}\label{June8-2}
\f{2\g}{1+\t_2}(1-a_2)\leq \b, ~~{\rm i.~e.}~~1+\t_2\geq
\f{2\g}{\b}(1-a_2).
\end{equation}
For $m\gg1$ large enough, it is sufficient to have the following
restriction
$$
1+\t_2> \f{\g}{\b}(1-\f{\alpha(1+\t_2)}{2}),
$$
 That is
\begin{equation}\label{June8-1}
(1+\t_2)(\frac{\b}{\g}+\frac{\alpha}{2})>1.
\end{equation}

Consequently, if
\begin{equation}\label{June21-3}
1<\alpha\le 2, \quad \beta>0, \quad \gamma>1,
\end{equation}
 we can choose $0\le
\frac{2}{\alpha}-1<\t_2<1$ such that \eqref{June8-1} and hence
\eqref{June8-2} hold true for $m\gg1$.
If
\begin{equation}\label{June21-4}
0<\alpha\le 1, \quad \beta>\frac12, \quad 1<\g\le 2\b,
\end{equation}
 we can choose $\max\{
\frac{\g}{\b}-1,0\}<\t_2<1$ such that \eqref{June8-1} and hence
\eqref{June8-2} hold true for large $m\gg1$.

 Then it follows from \eqref{Im2} that
\begin{equation}\label{Im2-1}
\|\r\|_{p_2\g}^{\f{2\g}{1+\t_2}}\leq C(\|\r\|^{\b}_{2m\b+1}+1)
\end{equation}
with the positive constant $C$ independent of $m$.

Now one can compute that
\begin{equation}\label{WE3}
\begin{array}{ll}
\di |-\int \mu\nabla\f{|u|^2}{2}\cdot\nabla(|x|^\a)dx|=\mu\a|\int u\cdot\nabla u\cdot x|x|^{\a-2}dx|\\
\di\quad \leq \mu\a\||x|^{\f\a2}\nabla u\|_2\||x|^{\f\a2-1}
u\|_2\leq \f{\mu\a^2}{2}\||x|^{\f\a2}\nabla u\|_2^2,
\end{array}
\end{equation}
where in the last inequality one has used the best constant $\f\a2$
for the Caffarelli-Kohn-Nirenberg inequality in Lemma \ref{lemma2}
(2). Similarly, it holds that
\begin{equation}\label{WE4}
\begin{array}{ll}
\di |-\int \mu({\rm div}u)u\cdot\nabla(|x|^\a)dx|=\mu\a|\int ({\rm div}u) |x|^{\a-2}u\cdot x dx|\\
\di\quad \leq \mu\a\||x|^{\f\a2}{\rm div}u\|_2\||x|^{\f\a2-1}
u\|_2\leq \f{\mu\a^2}{2}\||x|^{\f\a2}{\rm
div}u\|_2\||x|^{\f\a2}\nabla u\|_2.
\end{array}
\end{equation}
Then it follows that
\begin{equation}\label{WE5}
\begin{array}{ll}
\di |-\int \l(\r)({\rm div}u)u\cdot\nabla(|x|^\a)dx|=\a|\int \l(\r)({\rm div}u) |x|^{\a-2}u\cdot x dx|\\
\di\quad \leq \a\|\sqrt{\l(\r)}|x|^{\f\a2}{\rm div}u\|_2\|\sqrt{\l(\r)}\|_{p_3}\||x|^{\f\a2-1} u\|_{q_3}\\
\di \quad \leq C\|\sqrt{\l(\r)}|x|^{\f\a2}{\rm div}u\|_2\|\r^{\f{\b}{2}}\|_{p_3}\|\nabla u\|_2^{\t_3}\||x|^{\b_3} u\|_{\bar q_3}^{1-\t_3}\\
\di \quad \leq C\|\sqrt{\l(\r)}|x|^{\f\a2}{\rm div}u\|_2\|\r\|^{\f{\b}{2}}_{\f{\b p_3}{2}}\|\nabla u\|_2^{\t_3}\||x|^{\f\a2} \nabla u\|_{2}^{1-\t_3}\\
\di \quad\leq \sigma\|\sqrt{\l(\r)}|x|^{\f\a2}{\rm div}u\|_2^2+C_\sigma\|\r\|^\b_{\f{\b p_3}{2}}\|\nabla u\|_2^{2\t_3}\||x|^{\f\a2} \nabla u\|_{2}^{2(1-\t_3)}\\
\di \quad\leq \sigma\big[\|\sqrt{\l(\r)}|x|^{\f\a2}{\rm
div}u\|_2^2+\||x|^{\f\a2} \nabla
u\|_{2}^{2}\big]+C_\sigma\|\r\|^{\f{\b}{\t_3}}_{\f{\b
p_3}{2}}\|\nabla u\|_2^{2},
\end{array}
\end{equation}
By the H${\rm\ddot{o}}$lder inequality and the
Caffarelli-Kohn-Nirenberg inequality in  Lemma \ref{lemma2} (1), the
positive constants $p_3>2, q_3>2, \bar q_3>1,\b_3>0, \t_3\in(0,1)$
in the above inequality \eqref{WE5} satisfying
$$
\f{1}{p_3}+\f{1}{q_3}=\f12,
$$
$$
\f{1}{q_3}+\f{\f{\a}{2}-1}{2}=\t_3(\f12+\f{0-1}{2})+(1-\t_3)(\f{1}{\bar
q_3}+\f{\b_3}{2}),
$$
and
$$
\f{1}{\bar q_3}+\f{\b_3}{2}=\f{1}{2}+\f{\f{\a}{2}-1}{2}.
$$
The combination of the above three equalities yields that
\begin{equation}\label{W-p3}
p_3=\f{4}{\a\t_3}.
\end{equation}
Note that one should choose the parameters $\a,\t_3\in(0,1)$ such
that $p_3>2$ in \eqref{W-p3}. By the interpolation inequality, it
holds that
\begin{equation}\label{Im3}
\|\r\|_{\f{p_3\b}{2}}^{\f{\b}{\t_3}}\leq
\|\r\|_{1}^{\f{\b}{\t_3}a_3}\|\r\|^{\f{\b}{\t_3}(1-a_3)}_{2m\b+1},
\end{equation}
with $a_3\in(0,1)$ satisfying
$$
\f{a_3}{1}+\f{1-a_3}{2m\b+1}=\f{2}{p_3\b}=\f{\a\t_3}{2\b},
$$
which implies that
$$
a_3=\f{\a\t_3}{2\b}+\f{\f{\a\t_3}{2\b}-1}{2m\b}.
$$
The following restriction should be imposed to \eqref{Im3}
$$
\f{\b}{\t_3}(1-a_3)\leq \b, ~~{\rm i.~e.}~~\t_3\geq (1-a_3).
$$
For definiteness, one can choose $\t_3\in (0,1)$
$$
\t_3\geq 1-\f{\a\t_3}{2\b}
$$
if $m$ is sufficiently large. Then it follows from \eqref{Im2} that
\begin{equation}\label{Im3-1}
\|\r\|_{\f{p_3\b}{2}}^{\f{\b}{\t_3}}\leq C(\|\r\|^{\b}_{2m\b+1}+1)
\end{equation}
with the positive constant $C$ independent of $m$.

Substituting \eqref{WE1}, \eqref{Im1-1}, \eqref{WE2}, \eqref{Im2-1},
\eqref{WE3}, \eqref{WE4}, \eqref{WE5} and \eqref{Im3-1} into
\eqref{e3} yields that
\begin{equation}\label{e4}
\begin{array}{ll}
\di \f{d}{dt}\int |x|^\a(\r\f{|u|^2}{2}+\f{\r^\g}{\g-1})(t,x)dx+\||x|^{\f\a2}\sqrt{\l(\r)}{\rm div} u\|_2^2(t) \\
\di +\mu\Big[(1-\f{\a^2}{2})\||x|^{\f\a2}\nabla u\|_2^2-\f{\a^2}{2}\||x|^{\f\a2}\nabla u\|_2\||x|^{\f\a2}{\rm div} u\|_2+\||x|^{\f\a2}{\rm div} u\|_2^2\Big]\\
\di \leq \sigma\big[\|\sqrt{\l(\r)}|x|^{\f\a2}{\rm
div}u\|_2^2+3\||x|^{\f\a2} \nabla
u\|_{2}^{2}\big]+C_\sigma\big(\|\r\|^\b_{2m\b+1}+1\big)(\|\nabla
u\|_2^{2}+1).
\end{array}
\end{equation}
The determinant of the quadratic term in the second line of
\eqref{e4} can be calculated by
$$
\Delta=\f{\a^4}{4}-4(1-\f{\a^2}{2})=\f14\big(\a^4+8\a^2-16\big).
$$
Therefore, if the weight $\a$ satisfies
\begin{equation}\label{alpha}
0<\a^2<4(\sqrt2-1),
\end{equation}
then the determinant $\Delta<0,$ and thus there exists a positive
constant $C_\a$ such that
\begin{equation}\label{quad}
(1-\f{\a^2}{2})\||x|^{\f\a2}\nabla
u\|_2^2-\f{\a^2}{2}\||x|^{\f\a2}\nabla u\|_2\||x|^{\f\a2}{\rm div}
u\|_2+\||x|^{\f\a2}{\rm div} u\|_2^2\geq
C^{-1}_\a\Big[\||x|^{\f\a2}\nabla u\|_2^2+\||x|^{\f\a2}{\rm div}
u\|_2^2\Big].
\end{equation}
Substituting \eqref{quad} into \eqref{e4}, choosing $\sigma$
suitably small in \eqref{e4} and noting that
 $$
 \int \r_0^\g |x|^\a dx\leq \|\r_0|x|^\a\|_{1}\|\r_0^{\g-1}\|_\i\leq C \|\r_0(1+|x|^{\a_1})\|_1\|\r_0\|^{\g-1}_{W^{2,q}(\mathbb{R}^2)}\leq C,
 $$
 yield the estimate \eqref{We-l1} in Lemma \ref{lemma-wee}. The restrictions of $\alpha$ \eqref{June21-1} and \eqref{June21-2}  follow from \eqref{June21-3}, \eqref{June21-4} and \eqref{alpha}. $\hfill\Box$

\underline{Step 4. Density estimates:}

Applying the operator $div$
to the momentum equation $\eqref{CNS}_2$, it holds that
\begin{equation}\label{vis-f}
[{\rm div}(\r u)]_t+{\rm div}[{\rm div}(\r u\otimes u)]=\Delta F.
\end{equation}
Consider the following two elliptic problems on the whole space $\mathbb{R}^2$:
\begin{equation}\label{xi}
\Delta \x={\rm div}(\r u),
\end{equation}
\begin{equation}\label{eta}
\Delta \eta={\rm div}[{\rm div}(\r u\otimes u)],
\end{equation}
both with the boundary conditions $\x,\eta\rightarrow 0$ as
$|x|\rightarrow\i$. By the elliptic estimates and H${\rm\ddot{o}}$lder inequality, it
holds that
\begin{Lemma}\label{lemma4}
\begin{itemize}
\item[(1)] $\|\nabla\x\|_{2m}\leq Cm\|\r\|_{\f{2mk}{k-1}}\|u\|_{2mk},$
for any $k>1,m\geq1;$
\item[(2)] $\|\nabla\x\|_{2-r}\leq C\|\sqrt\r u\|_{2}\|\r\|^{\f12}_{\f{2-r}{r}}\leq C\|\r\|^{\f12}_{\f{2-r}{r}},$
for any $0< r<1;$
\item[(3)] $\|\eta\|_{2m}\leq Cm\|\r\|_{\f{2mk}{k-1}}\|u\|^2_{4mk},$
for any $k>1,m\geq1;$
\end{itemize} where $C$ are positive constants independent of $m,k$
and $r$.
\end{Lemma}
{\bf Proof:} (1) By the elliptic estimates to the equation
\eqref{xi} and then using the H${\rm\ddot{o}}$lder inequality, one has for any $k>1,m\geq1$,
$$
\|\nabla\x\|_{2m}\leq C m\|\r u\|_{2m}\leq
Cm\|\r\|_{\f{2mk}{k-1}}\|u\|_{2mk}.
$$
Similarly, the statements (2) and (3) can be proved. $\hfill\Box$

Based on Lemmas \ref{lemma1}-\ref{lemma3} and Lemma \ref{lemma4}, it holds that

\begin{Lemma}\label{lemma5}
\begin{itemize}
\item[(1)] $\|\x\|_{2m}\leq Cm^{\f12}\|\nabla\x\|_{\f{2m}{m+1}}\leq Cm^{\f12}\|\r\|_m^{\f12},$
for any $m\geq2;$
\item[(2)] $\|u\|_{2m}\leq C m^{\f12}\|\nabla u\|_{2}^{1-\f{1}{m\a}}\||x|^{\f\a2}\nabla u\|_2^{\f{1}{m\a}},$
for any $m+1\geq\f{4}{\a};$
\item[(3)] $\|\nabla\x\|_{2m}\leq C m^{\f32}k^{\f12}\|\r\|_{\f{2mk}{k-1}}\|\nabla u\|_2^{1-\f{2}{mk\a}}\||x|^{\f\a2}\nabla u\|_2^{\f{2}{mk\a}},$
for any $k>1, m+1\geq\f{4}{\a};$
\item[(4)] $\|\eta\|_{2m}\leq C m^2k\|\r\|_{\f{2mk}{k-1}}\|\nabla u\|_2^{2-\f{2}{mk\a}}\||x|^{\f\a2}\nabla u\|_2^{\f{2}{mk\a}},$
for any $k>1, m+1\geq\f{4}{\a};$
\end{itemize} where $C$ are positive constants independent of $m,k$.
\end{Lemma}
{\bf Proof:} (1) By Lemma \ref{lemma2} and Lemma \ref{lemma4} (2),
it holds that
$$
\|\x\|_{2m}\leq Cm^{\f12}\|\nabla\x\|_{\f{2m}{m+1}}\leq
Cm^{\f12}\|\sqrt\r u\|_2\|\r\|_m^{\f12}\leq Cm^{\f12}\|\r\|_m^{\f12},
$$
where in the last inequality one has used the elementary energy
estimates \eqref{Eee}.

(2) If $m+1>\f{4}{\a}$, then by interpolation inequality and
Caffarelli-Kohn-Nirenberg inequality, it holds that
\begin{equation}\label{ine}
\|u\|_{m+1}\leq \|u\|^\t_{2m}\|u\|_{\f{4}{\a}}^{1-\t}\leq C\|u\|^\t_{2m}\||x|^{\f\a2}\nabla u\|_{2}^{1-\t}
\end{equation}
where
$$
\t=\f{\f{1}{m+1}-\f\a4}{\f{1}{2m}-\f\a4}.
$$
Then it follows from Lemma \ref{lemma1} (2) and \eqref{ine} that
$$
\|u\|_{2m}\leq Cm^{\f14}\|\nabla u\|_2^{\f12-\f{1}{2m}}\|u\|_{m+1}^{\f12+\f{1}{2m}}\leq Cm^{\f14}\|\nabla u\|_2^{\f12-\f{1}{2m}}\|u\|_{2m}^{(\f12+\f{1}{2m})\t}\||x|^{\f\a2}\nabla u\|_{2}^{(\f12+\f{1}{2m})(1-\t)},
$$
which implies Lemma \ref{lemma5} (2) immediately.

The assertions (3) and (4) in Lemma \ref{lemma5} are the
direct consequences of Lemma \ref{lemma5} (2) and Lemma
\ref{lemma4} (1), (3), respectively. Thus the proof of Lemma \ref{lemma5} is completed. $\hfill\Box$

Substituting \eqref{xi} and \eqref{eta} into \eqref{vis-f} yields
that
$$
\Delta\Big(\x_t+\eta-F\Big)=0,
$$
which implies that
\begin{equation*}
\x_t+\eta-F=0.
\end{equation*}
It follows from the definition \eqref{flux} of the effective viscous flux $F$ that
\begin{equation*}
\x_t-(2\mu+\l(\r)){\rm div}u+P(\r)+\eta=0.
\end{equation*}
Then the continuity equation $\eqref{CNS}_1$ yields that
\begin{equation*}
\x_t+\f{2\mu+\l(\r)}{\r}(\r_t+u\cdot\nabla \r)+P(\r)+\eta=0.
\end{equation*}
Define
\begin{equation}\label{theta}
\nu(\r)=\int_1^\r\f{2\mu+\l(s)}{s}ds=2\mu\ln\r+\f{1}{\b}(\r^\b-1).
\end{equation}
Then we obtain a new transport equation
\begin{equation}\label{transport-e}
(\x+\nu(\r))_t+u\cdot\nabla(\x+\nu(\r))+P(\r)+\eta-u\cdot\nabla\x=0,
\end{equation}
which is crucial in the following Lemma for the higher integrability of the density function.
\begin{Lemma}\label{lemma-rho}
  For any $k\geq1,$ it holds that
  \begin{equation}\label{density-e}
    \sup_{t\in[0,T]}\|\r(t,\cdot)\|_k\leq C(M)~ k^{\f{2}{\b-1}}.
  \end{equation}
\end{Lemma}
{\bf Proof:} Multiplying the equation \eqref{transport-e} by
$\r[(\x+\nu(\r))_+]^{2m-1}$ with $m$ sufficiently large integer, here and in what follows, the
notation $(\cdots)_+$ denotes the positive part of $(\cdots)$, one can get that
\begin{equation}\label{trans-e1}
\begin{array}{ll}
\di \f{1}{2m}\f{d}{dt}\int\r[(\x+\nu(\r))_+]^{2m}dx+\int\r
P(\r)[(\x+\nu(\r))_+]^{2m-1}dx=-\int\r \eta[(\x+\nu(\r))_+]^{2m-1}dx\\
\di  +\int\r u\cdot\nabla \x[(\x+\nu(\r))_+]^{2m-1}dx.
\end{array}
\end{equation}
Denote
\begin{equation*}
f(t)=\big\{\int\r[(\x+\nu(\r))_+]^{2m} dx\big\}^{\f{1}{2m}},\qquad
t\in[0,T].
\end{equation*}
Now we estimate the two terms on the right hand side of
\eqref{trans-e1}. First, it holds that
\begin{equation}\label{rho-1}
\begin{array}{ll}
\di |-\int\r \eta[(\x+\nu(\r))_+]^{2m-1}dx|\leq
\int\r^{\f{1}{2m}}|\eta|\big[\r(\x+\nu(\r))^{2m}_+\big]^{\f{2m-1}{2m}}dx\\
\qquad\di\leq
\|\r\|_{2m\b+1}^{\f{1}{2m}}\|\eta\|_{2m+\f{1}{\b}}\|\r(\x+\nu(\r))^{2m}_+\|_1^{\f{2m-1}{2m}}\\
\qquad\di\leq
C\|\r\|_{2m\b+1}^{\f{1}{2m}}(m+\f{1}{2\b})^2k\|\r\|_{\f{2(m+\f{1}{2\b})k}{k-1}}\|\nabla u\|_2^{2-\f{2}{(m+\f{1}{2\b})k\a}}\||x|^{\f\a2}\nabla u\|_2^{\f{2}{(m+\f{1}{2\b})k\a}}f(t)^{2m-1}\\
\qquad\di\leq
Cm^2\|\r\|_{2m\b+1}^{1+\f{1}{2m}}\|\nabla u\|_2^{2-\f{4(\b-1)}{\a(2m\b+1)}}\||x|^{\f\a2}\nabla u\|_2^{\f{4(\b-1)}{\a(2m\b+1)}}f(t)^{2m-1},
\end{array}
\end{equation}
where in the last
inequality we have taken $k=\f{\b}{\b-1}.$
Next, for $\f{1}{2m\b+1}+\f1p+\f1q=1$ with $p,q\geq 1$, one has
\begin{equation}\label{rho-2}
\begin{array}{ll}
\di |\int\r u\cdot\nabla\x[(\x+\nu(\r))_+]^{2m-1}dx|\leq
\int\r^{\f{1}{2m}}|u||\nabla\x|\big[\r(\x+\nu(\r))^{2m}_+\big]^{\f{2m-1}{2m}}dx\\
\qquad\di\leq
\|\r\|_{2m\b+1}^{\f{1}{2m}}\|u\|_{2mp}\|\nabla\x\|_{2mq}\|\r(\x+\nu(\r))^{2m}_+\|_1^{\f{2m-1}{2m}}\\
\qquad\di\leq
C\|\r\|_{2m\b+1}^{\f{1}{2m}}(mp)^{\f12}\|\nabla u\|_2^{1-\f{2}{mp\a}}\||x|^{\f\a2}\nabla u\|_2^{\f{2}{mp\a}}\\
\qquad\qquad\qquad\qquad\qquad
\qquad \di \cdot(mq)^{\f32}k^{\f12}\|\r\|_{\f{2mqk}{k-1}}\|\nabla u\|_2^{1-\f{2}{mqk\a}}\||x|^{\f\a2}\nabla u\|_2^{\f{2}{mqk\a}}f(t)^{2m-1}\\
\qquad\di\leq
Cm^2\|\r\|_{2m\b+1}^{1+\f{1}{2m}}\|\nabla u\|_2^{2-\f{4(\b-1)}{\a(2m\b+1)}}\||x|^{\f\a2}\nabla u\|_2^{\f{4(\b-1)}{\a(2m\b+1)}}f(t)^{2m-1},
\end{array}
\end{equation}
where in the third inequality one has chosen $p=q=\f{2m\b+1}{m\b}$
and $k=\f{\b}{\b-2}.$
Substituting \eqref{rho-1} and \eqref{rho-2} into
\eqref{trans-e1} yields that
\begin{equation*}
\begin{array}{ll}
\di\f{1}{2m}\f{d}{dt}(f^{2m}(t))+\int\r P(\r)[(\x+\nu(\r))_+]^{2m-1}dx\\
\di\leq
Cm^2\|\r\|_{2m\b+1}^{1+\f{1}{2m}}\|\nabla u\|_2^{2-\f{4(\b-1)}{\a(2m\b+1)}}\||x|^{\f\a2}\nabla u\|_2^{\f{4(\b-1)}{\a(2m\b+1)}}f(t)^{2m-1}.
\end{array}
\end{equation*}
Then it holds that
\begin{equation*}
\f{d}{dt}f(t)\leq
Cm^2\|\r\|_{2m\b+1}^{1+\f{1}{2m}}\|\nabla u\|_2^{2-\f{4(\b-1)}{\a(2m\b+1)}}\||x|^{\f\a2}\nabla u\|_2^{\f{4(\b-1)}{\a(2m\b+1)}}.
\end{equation*}
Integrating the above inequality over $[0,t]$ gives that
\begin{equation}\label{f}
f(t)\leq
f(0)+Cm^2\int_0^t\|\r\|_{2m\b+1}^{1+\f{1}{2m}}\|\nabla u\|_2^{2-\f{4(\b-1)}{\a(2m\b+1)}}\||x|^{\f\a2}\nabla u\|_2^{\f{4(\b-1)}{\a(2m\b+1)}}d\tau.
\end{equation}
Completely similar to the proof of Lemma 3.4 in \cite{JWX2}, we
obtain
\begin{equation}\label{f1}
f(t)\leq C\Big[1+m^2\int_0^t\|\r\|_{2m\b+1}^{1+\f{1}{2m}}\|\nabla
u\|_2^{2-\f{4(\b-1)}{\a(2m\b+1)}}\||x|^{\f\a2}\nabla
u\|_2^{\f{4(\b-1)}{\a(2m\b+1)}}d\tau\Big],
\end{equation}
and
\begin{equation}\label{rs2}
\begin{array}{ll}
\|\r\|_{2m\b+1}^\b(t)&\di\leq
C\Big[1+f(t)+m^{\f12}\|\r\|_{2m\b+1}^{\f12}(t)\Big]\\
&\di \leq
\f12\|\r\|_{2m\b+1}^\b(t)+C\Big[1+f(t)+m^{\f{\b}{2\b-1}}\Big].
\end{array}
\end{equation}
It follows from \eqref{f1}, \eqref{rs2} and Lemma \ref{lemma-wee}
that
\begin{equation*}
\begin{array}{ll}
\|\r\|_{2m\b+1}^\b(t)\di\leq C\Big[f(t)+m^{\f{\b}{2\b-1}}\Big]\\
\di\leq
C\Big[m^{\f{\b}{2\b-1}}+m^2\int_0^t\|\r\|_{2m\b+1}^{1+\f{1}{2m}}\|\nabla u\|_2^{2-\f{4(\b-1)}{\a(2m\b+1)}}\||x|^{\f\a2}\nabla u\|_2^{\f{4(\b-1)}{\a(2m\b+1)}}d\tau\Big]\\
\di\leq
C\Big[m^{\f{\b}{2\b-1}}+\int_0^t\||x|^{\f\a2}\nabla u\|_2^2(\tau)d\tau+m^{\f{2\a(2m\b+1)}{\a(2m\b+1)-2(\b-1)}}\int_0^t\|\r\|_{2m\b+1}^{(1+\f{1}{2m})\f{\a(2m\b+1)}{\a(2m\b+1)-2(\b-1)}}\|\nabla u\|_2^2(\tau)d\tau\Big]\\
\di \leq C\Big[m^{\f{\b}{2\b-1}}+\int_0^t\big(\|\r\|^\b_{2m\b+1}(\tau)+1\big)\big(\|\nabla u\|_2^2(\tau)+1\big)d\tau\\
\di \qquad\qquad\qquad\qquad\qquad\qquad\qquad\qquad\qquad +m^2\int_0^t\|\r\|_{2m\b+1}^{(1+\f{1}{2m})\f{\a(2m\b+1)}{\a(2m\b+1)-2(\b-1)}}\|\nabla u\|_2^2(\tau)d\tau\Big].
\end{array}
\end{equation*}
Applying Gronwall's inequality to the above inequality yields that
\begin{equation*}
\|\r\|_{2m\b+1}^\b(t)\leq
C\Big[m^{\f{\b}{2\b-1}}+m^2\int_0^t\|\r\|_{2m\b+1}^{(1+\f{1}{2m})\f{\a(2m\b+1)}{\a(2m\b+1)-2(\b-1)}}\|\nabla u\|_2^2(\tau)d\tau\Big].
\end{equation*}
Denote
$$
y(t)=m^{-\f{2}{\b-1}}\|\r\|_{2m\b+1}(t).
$$
Then it holds that
$$
\begin{array}{ll}
\di y^\b(t)&\di \leq
C\Big[m^{\f{\b(1-3\b)}{(2\b-1)(\b-1)}}+\int_0^ty(\tau)^{(1+\f{1}{2m})[1+\f{2(\b-1)}{\a(2m\b+1)-2(\b-1)}]}\|\nabla u\|_2^2(\tau)d\tau\Big]\\
&\di \leq C\Big[1+\int_0^t\big(y^\b(\tau)+1\big)\|\nabla u\|_2^2(\tau) d\tau\Big].
\end{array}
$$
So applying the Gronwall's inequality the above inequality yields that
$$
y(t)\leq C,\quad \forall t\in[0,T],
$$
that is,
$$
\|\r\|_{2m\b+1}(t)\leq Cm^{\f{2}{\b-1}},\quad \forall t\in[0,T].
$$
Equivalently,  \eqref{density-e} holds. Thus Lemma \ref{lemma-rho} is proved. $\hfill\Box$

\underline{Step 4: First-order derivative estimates of the
velocity.}

\begin{Lemma}\label{lemma-u-der}
  There exists a positive constant $C=C(M)$, such that
  \begin{equation*}
    \sup_{t\in[0,T]}\int(\mu\omega^2+\f{F^2}{2\mu+\l(\r)}) dx+\int_0^T\int\r(H^2+L^2)
     dxdt\leq C.
  \end{equation*}
\end{Lemma}
{\bf Proof:} The proof of the lemma is similar to that of Lemma 3.5
in \cite{JWX2} besides the weighted estimates. Multiplying the equation $\eqref{F-omega}_1$ by
$\mu\omega$, the equation $\eqref{F-omega}_2$ by
$\f{F}{2\mu+\l(\r)}$, respectively, and then summing and integrating
the resulted equations with respect to $x\in\mathbb{R}^2$, one has
\begin{equation}\label{ue1}
\begin{array}{ll}
\di\f12\f{d}{dt}\int(\mu\omega^2+\f{F^2}{2\mu+\l(\r)}) dx+\int\r(H^2+L^2) dx\\
\di \quad =-\f{\mu}{2}\int\omega^2{\rm div }u dx+\f12\int F^2({\rm
div}u)\Big[\r(\f{1}{2\mu+\l(\r)})^\prime-\f{1}{2\mu+\l(\r)}\Big]
dx\\
\di \quad +\int F({\rm div}
u)\Big[\r(\f{P(\r)}{2\mu+\l(\r)})^\prime-\f{P(\r)}{2\mu+\l(\r)}\Big]
dx-\int 2F(u_{1x_2}u_{2x_1}-u_{1x_1}u_{2x_2})dx.
\end{array}
\end{equation}
Set
$$
Z^2(t)=\int(\mu\omega^2+\f{F^2}{2\mu+\l(\r)}) dx,
$$
$$
\varphi^2(t)=\int\r(H^2+L^2)
dx=\int\frac{1}{\r}\big[(\mu\o_{x_1}+F_{x_2})^2+(-\mu\o_{x_2}+F_{x_1})^2\big]
dx.
$$
Then similar to the proof of Lemma 3.5 in \cite{JWX}, it  yields
that
\begin{equation}\label{ue6}
\begin{array}{ll}
\di \f12\f{d}{dt}Z^2(t)+\varphi(t)^2\leq
\s\varphi(t)^2+C_\s(Z(t)^2+1)^{2+\f{\v}{1-3\v}}\v^{\f{2}{1-\b}\f{1-\v}{1-3\v}}\\
\di \qquad\qquad\qquad\qquad\qquad +C\left[1+\int
\f{|F|^3}{2\mu+\l(\r)}dx+\int|F||\nabla u|^2 dx\right].
\end{array}
\end{equation}
Now it remains to estimate the terms $\di \int
\f{|F|^3}{2\mu+\l(\r)}dx$ and $\di\int|F||\nabla u|^2 dx $ on the
right hand side of \eqref{ue6}.
By Lemma \ref{lemma3}, for
$\v\in[0,\f12]$ and $\eta=\v$, it holds that
\begin{equation}\label{Fe1}
\|F\|_{2m}\leq Cm^{\f12}\|\nabla
F\|_{\f{2m}{m+\v}}^{1-s}\|F\|^s_{2(1-\v)},
\end{equation}
where $\di s=\f{(1-\v)^2}{m-\v(1-\v)}$ and the positive constant $C$
is independent of $m$ and $\v.$

Note that if $\v$ is sufficiently small, then it holds that
\begin{equation}\label{F3-1}
\begin{array}{ll}
\|F\|_{2(1-\v)}=\|(2\mu+\l(\r)){\rm div}u-P(\r)\|_{2(1-\v)}\\
\di ~~\leq 2\mu\|{\rm div} u\|_{2(1-\v)}+\|\l(\r){\rm div} u\|_{2(1-\v)}+\|P(\r)\|_{2(1-\v)}\\
\di ~~\leq 2\mu\|(1+|x|^{\f\a2}){\rm div} u\|_2\|(1+|x|^{\f\a2})^{-1}\|_{\f{2(1-\v)}{\v}}+\|\sqrt{\l(\r)}{\rm div} u\|_2\|\sqrt{\l(\r)}\|_{\f{2(1-\v)}{\v}}+C\\
\di ~~\leq C\Big[\|(1+|x|^{\f\a2}){\rm div} u\|_2+\|\sqrt{\l(\r)}{\rm div} u\|_2\|\r\|^{\f\b2}_{\f{\b(1-\v)}{\v}}+1\Big]\\
\di ~~\leq C\Big[\|(1+|x|^{\f\a2}){\rm div} u\|_2+\|\sqrt{\l(\r)}{\rm div} u\|_2\big(\f{\b(1-\v)}{\v}\big)^{\f{\b}{\b-1}}+1\Big]\\
\di~~\leq C\Big[\||x|^{\f\a2}{\rm div}
u\|_2+(Z(t)+1)(\v^{-\f{\b}{\b-1}}+1)\Big],
\end{array}
\end{equation}
where in the third inequality we have used the fact that
$$
\di \|(1+|x|^{\f\a2})^{-1}\|_{\f{2(1-\v)}{\v}}<+\i,
$$
provided that $\v$ is sufficiently small. By \eqref{Fe1},
\eqref{F3-1} and setting $\v=2^{-m}$ with $m$ sufficiently large, it
holds that
\begin{equation}\label{Fe}
\begin{array}{ll}
\di \|F\|_{2m}&\di \leq
Cm^{\f12}(\f{m+\v}{\v})^{\f{1-s}{\b-1}}\varphi(t)^{1-s}\Big[\||x|^{\f\a2}{\rm
div} u\|_2+(Z(t)+1)(\v^{-\f{\b}{\b-1}}+1)\Big]^s\\
&\di\leq
Cm^{\f12}(\f{m+\v}{\v})^{\f{1-s}{\b-1}}\varphi(t)^{1-s}\Big[\||x|^{\f\a2}{\rm
div} u\|^s_2+(Z(t)^s+1)(2^{\f{ms\b}{\b-1}}+1)\Big]\\
&\di\leq
Cm^{\f12}(\f{m}{\v})^{\f{1-s}{\b-1}}\varphi(t)^{1-s}\Big[\||x|^{\f\a2}{\rm
div} u\|^s_2+Z(t)^s+1\Big],
\end{array}
\end{equation}
with $\di s=\f{(1-\v)^2}{m-\v(1-\v)}$.   Then it follows that
\begin{equation}\label{F3}
\begin{array}{ll}
\di \int \f{|F|^3}{2\mu+\l(\r)}dx\di
=\int\f{|F|^{2-\f{1}{m-1}}}{(2\mu+\l(\r))^{1-\f{1}{2(m-1)}}}(\f{1}{2\mu+\l(\r)})^{\f{1}{2(m-1)}}|F|^{1+\f{1}{m-1}}
dx\\
\di \leq
\int\left(\f{|F|^2}{2\mu+\l(\r)}\right)^{1-\f{1}{2(m-1)}}|F|^{\f{m}{m-1}}
dx \leq \left(\int
\f{|F|^2}{2\mu+\l(\r)}dx\right)^{\f{2m-3}{2(m-1)}}\left(\int
|F|^{2m}dx\right)^{\f{1}{2(m-1)}}\\
\di =Z(t)^{\f{2m-3}{m-1}}\|F\|_{2m}^{\f{m}{m-1}}\leq
Cm^{\f{m}{2(m-1)}}(\f{m}{\v})^{\f{(1-s)m}{(\b-1)(m-1)}}Z(t)^{\f{2m-3}{m-1}}\varphi(t)^{\f{(1-s)m}{m-1}}\Big[\||x|^{\f\a2}{\rm div} u\|^s_2+Z(t)^s+1\Big]^{\f{m}{m-1}}\\
\di \leq
Cm^{\f12}(\f{m}{\v})^{\f{1}{\b-1}}Z(t)^{\f{2m-3}{m-1}}\varphi(t)^{\f{(1-s)m}{m-1}}\Big[\||x|^{\f\a2}{\rm div} u\|_2^{\f{ms}{m-1}}+Z(t)^{\f{ms}{m-1}}+1\Big]\\
\di \leq
\s\varphi(t)^2+C_\s~m(\f{m}{\v})^{\f{2}{\b-1}}\Big[(Z(t)^2+1)^{2+\f{1-ms}{m(1+s)-2}}+(Z(t)^2+1)^{2+\f{1-2ms}{m(1+s)-2}}\||x|^{\f\a2}{\rm div} u\|_2^{\f{2ms}{m(1+s)-2}}\Big]\\
\di \leq
\s\varphi(t)^2+C_\s~m(\f{m}{\v})^{\f{2}{\b-1}}\Big[(Z(t)^2+1)^{2+\f{1-ms}{m(1+s)-2}}\\
\qquad\qquad\qquad\qquad\qquad\qquad\qquad\qquad\di+(Z(t)^2+1)\||x|^{\f\a2}{\rm
div} u\|_2^2+(Z(t)^2+1)^{2+\f{1-ms}{m-2}}\Big]
\end{array}
\end{equation}
where one has used the fact that
$ms=\f{m(1-\v)^2}{m-\v(1-\v)}\rightarrow1$ with $\v=2^{-m}$ as
$m\rightarrow +\i$.

Furthermore, it holds that
\begin{equation}\label{F-nabla-u}
\begin{array}{ll}
\di \int|F||\nabla u|^2 dx&\di \leq \|F\|_{2m}\|\nabla
u\|_{\f{4m}{2m-1}}^2 \leq C\|F\|_{2m}\Big(\|{\rm div}
u\|_{\f{4m}{2m-1}}^2+\|\o\|_{\f{4m}{2m-1}}^2\Big)\\
&\di \leq
C\|F\|_{2m}\Big(\|\f{F}{2\mu+\l(\r)}\|_{\f{4m}{2m-1}}^2+\|\o\|_{\f{4m}{2m-1}}^2+1\Big).
\end{array}
\end{equation}
Note that
\begin{equation}\label{F-nabla-u1}
\begin{array}{ll}
\di
\|\f{F}{2\mu+\l(\r)}\|_{\f{4m}{2m-1}}^2&\di =\left(\int\f{|F|^{\f{2m(2m-3)}{(2m-1)(m-1)}}}{(2\mu+\l(\r))^{\f{4m}{2m-1}}}|F|^{\f{2m}{(2m-1)(m-1)}}dx\right)^{\f{2m-1}{2m}}\\
&\di \leq
\|F\|_{2m}^{\f{1}{m-1}}\left(\int\f{|F|^2}{(2\mu+\l(\r))^{\f{4(m-1)}{2m-3}}}dx\right)^{\f{2m-3}{2(m-1)}}\\
&\di \leq
C\|F\|_{2m}^{\f{1}{m-1}}\left(\int\f{|F|^2}{2\mu+\l(\r)}dx\right)^{\f{2m-3}{2(m-1)}}
= C\|F\|_{2m}^{\f{1}{m-1}}Z(t)^{\f{2m-3}{m-1}},
\end{array}
\end{equation}
and from  Lemma \ref{lemma1} (1), one has
\begin{equation}\label{F-nabla-u2}
\begin{array}{ll}
\di
\|\o\|_{\f{4m}{2m-1}}^2&\di \leq C\|\o\|_2^{2-\f{1-\v}{m(1-2\v)}}\|\nabla\o\|_{2(1-\v)}^{\f{1-\v}{m(1-2\v)}}\leq C Z(t)^{2-\f{1-\v}{m(1-2\v)}}\Big[\v^{\f{1}{1-\b}}\varphi(t)\Big]^{\f{1-\v}{m(1-2\v)}}\\
&\di \leq C
2^{\f{(1-\v)}{(\b-1)(1-2\v)}}Z(t)^{2-\f{1-\v}{m(1-2\v)}}\varphi(t)^{\f{1-\v}{m(1-2\v)}}
\leq C Z(t)^{2-\f{1-\v}{m(1-2\v)}}\varphi(t)^{\f{1-\v}{m(1-2\v)}}.
\end{array}
\end{equation}
Substituting \eqref{Fe} into \eqref{F-nabla-u1}, then
substituting the resulted \eqref{F-nabla-u1} and \eqref{F-nabla-u2}
into \eqref{F-nabla-u} give that
\begin{equation}\label{F-nabla-u3}
\begin{array}{ll}
&\di \int |F||\nabla u|^2dx \leq
C\Big\{m^{\f12}(\f{m}{\v})^{\f{1-s}{\b-1}}\varphi(t)^{1-s}\Big[\||x|^{\f\a2}{\rm
div} u\|^s_2+Z(t)^s+1\Big]\Big\}^{1+\f{1}{m-1}}Z(t)^{\f{2m-3}{m-1}}
\\
&\di~~~+Cm^{\f12}(\f{m}{\v})^{\f{1-s}{\b-1}}\varphi(t)^{1-s}\Big[\||x|^{\f\a2}{\rm
div} u\|^s_2+Z(t)^s+1\Big]\Big[Z(t)^{2-\f{1-\v}{m(1-2\v)}}\varphi(t)^{\f{1-\v}{m(1-2\v)}}+1\Big]\\
&\di\leq
Cm^{\f12}(\f{m}{\v})^{\f{(1-s)m}{(\b-1)(m-1)}}\varphi(t)^{\f{(1-s)m}{m-1}}Z(t)^{\f{2m-3}{m-1}}\Big[\||x|^{\f\a2}{\rm
div}
u\|_2^\f{ms}{m-1}+Z(t)^{\f{ms}{m-1}}+1\Big]\\
&\di~+Cm^{\f12}(\f{m}{\v})^{\f{1-s}{\b-1}}\varphi(t)^{1-s+\f{1-\v}{m(1-2\v)}}Z(t)^{2-\f{1-\v}{m(1-2\v)}}\Big[\||x|^{\f\a2}{\rm
div} u\|^s_2+Z(t)^s+1\Big]\\
&\di~+Cm^{\f12}(\f{m}{\v})^{\f{1-s}{\b-1}}\varphi(t)^{1-s}\Big[\||x|^{\f\a2}{\rm
div} u\|^s_2+Z(t)^s+1\Big]\\
&\di \leq
\s\varphi(t)^2+C_\s\bigg\{\Big[m^{\f12}(\f{m}{\v})^{\f{(1-s)m}{(\b-1)(m-1)}}Z(t)^{\f{2m-3}{m-1}}\Big(\||x|^{\f\a2}{\rm
div}
u\|_2^\f{ms}{m-1}+Z(t)^{\f{ms}{m-1}}+1\Big)\Big]^{\f{2(m-1)}{m(s+1)-2}}\\
&\di~+\Big[m^{\f12}(\f{m}{\v})^{\f{1-s}{\b-1}}Z(t)^{2-\f{1-\v}{m(1-2\v)}}\Big(\||x|^{\f\a2}{\rm
div} u\|^s_2+Z(t)^s+1\Big)\Big]^{\f{2}{1+s-\f{1-\v}{m(1-2\v)}}}\\
&\di~+\Big[m^{\f12}(\f{m}{\v})^{\f{1-s}{\b-1}}\Big(\||x|^{\f\a2}{\rm
div} u\|^s_2+Z(t)^s+1\Big)\Big]^{\f{2}{1+s}}\bigg\}\\
&\di \leq
\s\varphi(t)^2+C_\s~m(\f{m}{\v})^{\f{2}{\b-1}}\bigg[(1+Z(t)^2)^{2+\f{1-2ms}{m(s+1)-2}}\big(\||x|^{\f\a2}{\rm
div}
u\|^{\f{2ms}{m(s+1)-2}}_2+Z(t)^{\f{2ms}{m(s+1)-2}}+1\big)\\
&\di
~+(1+Z(t)^2)^{2+\f{1-2ms+(4ms-1)\v}{(1+s)m(1-2\v)-1+\v}}\big(\||x|^{\f\a2}{\rm
div}
u\|^{\f{2ms(1-2\v)}{(1+s)m(1-2\v)-1+\v}}_2+Z(t)^{\f{2ms(1-2\v)}{(1+s)m(1-2\v)-1+\v}}+1\big)\\
&\di~+\Big(\||x|^{\f\a2}{\rm
div} u\|^{\f{2s}{1+s}}_2+Z(t)^{\f{2s}{1+s}}+1\Big)\bigg]\\
&\di \leq
\s\varphi(t)^2+C_\s~m(\f{m}{\v})^{\f{2}{\b-1}}\bigg[(\||x|^{\f\a2}{\rm
div}
u\|^2_2+1)(Z(t)^2+1)+(1+Z(t)^2)^{2+\f{1-ms+(2ms-1)\v}{m(1+s)(1-2\v)-1+\v}}\\
&\di ~\qquad\qquad+(1+Z(t)^2)^{2+\f{1-ms-\v}{m(1-2\v)-1+\v}}
+(1+Z(t)^2)^{2+\f{1-ms}{m(s+1)-2}}+(1+Z(t)^2)^{2+\f{1-ms}{m-2}}\bigg].
\end{array}
\end{equation}
Substituting \eqref{F3} and \eqref{F-nabla-u3} into \eqref{ue6} and
choosing $\s$ sufficiently small yield that
\begin{equation}\label{Z0}
\begin{array}{ll}
\di \f12\f{d}{dt}(Z^2(t))+\f12\varphi(t)^2\leq
C~\v^{\f{2}{1-\b}}(Z(t)^2+1)^{2+\f{\v}{1-3\v}} \\
\di\qquad +Cm(\f{m}{\v})^{\f{2}{\b-1}}\Big[(\||x|^{\f\a2}{\rm div}
u\|^2_2+1)(Z(t)^2+1)+(1+Z(t)^2)^{2+\f{1-ms+(2ms-1)\v}{m(1+s)(1-2\v)-1+\v}}\\
\di\qquad\qquad
+(1+Z(t)^2)^{2+\f{1-ms-\v}{m(1-2\v)-1+\v}}+(1+Z(t)^2)^{2+\f{1-ms}{m(s+1)-2}}+(Z(t)^2+1)^{2+\f{1-ms}{m-2}}\Big].
\end{array}
\end{equation}
Note that $\lim_{m\rightarrow+\i}[2^m(1-ms)]=2$, and so $1-ms\sim
2\v$ as $m\rightarrow+\i$. Thus for $m$ sufficiently large, one has
$$
\f{1-ms+(2ms-1)\v}{m(1+s)(1-2\v)-1+\v}\sim
\f{2\v+\v(1-4\v)}{(m+1-2\v)(1-2\v)-1+\v}\leq 3\v,
$$
$$
\f{1-ms-\v}{m(1-2\v)-1+\v}\sim\f{\v}{m(1-2\v)-1+\v}\leq \v,
$$
$$
\f{1-ms}{m(s+1)-2}\sim \f{2\v}{1-2\v+m-2}=\f{2\v}{m-1-2\v}\leq 2\v,
$$
and
$$
\f{1-ms}{m-2}\sim\f{2\v}{m-2}\leq 2\v.
$$
Then \eqref{Z0} yields the following inequality for suitably large
$m$,
\begin{equation}\label{Z}
\f12\f{d}{dt}(Z^2(t))+\f12\varphi(t)^2\leq
Cm(\f{m}{\v})^{\f{2}{\b-1}}(1+Z(t)^2)^{2+3\v}.
\end{equation}
Note that
\begin{equation}\label{Z1}
\begin{array}{ll}
Z^2(t)&\di =\int(\mu\o^2+\f{F^2}{2\mu+\l(\r)}) dx\\
&\di \leq C\int[\mu\o^2+(2\mu+\l(\r))({\rm
div}u)^2+\f{P^2(\r)}{2\mu+\l(\r)})] dx\\
&\di\leq C\big(\phi(t)+\int P^2(\r)dx\big)\in L^1(0,T),
\end{array}
\end{equation}
where $\phi(t)$ is defined as in \eqref{July-19}.

 Applying the Gronwall's inequality to \eqref{Z} and using
\eqref{Z1} show that
\begin{equation*}
\f{1}{(1+Z^2(t))^{3\v}}-\f{1}{(1+Z^2(0))^{3\v}}+Cm\v(\f{m}{\v})^{\f{2}{\b-1}}\geq0.
\end{equation*}
Then we have the inequality
\begin{equation}\label{i1}
\f{1}{(1+Z^2(t))^{3\v}}\geq\f{1}{2(1+Z^2(0))^{3\v}},
\end{equation}
provided that
\begin{equation}\label{condition}
Cm\v(\f{m}{\v})^{\f{2}{\b-1}}\leq \f{1}{2(1+Z^2(0))^{3\v}}.
\end{equation}
This condition, i. e., \eqref{condition}, is satisfied if
\begin{equation}\label{condition-1}
Cm^{1+\f{2}{\b-1}}2^{-m(1-\f{2}{\b-1})}\leq \f12,
\end{equation}
since
$$
\begin{array}{ll}
\di Z^2(0)&\di=\int\Big[\mu(\o_0)^2+\f{(F_0)^2}{2\mu+\l(\r_0)}\Big] dx\\
&\di \leq C\Big[\|\nabla
^2u_0\|^2_{2}+\|\r_0\|_{W^{2,q}(\mathbb{R}^2)}^{2\b}\|\nabla
^2u_0\|^2_{2}+\|\r_0\|_{W^{2,q}(\mathbb{R}^2)}^{2\g}\Big]\leq C.
\end{array}
$$
Now if $\b>3,$ that is, $1-\f{2}{\b-1}>0$, then we can choose
sufficiently large $m>2$ to guarantee the condition
\eqref{condition-1}. Consequently, the inequality \eqref{i1} is
satisfied with  $\b>3$ and sufficiently large $m>2$. Then
\begin{equation*}\label{Z2}
Z^2(t)\leq 2^{2^{m-1}}(1+Z^2(0))-1\leq C,
\end{equation*}
and \begin{equation*}\label{Z3}
 \int_0^T\varphi(t)dt\leq C.
\end{equation*}
Thus the proof of Lemma \ref{lemma-u-der} is completed. $\hfill\Box$

\underline{Step 5: Weighted estimates for the density:}

The following Lemma \ref{density-wee} is used in estimating
the nonlinear terms \eqref{I9-1} and \eqref{I9-2}.
\begin{Lemma}\label{density-wee}
It holds that for $\a_1>\a$ with $\a$ being the weight in Lemma \ref{lemma-wee}
$$
\int\r |x|^{\a_1} dx\leq C.
$$
\end{Lemma}
\noindent{\bf Proof:} Multiplying the continuity equation
$\eqref{CNS}_1$ by $|x|^{\a_1}$ yields that
$$
(\r |x|^{\a_1})_t+{\rm div}(\r u|x|^{\a_1})-\r u\cdot\nabla(|x|^{\a_1})=0.
$$
Integrating the above equation with respect to $t,x$ over $[0,t]\times\mathbb{R}^2$ gives that
$$
\begin{array}{ll}
\di \int\r(t,x) |x|^{\a_1} dx&\di=\int\r_0(x) |x|^{\a_1} dx+\int_0^t\int \r u\cdot\nabla(|x|^{\a_1})dxd\tau\\
&\di \leq M+\a_1\int_0^T\int \r|u||x|^{\a_1-1}dxdt\\
&\di \leq M+\a_1\int_0^T\|u|x|^{\a_1-1}\|_p\|\r\|_{p_1}dt\\
&\di \leq M+C\int_0^T\||x|^{\f\a 2}\nabla u\|_2dt\leq C(M),
\end{array}
$$
where in the second and third inequalities we have used the H${\rm\ddot{o}}$lder inequality and Caffarelli-Kohn-Nirenberg inequality Lemma \ref{lemma2} (1), respectively, such that the positive constants $p, p_1, \a_1, \a$ satisfying the relations
$$
\f1p+\f1{p_1}=1,\qquad~\f{1}{p}+\f{\a_1-1}{2}=\f12+\f{\f\a 2-1}{2}=\f{\a}{4},
$$
which implies that
$$
\a_1=\f{\a}{2}+1-\f{2}{p}>\a,~~{\rm if}~{\rm we}~{\rm choose}~p>\f{4}{2-\a}.
$$
Thus the proof of Lemma \ref{density-wee} is completed. $\hfill\Box$

\underline{Step 6: Second order derivative estimates for the
velocity:}

\begin{Lemma}\label{lemma-u-sec}
  There exists a positive constant $C=C(M)$, such that
  \begin{equation*}
\sup_{t\in[0,T]}\|(1+|x|^{\f\a 2})\sqrt\r(H,L)\|_2^2(t)+\int_0^T\|(1+|x|^{\f\a2})\nabla(H,L)\|_2^2dt\\
 \leq C.
  \end{equation*}
\end{Lemma}
{\bf Proof:} Multiplying the equations $\eqref{H-L}_1$ and
$\eqref{H-L}_2$ by $H$ and $L$, respectively,  summing the resulted
equations together and then integrating  with respect to $x$ over
$\mathbb{R}^2$ yields that
\begin{equation}\label{HL1}
\begin{array}{ll}
\di
\f12\f{d}{dt}\int\r(H^2+L^2)dx+\int\mu(H_{x_1}-L_{x_2})^2+(2\mu+\l(\r))(H_{x_2}+L_{x_1})^2dx
\\
\di =\int\r(H^2+L^2){\rm div}u dx-\int\mu\o{\rm div}u (L_{x_2}-H_{x_1}) dx\\
\di~~ -\int\r(2\mu+\l(\r))\big[F(\f{1}{2\mu+\l(\r)})^\prime+(\f{P(\r)}{2\mu+\l(\r)})^\prime\big]({\rm div}u)(H_{x_2}+L_{x_1}) dx\\
\di~~-\int \big[H(u_{x_2}\cdot\nabla F+\mu
u_{x_1}\cdot\nabla\o)+L(u_{x_1}\cdot\nabla F-\mu u_{x_2}\cdot\nabla\o)\big]dx\\
\di~~+\int(2\mu+\l(\r))[(u_{1x_1})^2+2u_{1x_2}u_{2x_1}+(u_{2x_2})^2](H_{x_2}+L_{x_1})
dx.
\end{array}
\end{equation}
Set
\begin{equation*}
Y(t)=\left(\int\r(H^2+L^2)dx\right)^{\f12},
\end{equation*}
and
\begin{equation*}
\psi(t)=\left(\int\mu(H_{x_1}-L_{x_2})^2+(2\mu+\l(\r))(H_{x_2}+L_{x_1})^2dx\right)^{\f12}.
\end{equation*}
Note that
$$
\begin{array}{ll}
\di \int(|\nabla H|^2+|\nabla L|^2)dx&\di =\int
(H_{x_1}^2+H_{x_2}^2+L_{x_1}^2+L_{x_2}^2)dx\\
&\di =\int\big[(H_{x_1}-L_{x_2})^2+(H_{x_2}+L_{x_1})^2\big] dx\leq \f{1}{\mu}\psi^2(t).
\end{array}
$$
Then it follows from the elliptic system
\begin{equation}\label{HL-elliptic}
\mu\o_{x_1}+F_{x_2}=\r H,\qquad\qquad -\mu\o_{x_2}+F_{x_1}=\r L,
\end{equation}
 that
\begin{equation}\label{fact4}
\|\nabla(F,\o)\|_p\leq C\|\r(H,L)\|_p, \qquad\forall 1<p<+\i.
\end{equation}
Now we estimate the right hand side of
\eqref{HL1} term by term. First, by the H${\rm\ddot{o}} $lder
inequality, \eqref{fact4}, Caffarelli-Kohn-Nirenberg inequality in Lemma \ref{lemma2} (1) and the density estimate
\eqref{density-e},  it holds that
\begin{equation}\label{HL2}
\begin{array}{ll}
\di |\int\r(H^2+L^2){\rm div} u dx| =|\int\r(H^2+L^2) \f{F+P(\r)}{2\mu+\l(\r)}dx|\\
\di \leq\|\sqrt\r(H,L)\|_2\|(H,L)\|_{\f{4}{\a}}\|\f{\sqrt\r(F+P(\r))}{2\mu+\l(\r)}\|_{\f{4}{2-\a}} \leq C\|\sqrt\r(H,L)\|_2\|(H,L)\|_{\f{4}{\a}}(\|F\|_{\f{4}{2-\a}}+1)\\
\di \leq C\|\sqrt\r(H,L)\|_2\|(H,L)\|_{\f{4}{\a}}(\|\nabla F\|_{\f{4}{4-\a}}+1)\leq C\|\sqrt\r(H,L)\|_2\|(H,L)\|_{\f{4}{\a}}(\|\r(H,L)\|_{\f{4}{4-\a}}+1)\\
\di\leq C\|\sqrt\r(H,L)\|_2\||x|^{\f\a2}\nabla(L,H)\|_2(\|\sqrt\r(H,L)\|_2\|\r\|^{\f12}_{\f{2}{2-\a}}+1)\\
  \di  \leq \s_1\||x|^{\f\a2}\nabla(L,H)\|_2^2+C_{\s_1}(Y(t)^4+1)
\end{array}
\end{equation}
where $\s_1>0$ is a small constant to be determined and $C_{\s_1}$ is a positive constant depending on $\s_1$.
Second, direct estimates give
\begin{equation}\label{HL22}
\begin{array}{ll}
\di |-\int\mu\o{\rm div}u (L_{x_2}-H_{x_1}) dx| \leq \mu\left(\int
(L_{x_2}-H_{x_1})^2dx\right)^{\f12}\left(\int\o^2({\rm div}
u)^2dx\right)^{\f12}\\[2mm]
\di \leq \s\psi^2(t)+C_\s \int\o^2({\rm div} u)^2dx\leq \s\psi^2(t)+C_\s \|\o\|_4^2\|\f{F+P(\r)}{2\mu+\l(\r)}\|_4^2\\
\di \leq \s\psi^2(t)+C_\s \|\o\|_4^2(1+\|F\|_4^2) \leq \s \psi^2(t)+C_\s \|\nabla\o\|_{\f43}^2(\|\nabla F\|_{\f43}^2+1)\\
\di \leq \s \psi^2(t)+C_\s (\|\r(H,L)\|_{\f43}^4+1)\leq \s \psi^2(t)+C_\s (\|\sqrt\r(H,L)\|_2^4\|\r\|_{8}^2+1)\\
\di \leq  \s \psi^2(t)+C_\s(Y(t)^4+1).
\end{array}
\end{equation}
Similarly, one has
\begin{equation}\label{HL23}
\begin{array}{ll}
\di
|-\int\r(2\mu+\l(\r))\big[F(\f{1}{2\mu+\l(\r)})^\prime+(\f{P(\r)}{2\mu+\l(\r)})^\prime\big]{\rm
div}u(H_{x_2}+L_{x_1}) dx|\\
\di \leq \s\int
(2\mu+\l(\r))(H_{x_2}+L_{x_1})^2dx\\
\di \qquad\qquad+C_\s\int\r^2(2\mu+\l(\r))\big[F(\f{1}{2\mu+\l(\r)})^\prime+(\f{P(\r)}{2\mu+\l(\r)})^\prime\big]^2({\rm
div}u)^2 dx\\
\di \leq
\s\psi^2(t)+C_\s\int\r^2\big[F(\f{1}{2\mu+\l(\r)})^\prime+(\f{P(\r)}{2\mu+\l(\r)})^\prime\big]^2\f{|F|^2+P^2(\r)}{2\mu+\l(\r)} dx\\
\di \leq \s\psi^2(t)+C_\s (\|F\|_4^4+1)\leq \s\psi^2(t)+C_\s (\|\nabla F\|_{\f43}^4+1)\\
\di \leq \s\psi^2(t)+C_\s (\|\r(H,L)\|_{\f43}^4+1)\leq \s \psi^2(t)+C_\s (Y(t)^4+1).
\end{array}
\end{equation}
Next,
\begin{equation}\label{HL24}
\begin{array}{ll}
\di |-\int \big[H(u_{x_2}\cdot\nabla F+\mu
u_{x_1}\cdot\nabla\o)+L(u_{x_1}\cdot\nabla F-\mu u_{x_2}\cdot\nabla\o)\big]dx|\\
\di\leq C\int|(H,L)||\nabla u||\nabla(F,\o)| dx\leq C\|\nabla u\|_2\|(H,L)\|_{\f4\a}\|\nabla(F,\o)\|_{\f{4}{2-\a}}\\
\di\leq C\||x|^{\f\a2}\nabla(H,L)\|_2\|\r(H,L)\|_{\f{4}{2-\a}},
\end{array}
\end{equation}
where one has used the fact that
$$
\|\nabla u\|_2\leq C(\|{\rm div}u\|_2+\|\o\|_2)\leq
C(\|\f{F+P(\r)}{2\mu+\l(\r)}\|_2+\|\o\|_2)\leq C.
$$
Note that
\begin{equation}\label{HL25}
\begin{array}{ll}
\di \|\r(H,L)\|_{\f{4}{2-\a}}&\di=\left(\int\r^{\f{4}{2-\a}}|(H,L)|^{\f{4}{2-\a}}dx\right)^{\f{2-\a}{4}} =\left(\int\sqrt\r|(H,L)||(H,L)|^{\f{2+\a}{2-\a}}\r^{\f{6+\a}{2(2-\a)}}dx\right)^{\f{2-\a}{4}} \\
&\di\leq \|\sqrt\r(H,L)\|_2^{\f{2-\a}{4}}
\|(H,L)\|_{\f{p(2+\a)}{2-\a}}^{\f{2+\a}{4}}
\|\r\|_{\f{q(6+\a)}{2(2-\a)}}^{\f{6+\a}{8}}\\
&\di \leq C  \|\sqrt\r(H,L)\|_2^{\f{2-\a}{4}}
\Big[\|(H,L)\|_{\f4\a}^{\f{2+\a}{4}}+\|\nabla(H,L)\|_2^{\f{2+\a}{4}}\Big]\\
&\di \leq C  \|\sqrt\r(H,L)\|_2^{\f{2-\a}{4}}
\|(1+|x|^{\f\a2})\nabla(H,L)\|_2^{\f{2+\a}{4}},
\end{array}
\end{equation}
where $p,q>2$ satisfying $\f1p+\f1q=1$ and we have chosen $p>2$ such
that $\f{p(2+\a)}{2-\a}\geq \f4\a$ for $\a>0$ given in Lemma
\ref{lemma-wee}.

It follows from \eqref{HL24} and \eqref{HL25} that
\begin{equation}\label{HL26}
\begin{array}{ll}
\di |-\int \big[H(u_{x_2}\cdot\nabla F+\mu
u_{x_1}\cdot\nabla\o)+L(u_{x_1}\cdot\nabla F-\mu u_{x_2}\cdot\nabla\o)\big]dx|\\
\di\leq
CY(t)^{\f{2-\a}{4}}\|(1+|x|^{\f\a2})\nabla(H,L)\|_2^{\f{6+\a}{4}}\leq
\s\|\nabla(H,L)\|_2^2+\s_1
\||x|^{\f\a2}\nabla(H,L)\|_2^2+C_{\s,\s_1} Y(t)^2.
\end{array}
\end{equation}
Moreover,
\begin{equation}\label{HL27}
\begin{array}{ll}
\di|\int(2\mu+\l(\r))[(u_{1x_1})^2+2u_{1x_2}u_{2x_1}+(u_{2x_2})^2](H_{x_2}+L_{x_1})
dx|\\
\di\leq
\s\int (2\mu+\l(\r))(H_{x_2}+L_{x_1})^2dx+C_\s\int(2\mu+\l(\r))|\nabla u|^4dx\\
\di\leq
\s\psi(t)^2+C_\s\Big[\|\nabla u\|_4^4+\|\l(\r)\|_2^4\|\nabla u\|_8^4\Big]
\di\leq
\s\psi(t)^2+C_\s\Big[\|({\rm div} u,\o)\|_4^4+\|({\rm div} u,\o)\|_8^4\Big]\\
\di\leq
\s\psi(t)^2+C_\s\Big[\|(F,\o)\|_4^4+\|(F,\o)\|_8^4+1\Big]
\di\leq
\s\psi(t)^2+C_\s\Big[\|\nabla(F,\o)\|_{\f43}^4+\|\nabla(F,\o)\|_{\f85}^4+1\Big]\\
\di \leq \s\psi(t)^2+C_\s\Big[\|\r(H,L)\|_{\f43}^4+\|\r(H,L)\|_{\f85}^4+1\Big]
\di  \leq \s \psi(t)^2+C_\s (Y(t)^4+1).
\end{array}
\end{equation}
Substituting the estimates \eqref{HL2}-\eqref{HL23}, \eqref{HL26}
and \eqref{HL27} into \eqref{HL1}, one can arrive at
\begin{equation*}
\f12\f{d}{dt}(Y^2(t))+\psi^2(t)\leq
4\s\psi^2(t)+2\s_1\||x|^{\f\a2}\nabla(H,L)\|_2^2+C_{\s,\s_1}(1+Y^2(t))^2.
\end{equation*}
Choosing $4\s=\f12$, noting that $Y^2(t)=\varphi^2(t)\in L^1(0,T)$,
and then using Gronwall's inequality yield that
\begin{equation}\label{Y-e}
Y^2(t)+\int_0^t\psi^2(t)dt\leq Y^2(0)+2\s_1\int_0^t\||x|^{\f\a2}\nabla(H,L)\|_2^2d\tau+C_{\s_1}.
\end{equation}
By the compatibility
condition \eqref{cc}, one has
\begin{equation*}\label{id1}
Y^2(0)=\|\sqrt{\r_0}(H_0,L_0)\|_2^2=\|g\|_2^2\leq C.
\end{equation*}
This, together with \eqref{Y-e}, shows that
\begin{equation}\label{HL-11}
Y^2(t)+\int_0^t \psi^2(\tau)d\tau\leq 2\s_1\int_0^t\||x|^{\f\a2}\nabla(H,L)\|_2^2d\tau+C_{\s_1}.
\end{equation}
In order to close the estimates in Lemma \ref{lemma-u-sec}, we need
to carry out the following weighted estimates to  $\sqrt\r(H,L)$.
Note that
$$
\begin{array}{ll}
\di \mu(H_{x_1}-L_{x_2})^2+(2\mu+\l(\r))(H_{x_2}+L_{x_1})^2-{\rm div}\big\{(2\mu+\l(\r))(H_{x_2}+L_{x_1})(L,H)^t\big\}\\
\di-\big[\mu H(H_{x_1}-L_{x_2})\big]_{x_1}+\big[\mu L(H_{x_1}-L_{x_2})\big]_{x_2}\\
\di =\mu|\nabla(H,L)|^2+(\mu+\l(\r))(H_{x_2}+L_{x_1})^2-{\rm div}\big[\mu\nabla\big( \f{H^2+L^2}{2}\big)+(\mu+\l(\r))(H_{x_2}+L_{x_1})(L,H)^t\big].
\end{array}
$$
Similar to \eqref{HL1}, it follows from \eqref{H-L} that
\begin{equation}\label{HL-e1}
\begin{array}{ll}
\di \f{d}{dt}\int \f12 \r(H^2+L^2)|x|^{\a}dx+J(t)=\int\f12 (H^2+L^2)\r u\cdot\nabla(|x|^{\a})dx+\int \r (H^2+L^2)({\rm div}u)|x|^{\a}\big]dx\\
\di+\int\mu \o ({\rm div}u)\big[ H(|x|^{\a})_{x_1}- L(|x|^{\a} )_{x_2}\big]dx-\int\mu\o({\rm div}u) (L_{x_2}-H_{x_1})|x|^{\a}dx\\
\di  -\int \big[H(u_{x_2}\cdot\nabla F+\mu
u_{x_1}\cdot\nabla\o)+L(u_{x_1}\cdot\nabla F-\mu u_{x_2}\cdot\nabla\o)\big]|x|^{\a}dx\\
\di -\int\r(2\mu+\l(\r))\big[F(\f{1}{2\mu+\l(\r)})^\prime+(\f{P(\r)}{2\mu+\l(\r)})^\prime\big]({\rm div}u)(L,H)^t\cdot\nabla(|x|^{\a})dx\\
\di -\int\r(2\mu+\l(\r))\big[F(\f{1}{2\mu+\l(\r)})^\prime+(\f{P(\r)}{2\mu+\l(\r)})^\prime\big]({\rm div}u)(H_{x_2}+L_{x_1})|x|^{\a}dx\\
\di +\int (2\mu+\l(\r))[(u_{1x_1})^2+2u_{1x_2}u_{2x_1}+(u_{2x_2})^2](L,H)^t\cdot\nabla(|x|^{\a})dx\\
\di +\int (2\mu+\l(\r))[(u_{1x_1})^2+2u_{1x_2}u_{2x_1}+(u_{2x_2})^2](H_{x_2}+L_{x_1})|x|^{\a}dx:=\sum_{i=1}^9{I_i},
\end{array}
\end{equation}
where
\begin{equation}\label{J(t)}
\begin{array}{ll}
J(t)=&\di \int \Big\{|x|^{\a}\big[\mu|\nabla(H,L)|^2+(\mu+\l(\r))(H_{x_2}+L_{x_1})^2\big]\\
&\di\qquad -\big[\mu\nabla\big(\f{H^2+L^2}{2}\big)+(\mu+\l(\r))(H_{x_2}+L_{x_1})(L,H)^t\big]\cdot\nabla(|x|^{\a})\Big\}dx.
\end{array}
\end{equation}
First, $J(t)$ in \eqref{J(t)} can be estimated similarly as in \eqref{WE3}, \eqref{WE4}, \eqref{WE5} and \eqref{quad} if replacing $u$ by $(L,H)^t=\dot u$. Therefore, if $\a^2<4(\sqrt2-1)$, then there exists a positive constant $C_\a$ such that
\begin{equation}\label{JE}
J(t)\geq C^{-1}_\a\Big[\||x|^{\f\a2}\nabla (H,L)\|_2^2+\||x|^{\f\a2}(H_{x_2}+L_{x_1})\|_2^2+\|\sqrt{\l(\r)}|x|^{\f\a2}(H_{x_2}+L_{x_1})\|_2^2\Big]-C_\a\|\nabla(H,L)\|_2^2.
\end{equation}
Then the terms $I_i~(i=1,2,\cdots 9)$ on the right hand side of \eqref{HL-e1} will be estimated as follows. By the H${\rm \ddot{o}}$lder inequality, Caffarelli-Kohn-Nirenberg inequality in Lemma \ref{lemma2} (1) and Young inequality, it holds that
\begin{equation}\label{I1-e}
\begin{array}{ll}
\di |I_1|&\di \leq \f{\a}{2}|\int\r |(H,L)|^2|u||x|^{\a-1} dx|\leq \f{\a}{2}\|\sqrt\r (H,L)\|_2\|\sqrt\r\|_p\||x|^{\b_1}(H,L)\|_q\||x|^{\a-\b_1-1}u\|_r\\
\di &\di\leq C\|\sqrt\r (H,L)\|_2\|\nabla(H,L)\|_2^\t\||x|^{\f\a2}\nabla(H,L)\|_2^{1-\t}\||x|^{\f\a2}\nabla u\|_2\\
\di &\di\leq C\|\sqrt\r (H,L)\|_2\|(1+|x|^{\f\a2})\nabla(H,L)\|_2\||x|^{\f\a2}\nabla u\|_2\\
\di &\di\leq \s\|(1+|x|^{\f\a2})\nabla(H,L)\|_2^2+C_\s\|\sqrt\r (H,L)\|_2^2\||x|^{\f\a2}\nabla u\|^2_2
\end{array}
\end{equation}
where the positive constants in the above inequality \eqref{I1-e} satisfying that
$$
\f1p+\f1q+\f1r=\f12,\qquad
\f1q+\f{\b_1}{2}=\f{\a}{4}(1-\t), \qquad \f{1}{r}+\f{\a-\b_1-1}{2}=\f{\a}{4},
$$
which implies that
$$
p=\f{4}{\a\t}, ~~~{\rm with}~~\t\in(0,1).
$$
Then H${\rm \ddot{o}}$lder inequality gives that
\begin{equation}\label{I2}
\di |I_2|\leq \|\sqrt\r (H,L)|x|^{\f\a2}\|_2\|\sqrt\r|x|^{\b_1}\|_{2\g}\|{\rm div} u\|_p\||x|^{\f\a2-\b_1}(H,L)\|_q,
\end{equation}
where $\b_1>0$ is to be determined and $p,q>2$ satisfying that
$$
\f{1}{2\g}+\f{1}{p}+\f1q=\f12.
$$
Note that it follows from Lemma \ref{lemma-wee} and Lemma \ref{lemma-rho} that if $2\g\b_1=\a$, then
$$
\|\sqrt\r|x|^{\b_1}\|_{2\g}=\big(\int\r^\g|x|^{2\g\b_1} dx\big)^{\f{1}{2\g}}\leq C.
$$
Thus, if $\b_1=\f\a{2\g},$ and
$$
\f{1}{q}+\f{\f{\a}{2}-\f{\a}{2\g}}{2}=\f{\a}{4}(1-\t),~~{\rm with}~~\t \in(0,1).
$$
then by Caffarelli-Kohn-Nirenberg inequality and Sobolev inequality, it follows from \eqref{I2} that
\begin{equation}\label{I2-e}
\begin{array}{ll}
\di |I_2|&\di \leq C \|\sqrt\r (H,L)|x|^{\f\a2}\|_2\|\f{F+P(\r)}{2\mu+\l(\r)}\|_p\||x|^{\f\a2-\f{\a}{2\g}}(H,L)\|_q\\
&\di \leq C\|\sqrt\r (H,L)|x|^{\f\a2}\|_2\big(\|F\|_p+\|P(\r)\|_p\big)\|\nabla(H,L)\|_2^\t\||x|^{\f\a2}\nabla(H,L)\|_2^{1-\t}\\
&\di \leq C\|\sqrt\r (H,L)|x|^{\f\a2}\|_2\big(\|\nabla F\|_{\f{2p}{p+2}}+1\big)\|(1+|x|^{\f\a2})\nabla(H,L)\|_2\\
&\di \leq C\|\sqrt\r (H,L)|x|^{\f\a2}\|_2\big(\|\r(H,L)\|_{\f{2p}{p+2}}+1\big)\|(1+|x|^{\f\a2})\nabla(H,L)\|_2\\
&\di \leq C\|\sqrt\r (H,L)|x|^{\f\a2}\|_2\big(\|\sqrt\r(H,L)\|_{2}+1\big)\|(1+|x|^{\f\a2})\nabla(H,L)\|_2\\
&\di \leq \s\|(1+|x|^{\f\a2})\nabla(H,L)\|_2^2+ C_\s\|\sqrt\r (H,L)|x|^{\f\a2}\|^2_2\big(\|\sqrt\r(H,L)\|^2_{2}+1\big).
\end{array}
\end{equation}
Similarly, it follows that
\begin{equation}\label{I3}
\begin{array}{ll}
\di |I_3|&\di \leq C\int|\o||(H,L)||{\rm div} u||x|^{\a-1} dx\leq C\|{\rm div} u\|_2\||x|^{\b_1}\o\|_p\||x|^{\a-\b_1-1}(H,L)\|_q\\
&\di \leq C\|\nabla \o\|_{p_1}\|\nabla(H,L)\|_2^\t\||x|^{\f\a2}\nabla(H,L)\|_2^{1-\t},
\end{array}
\end{equation}
if one has $\t \in(0,1), p, p_1, q>2, \b_1>0$ and
$$
\f1p+\f1q=\f12,\qquad \f{1}{p}+\f{\b_1}{2}=\f{1}{p_1}-\f12,\qquad \f1q+\f{\a-\b_1-1}{2}=\f{\a}4(1-\t),
$$
which implies that
$$
p_1=\f{4}{2+\a(1+\t)}<2.
$$
Thus one can obtain from \eqref{I3} that
\begin{equation}\label{I3-e}
\begin{array}{ll}
\di |I_3|&\di \leq C\|\r(H,L) \|_{p_1}\|(1+|x|^{\f\a2})\nabla(H,L)\|_2\leq C\|\sqrt\r(H,L) \|_{2}\|(1+|x|^{\f\a2})\nabla(H,L)\|_2\\
&\di \leq \s\|(1+|x|^{\f\a2})\nabla(H,L)\|_2^2+C_\s\|\sqrt\r(H,L) \|_{2}^2.
\end{array}
\end{equation}
Then for $\b_1>0$ to be determined and $\f1p+\f1q=\f12$, it holds that
\begin{equation}\label{I4}
\begin{array}{ll}
\di |I_4|&\di \leq C\||x|^{\f\a2}\nabla(H,L)\|_2\||x|^{\b_1}{\rm div} u\|_p\||x|^{\f{\a}{2}-\b_1}\o\|_q\\
&\di \leq C\||x|^{\f\a2}\nabla(H,L)\|_2\Big[\||x|^{\b_1}F\|_p+\||x|^{\b_1}P(\r)\|_p\Big]\||x|^{\f{\a}{2}-\b_1}\o\|_q\\
&\di \leq C\||x|^{\f\a2}\nabla(H,L)\|_2\Big[\|\nabla F\|_{p_1}+1\Big]\|\nabla \o\|_{q_1},
\end{array}
\end{equation}
where we have chosen $\b_1=\f{\a}{p\g}$ and  by Caffarelli-Kohn-Nirenberg inequality, $p_1$ and $ q_1$ satisfy $1<p_i<2~(i=1,2)$ and
$$
\f{1}{p}+\f{\b_1}{2}=\f{1}{p_1}-\f12,\qquad \f{1}{q}+\f{\f{\a}{2}-\b}{2}=\f{1}{q_1}-\f12.
$$
Thus one can get from \eqref{I4} that
\begin{equation}\label{I4-e}
\begin{array}{ll}
\di |I_4|&\di \leq C\||x|^{\f\a2}\nabla(H,L)\|_2\Big[\|\r(H,L) \|_{p_1}+1\Big]\|\r(H,L) \|_{q_1}\\
&\di \leq C\||x|^{\f\a2}\nabla(H,L)\|_2\Big[\|\sqrt\r(H,L) \|_{2}+1\Big]\|\sqrt\r(H,L) \|_{2}\\
&\di \leq \s\||x|^{\f\a2}\nabla(H,L)\|_2^2+C_\s\Big[\|\sqrt\r(H,L) \|_{2}^4+1\Big].
\end{array}
\end{equation}
Now we estimate $I_5,$ which is a little more delicate. For $\b_1>0$ to be determined and $\f1p+\f1q=\f12,$ it holds that
\begin{equation}\label{I5}
\begin{array}{ll}
\di I_5= -\int \big[H(u_{x_2}\cdot\nabla F+\mu
u_{x_1}\cdot\nabla\o)+L(u_{x_1}\cdot\nabla F-\mu u_{x_2}\cdot\nabla\o)\big]|x|^{\a}dx\\
\di\qquad \leq C\int||x|^{\a-\b_1}(H,L)||\nabla u|||x|^{\b_1}\nabla (F,\o)| dx\\
 \di\qquad \leq C\|\nabla u\|_2\||x|^{\a-\b_1}(H,L)\|_p\||x|^{\b_1}\nabla (F,\o)\|_q\\
\qquad\di \leq C\||x|^{\a-\b_1}(H,L)\|_p \|\r(H,L)|x|^{\b_1}\|_q,
\end{array}
\end{equation}
where in the last inequality one has used the equalities \eqref{HL-elliptic} and  Lemma \ref{lemma3} (2) provided $\b_1$ satisfies
\begin{equation}\label{beta1}
0<\b_1<2(1-\f1q).
\end{equation}
Now if we also choose
\begin{equation}\label{beta1-2}
\f\a2<\b_1<\a,
\end{equation}
 then
\begin{equation}\label{I511}
\||x|^{\a-\b_1}(H,L)\|_p\leq
C\|\nabla(H,L)\|_2^{\t_1}\||x|^{\f\a2}\nabla(H,L)\|_2^{1-\t_1}\leq C
\|(1+|x|^{\f\a2})\nabla(H,L)\|_2,
\end{equation}
where
\begin{equation}\label{p-r}
\f{1}{p}+\f{\a-\b_1}{2}=\f{\a}{4}(1-\t_1),~~{\rm that}~{\rm is},~~ p=\f{4}{2\b_1-\a(1+\t_1)}.
\end{equation}
For $\b_2>0$ to be determined and for $\f1{p_1}+\f1{q_1}=\f12$, it holds that
$$
\begin{array}{ll}
\di \|\r(H,L)|x|^{\b_1}\|_q^q=\int\r^q|(H,L)|^q|x|^{\b_1q} dx\\
\di
\qquad\leq \|\sqrt\r(H,L)|x|^{\f\a2}\|_2\||(H,L)|^{q-1}|x|^{\b_1q-\f{\a}{2}-\b_2}\|_{p_1}\|\r^{q-\f12}|x|^{\b_2}\|_{q_1}.
\end{array}$$
If we choose $\b_2=\f{\a}{q_1\g}$, then it holds that
$$
\begin{array}{ll}
\di
\|\r^{q-\f12}|x|^{\b_2}\|_{q_1}=(\int\r^{(q-\f12)q_1}|x|^{\b_2q_1}
dx)^{\f{1}{q_1}}\\
\qquad\qquad \quad \di =(\int\r^{(q-\f12)q_1}|x|^{\f{\a}{\g}}
dx)^{\f{1}{q_1}}\leq
C\|\r|x|^{\f{\a}{\g}}\|_\g^{\f{1}{q_1}}\|\r^{(q-\f12)q_1-1}\|_{\f{\g}{\g-1}}^{\f{1}{q_1}}\leq
C,
\end{array}$$
and thus
\begin{equation}\label{I512}
\begin{array}{ll}
\di \|\r(H,L)|x|^{\b_1}\|_q^q&\di \leq C\|\sqrt\r(H,L)|x|^{\f\a2}\|_2\||(H,L)|^{q-1}|x|^{\b_1q-\f{\a}{2}-\b_2}\|_{p_1}\\
&\di \leq C\|\sqrt\r(H,L)|x|^{\f\a2}\|_2\||(H,L)||x|^{\f{\b_1q-\f{\a}{2}-\b_2}{q-1}}\|^{q-1}_{p_1(q-1)}\\
&\di \leq C\|\sqrt\r(H,L)|x|^{\f\a2}\|_2\|\nabla(H,L)\|_{2}^{\t_2(q-1)}\||x|^{\f\a2}\nabla(H,L)\|_{2}^{(1-\t_2)(q-1)}\\
&\di \leq C\|\sqrt\r(H,L)|x|^{\f\a2}\|_2\|(1+|x|^{\f\a2})\nabla(H,L)\|_{2}^{q-1}.
\end{array}
\end{equation}
where
\begin{equation}\label{q-r}
\f{1}{p_1(q-1)}+\f{\b_1q-\f{\a}{2}-\b_2}{2(q-1)}=\f{\a}{4}(1-\t_2)
\end{equation}
It follows from \eqref{p-r} and \eqref{q-r} that
$$
q_1=(1+\f{\a}{2\g})\Big[\f{q-1}{2}+\f{\a}{4}\big(\t_1 q+\t_2(q-1)\big)\Big]^{-1}>2, ~~{\rm if}~q\rightarrow2+, \t_i\rightarrow0+,~(i=1,2).
$$
Note that $p$ is sufficiently large when $q\rightarrow 2+$ and thus the above restrictions \eqref{beta1} and \eqref{beta1-2} on $\b_1$ when estimating $I_5$ could be satisfied.
Substituting \eqref{I511} and \eqref{I512}  into \eqref{I5} yields that
\begin{equation}\label{I5-e}
\begin{array}{ll}
\di |I_{5}|\leq  C\|(1+|x|^{\f\a2})\nabla(H,L)\|_2\Big[\|\sqrt\r(H,L)|x|^{\f\a2}\|^{\f1q}_2\|(1+|x|^{\f\a2})\nabla(H,L)\|_2^{1-\f1q}
+\|\sqrt\r(H,L)\|_2\Big]\\
\qquad\di \leq \s\|(1+|x|^{\f\a2})\nabla(H,L)\|_2^2+C_\s\|\sqrt\r(H,L)(1+|x|^{\f\a2})\|^2_2.
\end{array}
\end{equation}
Then, for $\b_1>0$ to be determined and  for $\f1p+\f1q=\f12,$ it holds that
\begin{equation}\label{I6}
\begin{array}{ll}
\di |I_6|\leq C\|{\rm div} u\|_2\|(H,L)|x|^{\b_1}\|_p\Big[\|\r(2\mu+\l(\r))F(\f{1}{2\mu+\l(\r)})^\prime|x|^{\a-\b_1-1}\|_q\\
\di \qquad\qquad\qquad\qquad\qquad\qquad\qquad+\|\r(2\mu+\l(\r))(\f{P(\r)}{2\mu+\l(\r)})^\prime|x|^{\a-\b_1-1}\|_q\Big]\\
\di \leq C\||x|^{\f\a2}\nabla(H,L)\|_2\Big[\|F|x|^{\a-\b_1-1}\|_q+\|P(\r)|x|^{\a-\b_1-1}\|_q\Big],
\end{array}
\end{equation}
where
$$
\f1p+\f{\b_1}{2}=\f{\a}{4}.
$$
Furthermore, for $1<q_1<2$ and
 $$
\f{1}{q}+\f{\a-\b_1-1}{2}=\f{1}{q_1}-\f12,
$$ it follows that
$$
\|F|x|^{\a-\b_1-1}\|_q\leq C\|\nabla F\|_{q_1}\leq C\|\r(H,L)\|_{q_1}\leq C\|\sqrt\r(H,L)\|_2,
$$
and
$$
\|P(\r)|x|^{\a-\b_1-1}\|_q=\big(\int\r^{\g q}|x|^{(\a-\b_1-1)q} dx\big)^{\f1q}\leq C\|\r|x|^{(\a-\b_1-1)q}\|_{\g}\|\r^{\g q-1}\|_{\f{\g}{\g-1}}\leq C,
$$
provided $\b_1$ is chosen such that
$$
(\a-\b_1-1)q\g\leq \a.
$$
Therefore, from \eqref{I6}, it holds that
\begin{equation}\label{I6-e}
\begin{array}{ll}
\di |I_6|\leq  C\||x|^{\f\a2}\nabla(H,L)\|_2\Big[\|\sqrt\r(H,L)\|_2+1\Big]\leq \s \||x|^{\f\a2}\nabla(H,L)\|_2^2+C_\s \Big[\|\sqrt\r(H,L)\|^2_2+1\Big].
\end{array}
\end{equation}
Then
\begin{equation}\label{I7}
\begin{array}{ll}
\di |I_7|\leq C\|(1+\l(\r))(H_{x_2}+L_{x_1})|x|^{\f\a2}\|_2\|(|F|+P(\r))|{\rm div} u||x|^{\f\a2}\|_2\\
\di \leq C\|(1+\l(\r))(H_{x_2}+L_{x_1})|x|^{\f\a2}\|_2\|(|F|+P(\r))||x|^{\f\a2}\|_2\\
\end{array}
\end{equation}

Next, for $p=\f{4}{2-\a}$ and $q=\f{4}{2+\a}$ satisfying $\f1p+\f1q=1,$ one has
\begin{equation}\label{I8}
\begin{array}{ll}
\di |I_8|\leq C\|(H,L)|x|^{\a-1}\|_p\|(2\mu+\l(\r))|\nabla u|^2\|_q\\
\di \leq C\|\nabla(H,L)|x|^{\f\a2}\|_2\Big[\|\nabla u\|^2_{2q}+\|\l(\r)\|_{2q}^2\|\nabla u\|_{4q}^2\Big]\\
\di \leq C\|\nabla(H,L)|x|^{\f\a2}\|_2\Big[\|\nabla u\|_{4q}^2+1\Big]\leq C\|\nabla(H,L)|x|^{\f\a2}\|_2\Big[\|\o\|_{4q}^2+\|{\rm div} u\|_{4q}^2+1\Big]\\
\leq C\|\nabla(H,L)|x|^{\f\a2}\|_2\Big[\|(F,\o)\|_{4q}^2+1\Big]\leq C\|\nabla(H,L)|x|^{\f\a2}\|_2\Big[\|\nabla(F,\o)\|_{\f{4q}{2q+1}}^2+1\Big]\\
\leq C\|\nabla(H,L)|x|^{\f\a2}\|_2\Big[\|\sqrt\r(H,L)\|_2^2+1\Big]\leq \s\|\nabla(H,L)|x|^{\f\a2}\|^2_2+C_\s\Big[\|\sqrt\r(H,L)\|_2^4+1\Big].
\end{array}
\end{equation}
Finally, it holds that
\begin{equation}\label{I9}
\begin{array}{ll}
\di |I_9|\leq \s\int(1+\l(\r))(H_{x_2}+L_{x_1})^2|x|^\a dx+C_\s\int(1+\l(\r))|x|^\a|\nabla u|^4 dx.
\end{array}
\end{equation}
At the same time, it follows from Lemma \ref{lemma3} (2) and Lemma \ref{density-wee} that
\begin{equation}\label{I9-1}
\begin{array}{ll}
\di\int|x|^\a|\nabla u|^4 dx&\di =\||x|^{\f\a 4}\nabla u\|_4^4\leq C\Big[\||x|^{\f\a 4}{\rm div} u\|_4^4+\||x|^{\f\a 4}\o\|_4^4\Big]\\
&\di\leq C\Big[\||x|^{\f\a 4}F\|_4^4+\||x|^{\f\a 4}P(\r)\|_4^4+\||x|^{\f\a 4}\o\|_4^4\Big]\\
&\di \leq C\Big[\|\nabla(F,\o)\|^4_{\f{8}{\a+6}}+\|\r^\f{\a}{\a_1}|x|^\a\|_{\f{\a_1}{\a}}\|\r^{4\g-\f{\a}{\a_1}}\|_{\f{\a_1}{\a_1-\a}}\Big]\\
&\di \leq C\Big[\|\r(H,L)\|^4_{\f{8}{\a+6}}+\|\r|x|^{\a_1}\|^{\f{\a}{\a_1}}_1\Big] \leq C\Big[\|\sqrt\r(H,L)\|^4_2+1\Big],
\end{array}
\end{equation}
and for $p,q>1$ satisfying $\f1p+\f1q+\f{\a}{\a_1}=1,$ one has
\begin{equation}\label{I9-2}
\begin{array}{ll}
\di\int\l(\r)|x|^\a|\nabla u|^4 dx&\di \leq \|\r^\f{\a}{\a_1}|x|^\a\|_{\f{\a_1}{\a}}\|\r^{\b-\f{\a}{\a_1}}\|_p\||\nabla u|^4\|_q = \|\r|x|^{\a_1}\|^{\f{\a}{\a_1}}_1\|\r\|^{\b-\f{\a}{\a_1}}_{p(\b-\f{\a}{\a_1})}\|\nabla u\|^4_{4q}\\
&\di \leq C\|\nabla u\|^4_{4q}\leq C\Big[\|{\rm div} u\|^4_{4q}+\|\o\|^4_{4q}\Big]\leq C\Big[\|F\|^4_{4q}+\|\o\|^4_{4q}+1\Big]\\
&\di\leq C\Big[\|\nabla(F,\o)\|^4_{\f{4q}{2q+1}}+1\Big] \leq C\Big[\|\r(H,L)\|^4_{\f{4q}{2q+1}}+1\Big]\\
 &\di\leq C\Big[\|\sqrt\r(H,L)\|^4_2+1\Big].
\end{array}
\end{equation}
Substituting \eqref{I9-1} and \eqref{I9-2} into \eqref{I9} yields that
\begin{equation}\label{I9-e}
\begin{array}{ll}
\di |I_9|\leq \s\int(1+\l(\r))(H_{x_2}+L_{x_1})^2|x|^\a dx+C_\s\Big[\|\sqrt\r(H,L)\|^4_2+1\Big].
\end{array}
\end{equation}
Substituting the estimates \eqref{JE}, \eqref{I1-e},  \eqref{I2-e}, \eqref{I3-e}, \eqref{I4-e}, \eqref{I5-e}, \eqref{I6-e}, \eqref{I7}, \eqref{I8} and \eqref{I9-e} into \eqref{HL-e1}, then integrating the resulted inequality with respect to $t$ over $[0,t]$, and noting that
$$
\|\sqrt{\r_0}(H_0,L_0)|x|^{\f\a2}\|_2^2=\|g|x|^{\f\a2}\|_2^2,
$$
 it holds that
\begin{equation}\label{g1}
\|\sqrt\r(H,L)|x|^{\f\a2}\|_2^2(t)+\int_0^t\||x|^{\f\a2}\nabla(H,L)\|_2^2dt\leq C\int_0^t\|\nabla(H,L)\|_2^2dt+C,
\end{equation}
which, together with the estimate \eqref{HL-11} and choosing $\s,\s_1$ suitably small, completes the proof of Lemma \ref{lemma-u-sec}. $\hfill\Box$

\underline{Step 7. Upper bound of the density:} We are now ready to derive the upper
bound for the density in the super-norm. First, one has

\begin{Lemma}\label{lemma-F-o}
  It holds that
  \begin{equation*}
\int_0^T\|(F,\o)\|_\i^3dt\leq C(M).
  \end{equation*}
\end{Lemma}
{\bf Proof:} By \eqref{fact4} with $p=3$, one has for $\f1{p_1}+\f1{q_1}=\f12,$
 \begin{equation}\label{Fo1}
\begin{array}{ll}
\di \int_0^T\|\nabla(F,\o)\|_3^3dt&\di \leq C\int_0^T\|\r(H,L)\|_3^3dt =C\int_0^T\int\r^3|(H,L)|^3dxdt\\
&\di =C\int_0^T\int\sqrt\r|(H,L)||(H,L)|^2\r^{\f52}dxdt\\
&\di \leq C\int_0^T\|\sqrt\r(H,L)\|_2\|(H,L)\|_{2p_1}^2\|\r\|^{\f52}_{\f{5q_1}{2}}dt.
\end{array}
\end{equation}
By Caffarelli-Kohn-Nirenberg inequality in Lemma \ref{lemma2} (1), it holds that
\begin{equation}\label{Fo2}
\begin{array}{ll}
\di \|(H,L)\|_{2p}\leq C\|\nabla(H,L)\|_2^\t\||x|^{\b_1}(H,L)\|_{p_1}^{1-\t}\\
\qquad \di \leq C\|\nabla(H,L)\|_2^\t\||x|^{\f\a2}\nabla(H,L)\|_{2}^{1-\t}\leq C\|(1+|x|^{\f\a2})\nabla(H,L)\|_{2},
\end{array}
\end{equation}
with $p>2$ and $\t\in(0,1)$ satisfying
$$
\f{1}{2p}=(1-\t)(\f{1}{p_1}+\f{\b_1}{2})=(1-\t)(\f{1}{2}+\f{\f{\a}{2}-1}{2})=\f{\a(1-\t)}{4}.
$$
Substituting \eqref{Fo2} into \eqref{Fo1} yields that
\begin{equation*}
\di \int_0^T\|\nabla(F,\o)\|_3^3dt \leq C\int_0^T\|(1+|x|^{\f\a2})\nabla(H,L)\|_{2}^2dt\leq C,
\end{equation*}
which, combined with the estimate
\begin{equation*}
\begin{array}{ll}
\di \int_0^T\|(F,\o)\|_3^3dt&\di \leq C\int_0^T\|\nabla(F,\o)\|_{\f65}^3dt\leq C\int_0^T\|\r(H,L)\|_{\f65}^3dt\\
&\di \leq C\int_0^T\|\sqrt\r(H,L)\|_{2}^3\|\sqrt\r\|_3^3dt\\
&\di \leq C\sup_{t\in[0,T]}\|\sqrt\r(H,L))\|_{2}^2\int_0^T(\|\sqrt\r(H,L)\|_{2}^2+1)dt\leq C,
\end{array}
\end{equation*}
yields that
\begin{equation}\label{F-omega-infty}
\int_0^T\|(F,\o)\|_\i^3dt\leq
\int_0^T\|(F,\o)\|_{W^{1,3}(\mathbb{R}^2)}^3dt\leq C.
\end{equation}
The proof of Lemma \ref{lemma-F-o} is finished. $\hfill\Box$

With Lemma \ref{lemma-F-o} in hand, we can obtain the uniform upper bound for the density.

\begin{Lemma}\label{upper-b}
  It holds that
  \begin{equation*}
\r(t,x)\leq C(M),\qquad\forall (t,x)\in [0,T]\times \mathbb{R}^2.
  \end{equation*}
\end{Lemma}
{\bf Proof:} From the continuity equation $\eqref{CNS}_1$, we have
\begin{equation*}
\nu(\r)_t+u\cdot\nabla\nu(\r)+P(\r)+F=0,
\end{equation*}
where $\nu(\r)$ is defined in \eqref{theta}.

Along the particle path $\vec{X}(\tau;t,x)$ through the point
$(t,x)\in[0,T]\times\mathbb{R}^2$ defined by
\begin{equation*}
\left\{
\begin{array}{ll}
\di \f{d\vec{X}(\tau;t,x)}{d\tau}=u(\tau,\vec{X}(\tau;t,x)),\\
 \di
\vec{X}(\tau;t,x)|_{\tau=t}=x,
\end{array}
\right.
\end{equation*}
there holds the following ODE
\begin{equation*}
\f{d}{d\tau}\nu(\r)(\tau,\vec{X}(\tau;t,x))=-P(\r)(\tau,\vec{X}(\tau;t,x))-F(\tau,\vec{X}(\tau;t,x)),
\end{equation*}
which is integrated over $[0,t]$ to yield that
\begin{equation}\label{theta-1}
\nu(\r)(t,x)-\nu(\r_0)(\vec{X}_0)=-\int_0^t(P(\r)+F)(\tau,\vec{X}(\tau;t,x))d\tau,
\end{equation}
with $\vec{X}_0=\vec{X}(\tau;t,x)|_{\tau=0}$.

It follows from \eqref{theta-1} that
\begin{equation*}
2\mu\ln\f{\r(t,x)}{\r_0(\vec{X}_0)}+\f{1}{\b}\r^\b(t,x)+\int_0^tP(\r)(\tau,\vec{X}(\tau;t,x))d\tau=\f1\b\r_0(\vec{X}_0)^{\b}-\int_0^tF(\tau,\vec{X}(\tau;t,x))d\tau.
\end{equation*}
So
\begin{equation*}
2\mu\ln\f{\r(t,x)}{\r_0(\vec{X}_0)}\leq
\f1\b\|\r_0\|_\i^{\b}+\int_0^t\|F(\tau,\cdot)\|_\i d\tau\leq C,
\end{equation*}
which implies that
$$
\f{\r(t,x)}{\r_0(\vec{X}_0)}\leq C.
$$
Therefore, we have
\begin{equation*}
\r(t,x)\leq C, \qquad \forall (t,x)\in[0,T]\times\mathbb{R}^2.
\end{equation*}
Hence  Lemma \ref{upper-b} is proved. $\hfill\Box$

As an immediate consequence of the upper bound of the density, one has
\begin{Lemma}\label{lemma3.9}
  It holds that for any $1<p<\i$,
  \begin{equation}\label{nabla-F-omega}
    \int_0^T\big(\|{\rm div}u\|_\i^3+\|\nabla(F,\o)\|_p^2\big)dt\leq C(M).
  \end{equation}
  \end{Lemma}
{\bf Proof:} First, note that
\begin{equation}\label{div-u-infty}
\int_0^T\|{\rm div}u\|_\i^3dt\leq C\int_0^T
(\|F\|_\i^3+\|P(\r)\|_\i^3)dt\leq C.
\end{equation}
Then for $1<p\leq2$, it follows that
\begin{equation*}
\int_0^T\|\nabla(F,\o)\|_p^2dt\leq C \int_0^T\|\r(H,L)\|_p^2dt\leq C\int_0^T\|\sqrt\r(H,L)\|_2^2\|\sqrt\r\|^2_{\f{2p}{2-p}}dt\leq C.
\end{equation*}
For $\f4{\a}\leq p<\i$,
\begin{equation*}
\begin{array}{ll}
\di \int_0^T\|\nabla(F,\o)\|_p^2dt&\di \leq C\int_0^T\|\r(H,L)\|_p^2dt\\
&\di\leq C\int_0^T\|(H,L)\|_p^2dt \leq C\int_0^T\|(1+|x|^{\f{\a}{2}})\nabla(H,L)\|_2^2dt\leq C.
\end{array}
\end{equation*}
Thus Lemma \ref{lemma3.9} is proved. $\hfill\Box$

\section{Higher order estimates}
\setcounter{equation}{0}

Based on the basic estimates and bound of the density obtained in
Section 3, we can derive some uniform estimates on their higher
order derivatives. We start with estimates on first order
derivatives.

\begin{Lemma}\label{lemma4.1}
  It holds that for any $1< p<+\i$,
  \begin{equation*}
    \sup_{t\in[0,T]}\|(\nabla\r,\nabla P(\r))(t,\cdot)\|_p+\int_0^T\|\nabla u\|_\i^2(t)dt\leq C(M).
  \end{equation*}
  \end{Lemma}
{\bf Proof:} Applying the operator $\nabla$ to the continuity equation
$\eqref{CNS}_1$, one has
\begin{equation}\label{nabla-rho}
(\nabla\r)_t+\nabla
u\cdot\nabla\r+u\cdot\nabla(\nabla\r)+(\nabla\r) {\rm div}
u+\r\nabla({\rm div}u)=0.
\end{equation}
Multiplying the equation \eqref{nabla-rho} by
$p|\nabla\r|^{p-2}\nabla\r$ with $p\geq2$ implies that
\begin{equation*}
(|\nabla\r|^p)_t+{\rm div}(u|\nabla\r|^p)+(p-1)|\nabla\r|^p{\rm
div}u+p|\nabla\r|^{p-2}\nabla\r\cdot(\nabla
u\cdot\nabla\r)+p\r|\nabla\r|^{p-2}\nabla\r\cdot\nabla({\rm div}
u)=0.
\end{equation*}
Integrating the above equation with respect to $x$ over $\mathbb{R}^2$ gives that
\begin{equation}\label{100}
\begin{array}{ll}
\di \f{d}{dt}\|\nabla\r\|_p&\di \leq C\Big[\|\nabla
u\|_\i\|\nabla\r\|_p+\|\nabla{\rm div}u\|_p\Big]\leq C\Big[\|\nabla
u\|_\i\|\nabla\r\|_p+\|\nabla\big(\f{F+P(\r)}{2\mu+\l(\r)}\big)\|_p\Big]\\
&\di \leq C\Big[\Big(\|\nabla
u\|_\i+\|F\|_\i+1\Big)\|\nabla\r\|_p+\|\nabla F\|_p\Big].
\end{array}
\end{equation}
By \eqref{dot-u}, one has
\begin{equation}\label{elliptic-HL}
\mathcal{L}_\r u=\nabla P(\r)+\r \dot u=\nabla P(\r)+\r(L,H)^t.
\end{equation}
Thus the elliptic estimates yields that for any $\f{4}{\a}\leq p<\i,$
\begin{equation}\label{e5}
\begin{array}{ll}
\|\nabla^2u\|_{p}&\di \leq C\big[\|\nabla P(\r)\|_p+\|\r(L,H)\|_p\big]\\
&\di\leq C\big[\|\nabla \r\|_p+\|(L,H)\|_p\big]\leq C\big[\|\nabla \r\|_p+\|(1+|x|^{\f{\a}{2}})\nabla(L,H)\|_2\big].
\end{array}
\end{equation}
By Beal-Kato-Majda type inequality (see \cite{BKM}, \cite{hlx} and
\cite{Kazhikhov}) and \eqref{e5}, it holds that
\begin{equation}\label{103}
\begin{array}{ll}
\di\|\nabla u\|_\i&\di \leq C\big(\|{\rm
div}u\|_\i+\|\o\|_\i\big)\ln(e+\|\nabla^2u\|_{\f{4}{\a}})\\
&\di\leq C\big(\|{\rm
div}u\|_\i+\|\o\|_\i\big)\Big[\ln(e+\|\nabla\r\|_{\f{4}{\a}})+\ln(e+\|(1+|x|^{\f{\a}{2}})\nabla(H,L)\|_2)\Big].
\end{array}
\end{equation}
The combination of \eqref{100} with $p=\f{4}{\a}$ and \eqref{103} yields
that
\begin{equation*}
\begin{array}{ll}
\di\f{d}{dt}\|\nabla\r\|_{\f{4}{\a}}\leq C\|\nabla F\|_{\f{4}{\a}}\\
\di\qquad +C\Big[\big(\|{\rm
div}u\|_\i+\|\o\|_\i\big)\ln(e+\|\nabla(H,L)\|_2)+\|F\|_\i+1\Big]\|\nabla
\r\|_{\f{4}{\a}}\ln(e+\|\nabla\r\|_{\f{4}{\a}}).
\end{array}
\end{equation*}
By the estimates \eqref{F-omega-infty}, \eqref{div-u-infty},
\eqref{nabla-F-omega} and the Gronwall's inequality, it holds that
\begin{equation*}
\sup_{t\in[0,T]}\|\nabla\r\|_{\f{4}{\a}}\leq C,
\end{equation*}
which, together with \eqref{F-omega-infty}, \eqref{div-u-infty}, \eqref{e5} and
\eqref{103},  yields that
\begin{equation}\label{nabla-u-infty}
\int_0^T\|\nabla u\|_\i^2dt\leq C.
\end{equation}
Therefore, by \eqref{nabla-u-infty}, Lemma \ref{lemma-F-o}, Lemma
\ref{lemma3.9} and Gronwall inequality, one can derive from
\eqref{100} that
\begin{equation*}
\sup_{t\in[0,T]}\|\nabla\r\|_p\leq C(\|\nabla\r_0\|_p+1)\leq C,
\qquad\forall p\in(1,+\i).
\end{equation*}
Thus the proof of Lemma \ref{lemma4.1} is completed. $\hfill\Box$
\begin{Lemma}\label{lemma-u-w}
  It holds that
  \begin{equation*}
    \sup_{t\in[0,T]}\Big[\|u(t,\cdot)\|_{\f4\a}+\||x|^{\f\a2}\nabla u(t,\cdot)\|_2\Big]+\int_0^T\||x|^{\f\a2}\sqrt\r\dot u\|_{2}^2(t)dt\leq C(M).
  \end{equation*}
  \end{Lemma}
{\bf Proof:} The momentum equation $\eqref{CNS}_2$ can be rewritten
as
$$
\r\dot{u}+\nabla P(\r)=\mu\Delta u+\nabla((\mu+\l(\r)){\rm div} u).
$$
Multiplying the above equation by $\dot u|x|^\a$ with $\a$ being the
weight in Lemma \ref{lemma-wee} and integrating the resulted
equations with respect to $x$ over $\mathbb{R}^2$ give that
\begin{equation}\label{Ki}
\begin{array}{ll}
\di \f{d}{dt}\int \Big[\mu\f{|\nabla u|^2}{2}+(\mu+\l(\r))\f{({\rm
div} u)^2}{2}-P(\r){\rm div} u\Big]|x|^\a dx+\int \r|\dot u|^2|x|^\a
dx =\int \Big[P(\r)(\dot u-u({\rm div} u))\\[2mm]
\di\qquad+\mu u\f{|\nabla u|^2}{2}-\mu\nabla u\cdot\dot
u+(\mu+\l(\r))u\f{({\rm div} u)^2}{2}-(\mu+\l(\r))({\rm div} u)\dot
u\Big]\cdot \nabla (|x|^\a)dx\\
\di -\int \Big[\mu\sum_{i,j,k=1}^2\partial_{x_i}u_j\partial_{x_i}u_k
\partial_{x_k} u_j+\mu\f{|\nabla u|^2}{2}{\rm div} u-(\mu+\l(\r))({\rm div}
u)\sum_{i,j=1}^2\partial_{x_i} u_j\partial_{x_j}
u_i\\
\di\qquad -\r\l^\prime(\r)\f{({\rm div} u)^3}{2}\Big]|x|^\a dx+\int
\Big[P(\r)\sum_{i,j=1}^2\partial_{x_i} u_j\partial_{x_j}
u_i+(\g-1)P(\r)({\rm div} u)^2\Big]|x|^\a dx:=\sum_{i=1}^3 K_i.
\end{array}
\end{equation}
First, $K_1$ can be estimated as
\begin{equation}\label{K1}
\begin{array}{ll}
\di |K_1|\leq C\int \Big[P(\r)|\dot u|+P(\r)|u||{\rm div}
u|+|u||\nabla u|^2+|\nabla u||\dot u|\\
\di\qquad\qquad\qquad \qquad\qquad\qquad +(1+\l(\r))\big(|u|({\rm
div}
u)^2+|{\rm div} u||\dot u|\big)\Big]|x|^{\a-1}dx\\
\di\qquad \leq C\Big[\|P(\r)|x|^{\f{\a}{2}}\|_2\|\dot
u|x|^{\f{\a}{2}-1}\|_2+\|P(\r)\|_\i\|({\rm div}
u)|x|^{\f{\a}{2}}\|_2\|u|x|^{\f{\a}{2}-1}\|_2\\
\di \qquad\qquad +\big(1+\|\l(\r)\|_\i\big)\big(\|\nabla
u\|_\i\|\nabla u|x|^{\f{\a}{2}}\|_2\|u|x|^{\f{\a}{2}-1}\|_2+\|\nabla
u|x|^{\f{\a}{2}}\|_2\|\dot u|x|^{\f{\a}{2}-1}\|_2\big)\Big]\\
\di \qquad\leq C\Big[\|\nabla\dot u|x|^{\f{\a}{2}}\|_2+\|({\rm div}
u)|x|^{\f{\a}{2}}\|_2\|\nabla u|x|^{\f{\a}{2}}\|_2+\|\nabla
u\|_\i\|\nabla u|x|^{\f{\a}{2}}\|_2^2+\|\nabla
u|x|^{\f{\a}{2}}\|_2\|\nabla\dot u|x|^{\f{\a}{2}}\|_2\Big]\\
\di\qquad\leq C\Big[\big(\|\nabla u\|_\i+1\big)\|\nabla
u|x|^{\f{\a}{2}}\|_2^2+\|\nabla(H,L)|x|^{\f{\a}{2}}\|_2^2+1\Big].
\end{array}
\end{equation}
Then, it follows that
\begin{equation}\label{K2}
\begin{array}{ll}
\di |K_2|\leq C\int(1+\l(\r)) |\nabla u|^3 |x|^\a dx\leq C\|\nabla
u\|_\i\|\nabla u|x|^{\f{\a}{2}}\|_2^2,
\end{array}
\end{equation}
and
\begin{equation}\label{K3}
\begin{array}{ll}
\di |K_3|\leq C\int P(\r)|\nabla u|^2 |x|^\a dx\leq C\|\nabla
u|x|^{\f{\a}{2}}\|_2^2.
\end{array}
\end{equation}
Note that for sufficiently small constant $\s>0$, it holds that
\begin{equation}\label{K0}
\begin{array}{ll}
\di \int \Big[\mu\f{|\nabla u|^2}{2}+(\mu+\l(\r))\f{({\rm div}
u)^2}{2}-P(\r){\rm div} u\Big]|x|^\a dx\\
\di \geq\int \Big[\mu\f{|\nabla u|^2}{2}+(\mu+\l(\r))\f{({\rm div}
u)^2}{2}\Big]|x|^\a dx-\s\int ({\rm div} u)^2|x|^\a dx-C_\s\int
P^2(\r)|x|^\a dx\\
\di \geq\int \Big[\mu\f{|\nabla u|^2}{2}+(\mu+\l(\r))\f{({\rm div}
u)^2}{4}\Big]|x|^\a dx-C,
\end{array}
\end{equation}
if we choose $\s=\f{\mu}{4}$.

Substituting \eqref{K1}, \eqref{K2} and \eqref{K3}  into \eqref{Ki},
and integrating the resulted equation with respect to $t$ over
$[0,t]$, and then using \eqref{K0} and Gronwall inequality, it holds
that
$$
\sup_{t\in[0,T]}\|\nabla
u|x|^{\f{\a}{2}}\|_2^2(t)+\int_0^T\|\sqrt\r|\dot
u||x|^{\f{\a}{2}}\|_2^2(t) dt\leq C,
$$
which, together with the Caffarelli-Kohn-Nirenberg inequality,
completes the proof of Lemma \ref{lemma-u-w}. $\hfill\Box$

\begin{Lemma}\label{lemma-r-h2}
  It holds that for any $2\leq p<+\i$,
  \begin{equation*}
    \sup_{t\in[0,T]}\Big[\|u(t,\cdot)\|_\i+\|\nabla u\|_p+\|(\r_t,P_t)\|_p+\|(\nabla^2\r,\nabla^2 P(\r),\nabla^2u)\|_{2}\Big]+\int_0^T\|\nabla^3u\|_{2}^2dt\leq C.
  \end{equation*}
  \end{Lemma}
{\bf Proof:} By $L^2-$estimates to the elliptic system \eqref{elliptic-HL},
one has
\begin{equation}\label{u2}
\begin{array}{ll}
\di\sup_{t\in[0,T]}\|\nabla^2u\|_{2}&\di \leq C
\sup_{t\in[0,T]}\big(\|\nabla
P(\r)\|_2+\|\r(H,L)\|_2\big)\\
&\di \leq
C\sup_{t\in[0,T]}\big(\|\nabla P(\r)\|_2+\|\sqrt\r(H,L)\|_2\big)\leq
C.
\end{array}
\end{equation}
It follows from the interpolation theorem, Lemma \ref{lemma-u-w} and
\eqref{u2} that
\begin{equation}\label{u-infty}
\sup_{t\in[0,T]}\|u(t,\cdot)\|_\i\leq
C\sup_{t\in[0,T]}\|u(t,\cdot)\|_{\f4\a}^{\f{2}{2+\a}}\|\nabla^2u\|^{\f{\a}{2+\a}}_2\leq
C.
\end{equation}
For any $p\in [2,\i)$, by Sobolev embedding theorem, Lemma
\ref{lemma-u-der} and \eqref{u2}, it holds that
\begin{equation}\label{u-p}
\sup_{t\in[0,T]}\|\nabla u\|_p\leq C\sup_{t\in[0,T]}\|\nabla
u\|_{H^1}\leq C.
\end{equation}
Due to $\eqref{CNS}_1$, one can get $\r_t=-u\cdot\nabla\r-\r~{\rm
div}u$ and $P_t=-u\cdot\nabla P-\r P^\prime(\r)~{\rm div}u$, which,
together with the uniform upper bound of the density  and the
estimates in Lemma \ref{lemma4.1} and \eqref{u-infty}-\eqref{u-p},
yields that
\begin{equation*}\label{rho-t-p}
\sup_{t\in[0,T]}\|(\r_t,P_t)\|_p\leq C, \qquad\forall p\in[2,+\i).
\end{equation*}
Applying $\nabla^2$ to the continuity equation $\eqref{CNS}_1$, then
multiplying the resulted equation by $\nabla^2\r$, and then
integrating with respect to $x$ over $\mathbb{R}^2$, one can get that
\begin{equation}\label{nabla-rho-2}
\begin{array}{ll}
\di\f{d}{dt}\|\nabla^2\r\|_2^2&\di \leq C\Big[\|\nabla u\|_\i\|\nabla^2\r\|^2_2+\|\nabla\r\|_4\|\nabla^2 \r\|_2\|\nabla^2 u\|_4+\|\r\|_\i\|\nabla^2\r\|_2\|\nabla^3u\|_2\Big]\\
&\di\leq C\Big[\big(\|\nabla
u\|_\i+1)\|\nabla^2\r\|^2_2+\|\nabla^3u\|^2_2+1\Big],
\end{array}
\end{equation}
where one has used the fact that
\begin{equation}\label{u24}
\|\nabla^2 u\|_4\leq
C\|\nabla^2 u\|_2^{\f12}\|\nabla^3u\|_2^{\f12}\leq
C\|\nabla^3u\|_2^{\f12}.
\end{equation}
Similarly,
\begin{equation}\label{P2}
\f{d}{dt}\|\nabla^2P(\r)\|_2^2\leq C\Big[\big(\|\nabla
u\|_\i+1)\|\nabla^2P(\r)\|^2_2+\|\nabla^3u\|^2_2+1\Big].
\end{equation}
Note that \eqref{elliptic-HL} implies that
\begin{equation}\label{phi}
\mathcal{L}_\r (\nabla
u)=\nabla^2P(\r)+\nabla[\r(H,L)]+\nabla(\nabla\l(\r){\rm
div}u):=\Phi.
\end{equation}
Then the standard elliptic estimates and the estimate \eqref{u24}
give that
\begin{equation*}
\begin{array}{ll}
\di \|\nabla^3u\|_{2}&\di \leq C\|\Phi\|_2 \leq C\Big[\|\nabla^2P(\r)\|_2+\|\r\|_\i\|\nabla(H,L)\|_2+\|\nabla\r\|_{\f{4}{2-\a}}\|(H,L)\|_{\f4\a}\\
&\di
~~~~~~~\qquad\qquad+\big(\|\nabla^2\r\|_2+\|\nabla\r\|_4^2\big)\|{\rm
div}u\|_\i+\|\nabla\r\|_4\|\nabla^2u\|_4\Big]\\
&\di \leq C\Big[\|(\nabla^2P(\r),\nabla^2\r)\|_2\big(\|\nabla
u\|_\i+1\big)+\|(1+|x|^{\f\a2})\nabla(H,L)\|_2+\|\nabla^3u\|_2^{\f12}\Big].
\end{array}
\end{equation*}
Consequently,
\begin{equation}\label{nabla-u-3}
\di \|\nabla^3u\|_{2}\di \leq
C\Big[\|(\nabla^2P(\r),\nabla^2\r)\|_2\big(\|\nabla
u\|_\i+1\big)+\|(1+|x|^{\f\a2})\nabla(H,L)\|_2+1\Big].
\end{equation}
Substituting \eqref{nabla-u-3} into \eqref{nabla-rho-2} and
\eqref{P2} yields that
\begin{equation*}
\f{d}{dt}\|(\nabla^2\r,\nabla^2P(\r))\|_2^2\leq C\Big[\big(\|\nabla
u\|^2_\i+1\big)\|(\nabla^2\r,\nabla^2P(\r))\|^2_2+\|(1+|x|^{\f\a2})\nabla(H,L)\|^2_2+1\Big].
\end{equation*}
Then the Gronwall's inequality yields that
\begin{equation*}
\begin{array}{ll}
\di\|(\nabla^2\r,\nabla^2P(\r))\|_2^2(t)\di \leq
\Big[\|(\nabla^2\r_0,\nabla^2P_0)\|_2^2\\
\di\qquad\qquad
+C\int_0^T(\|(1+|x|^{\f\a2})\nabla(H,L)\|_2^2+1)dt\Big]e^{\di
C\int_0^T\big(\|\nabla u\|^2_\i+1\big)dt} \leq C,
\end{array}
\end{equation*}
which also implies that
\begin{equation*}
\sup_{t\in[0,T]}\|(\nabla^2\r,\nabla^2P(\r))\|_{2}(t)+
\int_0^T\|\nabla ^3u\|_{2}^2dt\leq C.
\end{equation*}
The proof of Lemma \ref{lemma-r-h2} is completed. $\hfill\Box$

\begin{Lemma}\label{nabla-ut-1}
It holds that for $\f4\a\leq p<\i,$
  \begin{equation*}
\sup_{t\in[0,T]}\|\sqrt\r
u_t(1+|x|^{\f\a2})\|_2^2(t)+\int_0^T\big(\|\nabla
u_t\|_{2}^2+\|u_t\|_p^2\big)dt\leq C(M).
\end{equation*}
\end{Lemma}
{\bf Proof:} First, it holds that
$$
\begin{array}{ll}
\di  \sup_{t\in[0,T]}\|\sqrt\r u_t(1+|x|^{\f\a2})\|_2^2(t) \leq
\sup_{t\in[0,T]}\Big[\|\sqrt\r \dot
u(1+|x|^{\f\a2})\|_2^2(t)+\|\sqrt\r u\cdot\nabla
u(1+|x|^{\f\a2})\|_2^2(t)\Big]\\
\qquad\di \leq \sup_{t\in[0,T]}\Big[\|\sqrt\r
(L,H)^t(1+|x|^{\f\a2})\|_2^2(t)+\|\sqrt\r u\|_\i^2\|\nabla
u(1+|x|^{\f\a2})\|_2^2(t)\Big]\leq C.
\end{array}
$$
Then, one can arrive at
$$
\int_0^T\|\nabla u_t\|_{2}^2dt\leq \int_0^T\Big[\|\nabla\dot
u\|_{2}^2+\|u\|_\i^2\|\nabla^2u\|_2^2+\|\nabla u\|_4^2\Big]dt\leq C.
$$
For any $\f4\a\leq p<+\i,$
$$
\begin{array}{ll}
\di \int_0^T\|u_t\|_p^2dt&\di \leq
\int_0^T(\|(H,L)\|_p^2+\|u\|_\i^2\|\nabla u\|_p^2)dt\\
&\di
\leq\int_0^T(\|(1+|x|^{\f\a2})\nabla(H,L)\|_2^2+\|u\|_\i^2\|\nabla
u\|_p^2)dt\leq C.
\end{array}
$$
Thus the proof of Lemma \ref{nabla-ut-1} is completed. $\hfill\Box$

\begin{Lemma}\label{lemma-rho-tt}
It holds that
  \begin{equation*}
\sup_{t\in[0,T]}\|(\r_t,P(\r)_t,\l(\r)_t)\|_{H^1}(t)+\int_0^T\|(\r_{tt},P(\r)_{tt},\l(\r)_{tt})\|_2^2
dt\leq C.
\end{equation*}
\end{Lemma}
{\bf Proof:} From the continuity equation, it holds that $\r_t=-u\cdot \nabla
\r-\r{\rm div}u$ and
$\r_{tt}=-u_t\cdot\nabla\r-u\cdot\nabla\r_t-\r_t{\rm div}u-\r{\rm
div}u_t,$ and thus
\begin{equation*}
\sup_{t\in[0,T]}\|\nabla\r_t\|_2(t)\leq
\sup_{t\in[0,T]}\Big[\|\nabla\r\|_4\|\nabla
u\|_4+\|u\|_\i\|\nabla^2\r\|_2+\|\r\|_\i\|\nabla^2u\|_2\Big]\leq C.
\end{equation*}
and
\begin{equation*}
\begin{array}{ll}
\di \|\r_{tt}\|_2^2 &\di \leq
\Big[\|u_t\|_{\f4\a}^2\|\nabla\r\|_{\f4{2-\a}}^2+\|u\|_\i^2\|\nabla\r_t\|_2^2+\|\r_t\|_4^2\|\nabla
u\|_4^2+\|\r\|_\i^2\|\nabla u_t\|_2^2\Big]\\
&\di\leq C(\|\nabla u_t\|_{2}^2+\|u_t\|_{\f4\a}^2+1).
\end{array}
\end{equation*}
Therefore, it holds that
\begin{equation*}
 \int_0^T\|\r_{tt}\|_2^2 dt\leq C\int_0^T(\|\nabla u_t\|_{2}^2+\|u_t\|_{\f4\a}^2+1)dt\leq
C.
\end{equation*}
Similarly, one has
\begin{equation*}
\sup_{t\in[0,T]}\|\nabla
(P(\r)_t,\l(\r)_t)\|_2(t)+\int_0^T\|(P(\r)_{tt},\l(\r)_{tt})\|_2^2
dt\leq C.
\end{equation*}
Thus the proof of Lemma \ref{lemma-rho-tt} is completed. $\hfill\Box$

\begin{Lemma}\label{nabla-ut}
  It holds that
  \begin{equation*}
  \begin{array}{ll}
\di \sup_{t\in[0,T]}\Big[t\|(\nabla u_t,\nabla \dot u)\|^2_{2}+\|(\r,P(\r))\|_{W^{2,q}(\mathbb{R}^2)}+\|(\nabla\r,\nabla P(\r))\|_\i\Big]\\
\di\qquad\qquad\qquad\qquad\qquad   +\int_0^Tt\Big[\|\sqrt\r u_{tt}\|_2^2(t)+\|\nabla^2u_t\|_{2}^2(t)\Big] dt\leq C,
 \end{array}\end{equation*}
where $q>2$ is given in Theorem \ref{theorem2}.
\end{Lemma}
{\bf Proof:} The estimates are similar to Lemma 4.5 in \cite{JWX2}
except noting that, by \eqref{dot-u},
$$
\nabla\dot u=\nabla(L,H)^t=\nabla u_t-\nabla (u\cdot\nabla u).
$$
Consequently, it holds that
\begin{equation}
\sup_{t\in[0,T]}\big[t\|\nabla\dot u\|^2_2(t)\big]\leq \sup_{t\in[0,T]}\big[t\|\nabla u_t\|^2_2(t)+t\|u\|_\i^2\|\nabla^2u\|_2^2+t\|\nabla u\|_4^2\big]\leq C.
\end{equation}
We omit the details and  the proof of Lemma \ref{nabla-ut} is
completed. $\hfill\Box$

Based on the estimates obtained so far, similar to Lemma 4.6 in
\cite{JWX2}, one has
\begin{Lemma}\label{lemma4.6}
  It holds that for any $0<\tau\leq T$,
  \begin{equation*}
    \sup_{t\in[0,T]}\big[t^2\|\sqrt\r u_{tt}\|_2^2(t)\big]+\int_0^Tt^2\|\nabla u_{tt}\|_2^2(t) dt\leq C.
  \end{equation*}
\end{Lemma}

Finaaly, we have the following
\begin{Lemma}\label{lemma4.7}
  It holds that
  \begin{equation*}
  \begin{array}{ll}
    \di\sup_{t\in[0,T]}\Big[t\|\nabla \dot u |x|^{\f{\a}{4}}\|^2_{2}+t\|(u_t,\dot u)\|^2_{\f{8}{\a}}+t\|\nabla^3u\|_2^2+t^2\|\nabla^2 u_t\|_2^2+t^2\|\nabla^3u\|_q^2\Big]\\
    \di \qquad\qquad\qquad\qquad\qquad\qquad\qquad  +\int_0^Tt\|\sqrt\r \dot u_{t}|x|^{\f{\a}{4}}\|_2^2(t) dt\leq C.
   \end{array} \end{equation*}
\end{Lemma}
{\bf Proof:} Applying the operator $\partial_t+{\rm div}(u\cdot)$ to
the equation $\eqref{m1}_i~(i=1,2)$ gives that
\begin{equation}
\begin{array}{ll}
\di \r\dot u_{it}+\r u\cdot\nabla\dot u_i=\mu\Delta \dot
u_{it}+\mu{\rm div}(u \Delta u_i)\\
\di \qquad\quad+\partial_{x_it}((\mu+\l(\r)){\rm div} u)+{\rm
div}\big[u\partial_{x_i}((\mu+\l(\r)){\rm div}
u)\big]-\partial_{x_it}P(\r)-{\rm div} (u\partial_{x_i}P(\r)).
\end{array}
\end{equation}
Multiplying the above equation by $\dot u_{it}|x|^{\f\a2}$ with
$\a>0$ to be determined, then summing over $i=1,2$, and then
integrating with respect to $x$ over $\mathbb{R}^2$ imply that
\begin{equation}\label{ut-w}
\begin{array}{ll}
\di \f{d}{dt}\int\big[\f{\mu}{2}|\nabla \dot
u|^2+\f{\mu+\l(\r)}{2}({\rm div}\dot u)^2\big]|x|^{\f\a2}
dx+\int\r|\dot u_t|^2|x|^{\f\a2} dx=-\int \r u\cdot \nabla \dot u\cdot
\dot u_t
|x|^{\f\a2}dx\\
\di \quad-\Big[\int\mu\dot u_{t}\cdot\dot
u_{x_j}\partial_{x_j}(|x|^{\f\a2})dx+\int(\mu+\l(\r))({\rm div}\dot
u)\dot u_t\cdot \nabla(|x|^{\f\a2}) dx\Big]\\
\di \quad+\int \l(\r)_t\f{({\rm div}
\dot u)^2}{2}|x|^{\f\a2}dx+\Big[\int\partial_{x_j}u\cdot\nabla u\cdot(\dot
u|x|^{\f\a2})_{x_jt}dx-\int({\rm div} u)\partial_{x_j}u\cdot(\dot
u|x|^{\f\a2})_{x_jt}dx\\
\di\quad +\int u_{x_j}\cdot\nabla(\dot
u|x|^{\f\a2})_t\cdot u_{x_j} dx-\int(\mu+\l(\r)+\r\l^\prime(\r))({\rm div} u)^2 {\rm
div}(\dot u|x|^{\f\a2})_t dx\\
\di\quad +\int(\mu+\l(\r))({\rm div} u)
u_{x_j}\cdot\nabla(\dot u_j|x|^{\f\a2})_t
dx+\int(\mu+\l(\r))\partial_{x_j} u_k\partial_{x_k} u_j{\rm
div}(\dot u|x|^{\f\a2})_t dx\Big]\\
\di\quad-\Big[\int P(\r)u_{x_j}\cdot\nabla(\dot
u_j|x|^{\f\a2})_t+\int (\g-1)P(\r)({\rm div}u) {\rm div}(\dot
u|x|^{\f\a2})_tdx\Big]:=\sum_{i=1}^5Q_i.
\end{array}
\end{equation}
First, it holds that
\begin{equation}
\label{Q1}
\begin{array}{ll}
\di |Q_1|&\di=|-\int \r u\cdot \nabla \dot u\cdot
\dot u_t
|x|^{\f\a2}dx|\leq \|\sqrt\r \dot u_t|x|^{\f{\a}{4}}\|_2\|\sqrt\r u\|_\i\|\nabla\dot u|x|^{\f{\a}{4}}\|_2\\
&\di\leq C\|\sqrt\r \dot u_t|x|^{\f{\a}{4}}\|_2\|\nabla \dot u(1+|x|^{\f{\a}{2}})\|_2 \leq \s\|\sqrt\r \dot u_t|x|^{\f{\a}{4}}\|^2_2+C_\s\|\nabla \dot u(1+|x|^{\f{\a}{2}})\|_2^2.
\end{array}
\end{equation}
Then, one can obtain
\begin{equation}
\label{Q2}
\begin{array}{ll}
\di Q_2=-\f{d}{dt}\Big[\int\mu\dot u\cdot\dot
u_{x_j}\partial_{x_j}(|x|^{\f\a2})dx+\int(\mu+\l(\r))({\rm div}\dot
u)\dot u\cdot \nabla(|x|^{\f\a2}) dx\Big]\\
\di~~+\int\mu\dot u\cdot\dot
u_{x_jt}\partial_{x_j}(|x|^{\f\a2})dx+\int(\mu+\l(\r))({\rm div}\dot
u)_t\dot u\cdot \nabla(|x|^{\f\a2}) dx+\int\l(\r)_t({\rm div}\dot
u)\dot u\cdot \nabla(|x|^{\f\a2}) dx\\
\di ~~\leq-\f{d}{dt}\Big[\int\mu\dot u\cdot\dot
u_{x_j}\partial_{x_j}(|x|^{\f\a2})dx+\int(\mu+\l(\r))({\rm div}\dot
u)\dot u\cdot \nabla(|x|^{\f\a2}) dx\Big]\\
\di~~\qquad+C\|\dot u|x|^{\f\a2-1}\|_2\|\nabla\dot
u_{t}\|_2+C\|\l(\r)_t\|_\i\|\dot u|x|^{\f\a2-1}\|_2\|\nabla\dot
u\|_2\\
\di ~~\leq-\f{d}{dt}\Big[\int\mu\dot u\cdot\dot
u_{x_j}\partial_{x_j}(|x|^{\f\a2})dx+\int(\mu+\l(\r))({\rm div}\dot
u)\dot u\cdot \nabla(|x|^{\f\a2}) dx\Big]\\
\di~~\qquad+C\|\nabla\dot u|x|^{\f\a2}\|_2\|\nabla(u_{tt}+ u_t\cdot \nabla u+u\cdot\nabla u_t)\|_2\\
\di~~\qquad+C(\|u\cdot\nabla\l(\r)\|_\i+\|\r\l^\prime(\r){\rm div} u\|_\i)\|\nabla\dot u|x|^{\f\a2}\|_2\|\nabla\dot
u\|_2\\
\di ~~\leq-\f{d}{dt}\Big[\int\mu\dot u\cdot\dot
u_{x_j}\partial_{x_j}(|x|^{\f\a2})dx+\int(\mu+\l(\r))({\rm div}\dot
u)\dot u\cdot \nabla(|x|^{\f\a2}) dx\Big]\\
\di~~+C\|\nabla\dot u|x|^{\f\a2}\|_2\Big[\|\nabla u_{tt}\|_2+\|\nabla u\|_\i(\|\nabla u_t\|_2+\|\nabla\dot
u\|_2)+\|\nabla^2 u_t\|_2+\|u_t\|_{\i}\|\nabla^2u\|_2+\|\nabla\dot
u\|_2\Big].
\end{array}
\end{equation}
Obviously,
\begin{equation}\label{Q3}
|Q_3|\leq C\|\l(\r)_t\|_\i\|({\rm div}
\dot u)|x|^{\f\a4}\|_2^2\leq C(1+\|{\rm div} u\|_\i)\|\nabla \dot u|x|^{\f\a4}\|_2^2.
\end{equation}
Now we estimate $Q_4$, which contains six integrals. For simplicity, only the first and the last terms in $Q_4$, denoted by $Q_4^1$ and $Q_4^6$, respectively, will be computed as follows. The others terms in $Q_4$ can be done similarly.
\begin{equation}\label{Q41}
\begin{array}{ll}
\di Q_4^1=\f{d}{dt}\int\partial_{x_j}u\cdot\nabla u\cdot(\dot
u|x|^{\f\a2})_{x_j}dx+\int\partial_{x_j}u_t\cdot\nabla u\cdot(\dot
u|x|^{\f\a2})_{x_j}dx+\int\partial_{x_j}u\cdot\nabla u_t\cdot(\dot
u|x|^{\f\a2})_{x_j}dx\\
\di\quad~ \leq \f{d}{dt}\int\partial_{x_j}u\cdot\nabla u\cdot(\dot
u|x|^{\f\a2})_{x_j}dx+\|\nabla u\|_\i\|\nabla u_t\|_2\Big[\|\nabla \dot u|x|^{\f{\a}{2}}\|_2+\|\dot u|x|^{\f{\a}{2}-1}\|_2\Big]\\
\di\quad~ \leq \f{d}{dt}\int\partial_{x_j}u\cdot\nabla u\cdot(\dot
u|x|^{\f\a2})_{x_j}dx+C\|\nabla u\|_\i\|\nabla u_t\|_2\|\nabla \dot u|x|^{\f{\a}{2}}\|_2\\
\di\quad~ \leq \f{d}{dt}\int\partial_{x_j}u\cdot\nabla u\cdot(\dot
u|x|^{\f\a2})_{x_j}dx+C\Big[\|\nabla \dot u|x|^{\f{\a}{2}}\|^2_2+\|\nabla u\|^2_\i\|\nabla u_t\|^2_2\Big],
\end{array}
\end{equation}
and
\begin{equation}\label{Q46}
\begin{array}{ll}
\di Q_4^6=\f{d}{dt}\int(\mu+\l(\r))\partial_{x_j} u_k\partial_{x_k} u_j{\rm
div}(\dot u|x|^{\f\a2}) dx\\
\di\quad~~~~+\int(\mu+\l(\r))\big(\partial_{x_j} u_k\partial_{x_k} u_j\big)_t{\rm
div}(\dot u|x|^{\f\a2}) dx+\int\l(\r)_t\partial_{x_j} u_k\partial_{x_k} u_j{\rm
div}(\dot u|x|^{\f\a2}) dx\\
\di\quad~ \leq \f{d}{dt}\int(\mu+\l(\r))\partial_{x_j} u_k\partial_{x_k} u_j{\rm
div}(\dot u|x|^{\f\a2}) dx\\
\di\qquad~~+\Big(\|\nabla u\|_\i\|\nabla u_t\|_2+\|\l(\r)_t\|_\i\|\nabla u\|_4^2\Big)\Big[\|\nabla \dot u|x|^{\f{\a}{2}}\|_2+\|\dot u|x|^{\f{\a}{2}-1}\|_2\Big]\\
\di\quad~ \leq \f{d}{dt}\int\partial_{x_j}u\cdot\nabla u\cdot(\dot
u|x|^{\f\a2})_{x_j}dx+C\Big[\|\nabla \dot u|x|^{\f{\a}{2}}\|^2_2+\|\nabla u\|^2_\i(1+\|\nabla u_t\|_2^2)\Big].
\end{array}
\end{equation}
It follows that
\begin{equation}\label{Q5}
\begin{array}{ll}
\di Q_5=-\f{d}{dt}\Big[\int P(\r)u_{x_j}\cdot\nabla(\dot
u_j|x|^{\f\a2})+\int (\g-1)P(\r)({\rm div}u) {\rm div}(\dot
u|x|^{\f\a2})dx\Big]\\
\di\qquad~~+\int P(\r)_tu_{x_j}\cdot\nabla(\dot
u_j|x|^{\f\a2})+\int (\g-1)P(\r)_t({\rm div}u) {\rm div}(\dot
u|x|^{\f\a2})dx\\
\di\qquad~~+\int P(\r)u_{x_jt}\cdot\nabla(\dot
u_j|x|^{\f\a2})+\int (\g-1)P(\r)({\rm div}u)_t {\rm div}(\dot
u|x|^{\f\a2})dx\\
\di\qquad \leq -\f{d}{dt}\Big[\int P(\r)u_{x_j}\cdot\nabla(\dot
u_j|x|^{\f\a2})+\int (\g-1)P(\r)({\rm div}u) {\rm div}(\dot
u|x|^{\f\a2})dx\Big]\\
\di\qquad~~+C\|P(\r)_t\|_\i\|\nabla u\|_2\|\nabla\dot u|x|^{\f\a2}\|_2+C\|P(\r)\|_\i\|\nabla u_t\|_2\|\nabla\dot u|x|^{\f\a2}\|_2\\
\di\qquad \leq -\f{d}{dt}\Big[\int P(\r)u_{x_j}\cdot\nabla(\dot
u_j|x|^{\f\a2})+\int (\g-1)P(\r)({\rm div}u) {\rm div}(\dot
u|x|^{\f\a2})dx\Big]\\
\di\qquad~~+C\Big[\|\nabla\dot u|x|^{\f\a2}\|^2_2+\|\nabla u\|_\i^2+\|\nabla u_t\|_2^2\Big].
\end{array}
\end{equation}
Substituting \eqref{Q1}, \eqref{Q2},  \eqref{Q3},  \eqref{Q41},  \eqref{Q46} and \eqref{Q5} into \eqref{ut-w} and choosing $\s=\f12$ gives that
\begin{equation}\label{dR}
\begin{array}{ll}
  \di \f{d}{dt}R(t)+\f12\| \sqrt\r \dot u_t|x|^{\f\a4}\|_2^2(t)\leq C\Big[(\|\nabla u\|^2_\i+1)(1+\|\nabla u_t\|_2^2+\|\nabla \dot u\|_2^2)+\|\nabla \dot u|x|^{\f\a2}\|_2^2\\
  \di ~~+(1+\|{\rm div} u\|_\i)\|\nabla \dot u|x|^{\f\a4}\|_2^2+\|\nabla \dot u|x|^{\f\a2}\|_2\|\nabla u_{tt}\|_2+\|\nabla^2 u_t\|^2_2+\|u_t\|_{\f4\a}^{\f{4}{\a+2}}\|\nabla^2u_t\|^{\f{2\a}{\a+2}}_{2}\|\nabla^2u\|^2_{2}\Big]\\
  \leq C\Big[(\|\nabla u\|^2_\i+1)(1+\|\nabla u_t\|_2^2+\|\nabla \dot u\|_2^2)+\|\nabla \dot u|x|^{\f\a2}\|_2^2+(1+\|{\rm div} u\|_\i)\|\nabla \dot u|x|^{\f\a4}\|_2^2\\
  \di ~~+\|\nabla \dot u|x|^{\f\a2}\|_2\|\nabla u_{tt}\|_2+\|\nabla^2 u_t\|^2_2+\|u_t\|_{\f4\a}^{2}\|\nabla^2u_t\|^{2}_{2}\Big],
\end{array}
\end{equation}
where
\begin{equation}\label{R(t)}
\begin{array}{ll}
\di R(t)=\f12\Big[\mu\|\nabla \dot u|x|^{\f\a4}\|_2^2+\mu\|({\rm div} \dot u)|x|^{\f\a4}\|_2^2+\|\sqrt{\l(\r)}({\rm div} \dot u)|x|^{\f\a4}\|_2^2\Big]\\
\di +\int\mu\dot u\cdot\dot
u_{x_j}\partial_{x_j}(|x|^{\f\a2})dx+\int(\mu+\l(\r))({\rm div}\dot
u)\dot u\cdot \nabla(|x|^{\f\a2}) dx-\int\partial_{x_j}u\cdot\nabla u\cdot(\dot
u|x|^{\f\a2})_{x_j}dx\\
\di -\int\partial_{x_j}u\cdot\nabla u\cdot(\dot
u|x|^{\f\a2})_{x_j}dx+\int({\rm div} u)\partial_{x_j}u\cdot(\dot
u|x|^{\f\a2})_{x_j}dx -\int u_{x_j}\cdot\nabla(\dot
u|x|^{\f\a2})\cdot u_{x_j} dx\\
\di +\int(\mu+\l(\r)+\r\l^\prime(\r))({\rm div} u)^2 {\rm
div}(\dot u|x|^{\f\a2}) dx-\int(\mu+\l(\r))({\rm div} u)
u_{x_j}\cdot\nabla(\dot u_j|x|^{\f\a2})
dx\\
\di  +\int P(\r)u_{x_j}\cdot\nabla(\dot
u_j|x|^{\f\a2})+\int (\g-1)P(\r)({\rm div}u) {\rm div}(\dot
u|x|^{\f\a2})dx.
\end{array}
\end{equation}
It can be computed that
\begin{equation}\label{I-R}
\int_0^T R(t) dt\leq C\int_0^T\Big[\|\nabla \dot u(1+|x|^{\f\a2})\|_2^2+\|\nabla u\|_4^4+\|\nabla u\|_2^2\Big]dt\leq C,
\end{equation}
and
\begin{equation}
|\int\mu\dot u\cdot\dot
u_{x_j}\partial_{x_j}(|x|^{\f\a2})dx|\leq \f{\mu\a}{2}\|\nabla \dot u|x|^{\f\a4}\|_2\|\dot u|x|^{\f\a4-1}\|_2\leq \f{\mu\a^2}{8}\|\nabla \dot u|x|^{\f\a4}\|_2^2,
\end{equation}
For $r>2$ and close to $2$  and satisfying $\f1s+\f1r=\f12$ and $\f1r+\f{\f\a4-1}{2}=\f\a8\t$ with $\t\in(0,1)$, it holds that
\begin{equation}
\begin{array}{ll}
\di |\int(\mu+\l(\r))({\rm div}\dot
u)\dot u\cdot \nabla(|x|^{\f\a2}) dx|\\
\di \leq \f{\mu\a}{2}\big[\|({\rm div}\dot u)|x|^{\f\a4}\|_2\|\dot u|x|^{\f\a4-1}\|_2+\|\sqrt{\l(\r)}({\rm div}\dot u)|x|^{\f\a4}\|_2\|\l(\r)\|_{s}\|\dot u|x|^{\f\a4-1}\|_{r}\big]\\
\di \leq \f{\mu\a^2}{8}\|({\rm div}\dot u)|x|^{\f\a4}\|_2\|\nabla \dot u|x|^{\f\a4}\|_2+C\|\sqrt{\l(\r)}({\rm div}\dot u)|x|^{\f\a4}\|_2\|\nabla\dot u\|_2^{1-\t}\|\nabla\dot u|x|^{\f\a4}\|_2^{\t}\\
\di \leq \f{\mu\a^2}{8}\|({\rm div}\dot u)|x|^{\f\a4}\|_2\|\nabla \dot u|x|^{\f\a4}\|_2+\s\big[\|\sqrt{\l(\r)}({\rm div}\dot u)|x|^{\f\a4}\|^2_2+\|\nabla\dot u|x|^{\f\a4}\|_2^2\big]+C_\s\|\nabla\dot u\|_2^2.
\end{array}
\end{equation}
By \eqref{phi}, it holds that
\begin{equation}\label{3u}
\begin{array}{ll}
\di \|\nabla^3u\|_{2}&\di \leq C\|\Phi\|_2 \leq C\Big[\|\nabla^2P(\r)\|_2+\|\r\|_\i\|\nabla\dot u\|_2+\|\nabla\r\|_{\f{8}{4-\a}}\|\dot u\|_{\f8\a}\\
&\di
~~~~~~~\qquad\qquad+\big(\|\nabla^2\r\|_4+\|\nabla\r\|_8^2\big)\|{\rm
div}u\|_4+\|\nabla\r\|_\i\|\nabla^2u\|_2\Big]\\
&\di \leq C\Big[1+\|\nabla\dot u\|_2+\|\nabla\dot u|x|^{\f\a4}\|_{2}\Big].
\end{array}
\end{equation}
Then the other terms in \eqref{R(t)} can be estimated as
\begin{equation}
\begin{array}{ll}
\di \Big|-\int\partial_{x_j}u\cdot\nabla u\cdot(\dot
u|x|^{\f\a2})_{x_j}dx -\int\partial_{x_j}u\cdot\nabla u\cdot(\dot
u|x|^{\f\a2})_{x_j}dx\\
\di +\int({\rm div} u)\partial_{x_j}u\cdot(\dot
u|x|^{\f\a2})_{x_j}dx-\int u_{x_j}\cdot\nabla(\dot
u|x|^{\f\a2})\cdot u_{x_j} dx\\
\di+\int(\mu+\l(\r)+\r\l^\prime(\r))({\rm div} u)^2 {\rm
div}(\dot u|x|^{\f\a2}) dx  -\int(\mu+\l(\r))({\rm div} u)
u_{x_j}\cdot\nabla(\dot u_j|x|^{\f\a2})
dx\\
\di +\int P(\r)u_{x_j}\cdot\nabla(\dot
u_j|x|^{\f\a2})+\int (\g-1)P(\r)({\rm div}u) {\rm div}(\dot
u|x|^{\f\a2})dx\Big|\\
 \di \leq C\|\nabla \dot u|x|^{\f\a4}\|_2\big(\|\nabla u|x|^{\f\a4}\|_2+\| u|x|^{\f\a4-1}\|_2\big)\big(\|\nabla u\|_\i+1\big)\\
\di \leq C\|\nabla \dot u|x|^{\f\a4}\|_2\|\nabla u|x|^{\f\a4}\|_2\big(\|\nabla u\|_\i+1\big)\\
\di\leq C\|\nabla \dot u|x|^{\f\a4}\|_2\|\nabla u(1+|x|^{\f\a2})\|_2\big(\|\nabla u\|^{\f12}_{2}\|\nabla^3 u\|_2^{\f12}+1\big) \leq C\|\nabla \dot u|x|^{\f\a4}\|_2\big(\|\nabla^3 u\|_2^{\f12}+1\big)\\
\di \leq C\|\nabla \dot u|x|^{\f\a4}\|_2\big(\|\nabla\dot u\|^{\f12}_2+\|\nabla\dot u|x|^{\f\a4}\|_{2}^{\f12}+1\big)
\leq \s\|\nabla \dot u|x|^{\f\a4}\|_2^2+C_\s\big(\|\nabla\dot u\|^2_2+1\big).
\end{array}
\end{equation}
Thus it holds that
\begin{equation}\label{R(t)2}
\begin{array}{ll}
\di R(t)&\di\geq  \f\mu2\Big[(1-\f{\a^2}{4})\|\nabla \dot u|x|^{\f\a4}\|_2^2+\|({\rm div} \dot u)|x|^{\f\a4}\|_2^2-\f{\a^2}{4}\|({\rm div}\dot u)|x|^{\f\a4}\|_2\|\nabla \dot u|x|^{\f\a4}\|_2\Big]\\
&\di\qquad +(\f12-\s)\|\sqrt{\l(\r)}({\rm div} \dot u)|x|^{\f\a4}\|_2^2-2\s\|\nabla \dot u|x|^{\f\a4}\|_2^2-C_\s\big(\|\nabla\dot u\|^2_2+1\big)\\
\di &\di \geq C^{-1}\Big[\|\nabla \dot u|x|^{\f\a4}\|_2^2+\|({\rm div} \dot u)|x|^{\f\a4}\|_2^2+\|\sqrt{\l(\r)}({\rm div} \dot u)|x|^{\f\a4}\|_2^2\Big]-C\big(\|\nabla\dot u\|^2_2+1\big)
\end{array}
\end{equation}
where in the last inequality one has used the fact that the quadratic term in the square brackets is positively definite when $\a^2<4(\sqrt2-1)$ and one has chosen $\s$ suitably small.

Multiplying the inequality \eqref{dR} by $t$ and then integrating the resulting inequality with respect to $t$ over $[\tau,t_1]$ with both $\tau,t_1\in[0,T]$ give that
\begin{equation}\label{dR1}
\begin{array}{ll}
  \di t_1R(t_1)+\f12\int_{\tau}^{t_1}t\| \sqrt\r \dot u_t|x|^{\f\a4}\|_2^2(t)dt\leq \tau R(\tau)+C\int_{\tau}^{t_1}\Big[t(\|\nabla u\|^2_\i+1)(1+\|\nabla u_t\|_2^2+\|\nabla \dot u\|_2^2)\\
  \di \qquad+\|\nabla \dot u|x|^{\f\a2}\|_2^2+t^2\|\nabla u_{tt}\|_2^2+t\|\nabla^2 u_t\|^2_2+t\|u_t\|_{\f4\a}^2\|\nabla^2u\|^2_{\f{4}{2-\a}}+R(t)\Big]dt\\
  \di\leq \tau R(\tau)+C\int_{\tau}^{t_1}\Big[\|\nabla u\|^2_\i+1+\|\nabla \dot u|x|^{\f\a2}\|_2^2+t^2\|\nabla u_{tt}\|_2^2+t\|\nabla^2 u_t\|^2_2+\|u_t\|_{\f4\a}^2+R(t)\Big]dt.
\end{array}
\end{equation}
By \eqref{I-R}, there exists a subsequence $\tau_k$ such that
\begin{equation*}
\tau_k\rightarrow 0,\qquad \tau_kR(\tau_k)\rightarrow 0,\qquad {\rm as}~~k\rightarrow+\i.
\end{equation*}
Take $\tau=\tau_k$ in \eqref{dR1}, then $k\rightarrow+\i$ and using Gronwall inequality, one can obtain
\begin{equation*}
\sup_{t\in[0,T]}\big[tR(t)\big]+\int_{0}^{T}\| \sqrt\r \dot u_t|x|^{\f\a4}\|_2^2(t)dt\leq C,
\end{equation*}
which, combined with \eqref{R(t)2} and Lemma \ref{nabla-ut}, yields that
\begin{equation*}
\sup_{t\in[0,T]}\big[t\|\nabla \dot u|x|^{\f\a4}\|_{2}^2(t)+t\|(u_t,\dot u)\|^2_{\f8\a}(t)\big]+\int_{0}^{T}\| \sqrt\r \dot u_t|x|^{\f\a4}\|_2^2(t)dt\leq C.
\end{equation*}
Finally, by \eqref{phi}, it holds that
$$
\begin{array}{ll}\|\nabla^3 u\|_2\leq C\big[\|(\nabla^2 \r,\nabla^2P(\r))\|_2(\|\nabla u\|_\i+1)+\|\nabla \dot u\|_2+\|\nabla \r\|_{\f{8}{4-\a}}\|\dot u\|_{\f8\a}+\|\nabla
\r\|_\i\|\nabla^2u \|_2\big]\\
\di \qquad\quad~\leq C\big[\|\nabla u\|_2^{\f12}\|\nabla^3u\|_2^{\f12}+\|\nabla \dot u\|_2+\|\dot u\|_{\f8\a}+1\big]\leq C\big[\|\nabla^3u\|_2^{\f12}+\|\nabla \dot u\|_2+\|\dot u\|_{\f8\a}+1\big],
\end{array}
$$
which implies that
\begin{equation*}
\di \sup_{t\in[0,T]}\big[t\|\nabla^3 u\|_2^2\big]\leq C\sup_{t\in[0,T]}\big[t\|\nabla \dot u\|_2^2+t\|u_t\|_{\f8\a}^2(t)+1\big]\leq C.
\end{equation*}
Then it holds that
\begin{equation*}
\sup_{t\in[0,T]}\big[t^2\|\nabla^2 u_t\|_2^2(t)\big]\leq C\sup_{t\in[0,T]}\Big[t^2\|\sqrt\r u_{tt}\|_2^2(t)+t^2\|u_t\|^2_{\f8\a}(t)+t^2\|\nabla u_t\|_{2}^2+t^2\|\nabla^3 u\|_{2}^2+1\Big]\leq C.
\end{equation*}
and  one can obtain
\begin{equation*}
\sup_{t\in[0,T]}\big[t^2\|\nabla u\|_{W^{2,q}}^2(t)\big]\leq C\sup_{t\in[0,T]}\Big[t^2\|\nabla u\|_{H^2}^2+t^2\|\nabla u_t\|_{H^1}^2+1\Big]\leq C.
\end{equation*}
So the proof of Lemma \ref{lemma4.7} is completed. $\hfill\Box$

\section{The proof of main results}

In this section, we give the proof of our main results.

 {\bf The
proof of Theorem \ref{theorem2}.} We first show that $(\r,u)$ is a
classical solution to \eqref{CNS} if $(\r,u)$ satisfies
\eqref{Jan-16-1}. Since $u\in L^2(0,T;L^{\f4\a}\cap
D^3(\mathbb{R}^2))$ and $u_t\in L^2(0,T;L^{\f4\a}\cap
D^1(\mathbb{R}^2))$, so the Sobolev's embedding theorem implies that
$$
u\in C([0,T];L^{\f4\a}\cap D^2(\mathbb{R}^2))\hookrightarrow C([0,T]\times\mathbb{R}^2).
$$
Then it follows from $(\r,P(\r))\in L^\i(0,T;W^{2,q}(\mathbb{R}^2))$ and $(\r,P(\r))_t\in L^\i(0,T;H^1(\mathbb{R}^2))$ that $(\r,P(\r))\in C([0,T];W^{1,q}(\mathbb{R}^2))\cap C([0,T];W^{2,q}(\mathbb{R}^2)-weak)$. This and  Lemma \ref{nabla-ut} then
imply that
$$
(\r,P(\r))\in C([0,T];W^{2,q}(\mathbb{R}^2)).
$$
Since for any $\tau\in(0,T)$,
$$
(\nabla u,\nabla^2 u)\in L^\i(\tau,T;W^{1,q}(\mathbb{R}^2)\cap
L^2(\mathbb{R}^2)),\qquad (\nabla u_t,\nabla^2 u_t)\in
L^\i(\tau,T;L^2(\mathbb{R}^2)).
$$
Therefore,
$$
(\nabla u,\nabla^2u)\in C([\tau,T]\times\mathbb{R}^2),
$$
Due to the fact that
$$
\nabla (\r,P(\r))\in C([0,T];W^{1,q}(\mathbb{R}^2))\hookrightarrow C([0,T]\times\mathbb{R}^2)
$$
and the continuity equation $\eqref{CNS}_1$, it holds that
$$
\r_t=u\cdot\nabla \r+\r{\rm div}u\in C([\tau,T]\times\mathbb{R}^2).
$$
It follows from the momentum equation $\eqref{CNS}_2$ that
$$
\begin{array}{ll}
  (\r u)_t&\di =\mathcal{L}_{\r}u-{\rm div}(\r u\otimes u)-\nabla P(\r)\\
  &\di =\mu \Delta u+(\mu+\l(\r))\nabla({\rm div} u)+({\rm div}u)\nabla\l(\r)+\r u\cdot\nabla u+\r u{\rm div} u+(u\cdot\nabla\r) u-\nabla P(\r)\\
  &\di \in C([\tau,T]\times\mathbb{R}^2).
\end{array}
$$

Then the proof of Theorem \ref{theorem2} follows from Lemma
\ref{lemma0} which is about the local well-posedness of  the
classical solution and the global (in time) a priori estimates in
Sections 3-4. In fact, by Lemma \ref{lemma0}, there exists a local
classical solution $(\r,u)$ on the time interval $(0,T_*]$ with
$T_*>0$.   Now let $T^*$ be the maximal existing time of the
classical solution $(\r,u)$ in Lemma \ref{lemma0}. Then obviously
one has  $T^*\geq T_*$. Now we claim that $T^*> T$ with $T>0$ being
any fixed positive constant given in Theorem \ref{theorem2}.
Otherwise, if $T^*\le T$, then all the a priori estimates in
Sections 3-4 hold with $T$ being replaced by $T^*$. In particular,
from the inequality \eqref{g1}, it holds that
$$(1+|x|^{\f\a2})\sqrt\r\dot u\in C([0,T^*];L^2(\mathbb{R}^2)).$$
Therefore, it follows from a priori estimates in Sections 3-4 that
$(\r,u)(x,T^*)$ satisfy \eqref{in-d1} and the compatibility
condition \eqref{cc} at time $t=T^*$ with $g(x)=\sqrt\r \dot
u(x,T^*)$. By using Lemma \ref{lemma0} again, there exists a
$T_1^*>0$ such that the classical solution $(\r,u)$ in Lemma
\ref{lemma0} exist on $(0,T^*+T_1^*]$, which contradicts with $T^*$
being the maximal existing time of the classical solution $(\r,u)$.
Thus it holds that $T^*> T$, and the proof of Theorem \ref{theorem2}
is completed. $\hfill\Box$

\vspace{3mm}

 {\bf The proof of Theorem \ref{theorem}.}
 Based on Theorem  \ref{theorem2}, one can
prove Theorem \ref{theorem} easily as follows. Since
$$
\r_0\in H^3(\mathbb{R}^2)\hookrightarrow W^{2,q}(\mathbb{R}^2)
$$
for any $2<q<+\i$, it follows that under the conditions of Theorem \ref{theorem}, Theorem \ref{theorem2} holds for any $2<q<+\i$.
Thus, we need only to prove the higher order regularity presented  in Theorem \ref{theorem}. $\hfill\Box$

\begin{Lemma}\label{new}
  It holds that
$$
  \begin{array}{ll}
    \di\sup_{t\in[0,T]}\Big[\|\sqrt\r \nabla^3u\|_{2}(t)+\|(\r,P(\r),\l(\r))\|_{H^3}(t)\Big]+\int_0^T\|\nabla^4u\|_{2}^2 dt\leq C.
   \end{array} $$
\end{Lemma}

The proof of Lemma \ref{new} is completely similar to Lemma 6.1 in
\cite{JWX2}. We
omit the details here for simplicity and the proof Theorem \ref{theorem} is
complete. $\hfill\Box$

\section*{Acknowledgments}
The results of this paper were presented in International Workshop on
Partial Differential Equations and Mathematical Physics in Jiangsu
University, China, from May 25 to May 28, 2012, The 14th International on
Hyperbolic Problem: Theory, Numerics, Applications in University of
Padova, Italy, from June 25 to June 29, 2012, and the conference "Dynamics of Nonlinear Dispersive and
Fluid Mechanics Equations" in Peking University, China, from June 25 to July
13, 2012. The authors would like to thank the organizers of these
conferences for their kind invitations and giving chance to present
our results.

\end{document}